\documentclass[11pt]{article}

\usepackage[margin=1.05in]{geometry}
\usepackage{amsmath,amssymb,amsthm,mathtools,bbm,mathrsfs}
\usepackage{enumitem}
\usepackage[hidelinks]{hyperref}
\usepackage{microtype}
\usepackage{dsfont}
\usepackage{color}
\usepackage{tikz,tikz-cd}
\usepackage{subcaption}
\usepackage{needspace}
\usepackage{comment}

\allowdisplaybreaks
\numberwithin{equation}{section}

\newtheorem{innercustomthm}{Theorem}
\newenvironment{customthm}[1]
  {\renewcommand\theinnercustomthm{#1}\innercustomthm}
  {\endinnercustomthm}

\newtheorem{theorem}{Theorem}[section]
\newtheorem{proposition}[theorem]{Proposition}
\newtheorem{lemma}[theorem]{Lemma}
\newtheorem{corollary}[theorem]{Corollary}
\theoremstyle{definition}
\newtheorem{definition}[theorem]{Definition}
\newtheorem{remark}[theorem]{Remark}

\newcommand{\R}{\mathbb R}
\newcommand{\C}{\mathbb C}
\newcommand{\E}{\mathbb E}

\newcommand{\Z}{\mathbb Z}
\newcommand{\PP}{\mathbb P}

\newcommand{\don}{\mathds 1}
\newcommand{\N}{\mathbb N}

\newcommand{\U}{\mathbb U}

\newcommand{\GFF}{\mathrm{GFF}}

\renewcommand{\d}{\mathrm d}
\newcommand{\dd}{\mathrm d}

\newcommand{\ii}{\mathbf i}

\newcommand{\eps}{\varepsilon}

\newcommand{\cA}{\mathcal A}
\newcommand{\AB}{\mathcal{A}^\beta}

\newcommand{\sK}{\mathscr K}
\newcommand{\bH}{\mathbf H}
\newcommand{\bI}{\mathbf I}
\newcommand{\bF}{\mathbf F}
\newcommand{\bb}{\boldsymbol{\upsilon}}
\newcommand{\bu}{\mathbf u}
\newcommand{\bp}{\mathbf p}
\newcommand{\bL}{\mathbf L}
\newcommand{\bLC}{\mathbf {L}^{\mathbf c}}

\newcommand{\bdel}{\boldsymbol{\delta}}
\newcommand{\bDel}{\mathbf\Delta}
\newcommand{\ttt}{\boldsymbol{\tau}}

\newcommand{\bk}{\mathbf k}
\newcommand{\bQ}{\mathbf Q}
\newcommand{\fB}{\mathfrak B}

\newcommand{\llbracket}{\mathopen{[\mkern-3mu[}}
\newcommand{\rrbracket}{\mathclose{]\mkern-3mu]}}

\newcommand{\Cov}{\operatorname{Cov}}

\title{Gaussian free field limit for the Airy$_\beta$ line ensemble}
\author{Vadim Gorin \and Lingfu Zhang}
\date{\today}

\begin{document}
\maketitle

\begin{abstract}
We prove that the centered height function of the Airy$_\beta$ line ensemble converges for each $\beta>0$, under rescaling of space and time variables, towards the Gaussian free field in the half-plane. The conformal structure for GFF involves the spatial change of coordinate $x\mapsto \sqrt{-x}$. The $\beta$--dependence is isolated in the field rescaling by $\sqrt{\beta}$. For the proofs we study the asymptotic expansions of the moments of Laplace transforms of the Airy$_\beta$ line ensemble developed in our previous work with Xu.
\end{abstract}

\setcounter{tocdepth}{2}
\tableofcontents

\section{Introduction}
\label{sec:result}

\subsection{Overview}

The Airy$_\beta$ line ensemble and the two-dimensional Gaussian free field are two
universal objects which appear in the asymptotic study of random matrices
and random surfaces, respectively.

For arbitrary $\beta>0$, the Airy$_\beta$ line ensemble 
is a stationary random collection of ordered continuous curves
\[
 \cA_1^\beta(\ttt)\ge\cA_2^\beta(\ttt)\ge\cdots,
 \qquad \ttt\in\R.
\]
At a fixed time $\ttt$, its particles
are distributed as the Airy$_\beta$ point process, which is a universal scaling limit for the largest eigenvalues
in models of self-adjoint random matrices \cite{BEY,BFG,KRV} (here the classical values $\beta=1,2,4$ correspond to real, complex, and quaternionic matrices), and describes the spectrum of the stochastic Airy
operator \cite{ES,RRV}. More generally, the random curves $\{\cA_i^\beta(\ttt)\}_{i=1,2,\dots; \ttt\in\R}$ are constructed in \cite{GXZ} at each $\beta>0$ as the common edge scaling limit of the Dyson Brownian motion and
of the Gaussian $\beta$ corners process, see also \cite{Landon,HZ,XuJack} for related results. At $\beta=2$, the Airy$_2$ line ensemble is a universal object for edge fluctuations in random matrix theory, 2d statistical mechanics, and for growth models in the KPZ universality class, see \cite{PS,CorwinKPZ,CH,AH,AHCharacterization,GFB} for background. Universal properties are also expected, but less understood at other values of $\beta$, see, e.g.\ \cite{Sodin} for a discussion of the $\beta=1$ case. The parameter $\beta$ controls the strength of repulsion between particles. Its effect on the trajectories is visible in Figures \ref{Figure_Airy_surface_1}, \ref{Figure_Airy_surface_2}, \ref{Figure_Airy_surface_8}. 

On the other hand, the two-dimensional Gaussian free field (GFF) with zero Dirichlet boundary conditions is an
analogue of the Brownian bridge. It is defined on planar domains admitting a Dirichlet Green
function and is conformally invariant. We refer to \cite{Sheffield,WernerPowell2022GFF} for a general introduction to the Gaussian free field. This object is a universal scaling limit for a wide class of models of random surfaces. An informal definition treats GFF as a random function $f$ with probability density proportional to $\exp\bigg(-\frac{1}{2} \iint ||\nabla f||^{2}\,dx\,dy\bigg)$, which predicts its asymptotic appearance in the models of random height functions, in which probability weights can be approximated by an exponentiated integral of a function of the gradient of the height. In many situations this prediction has been confirmed by rigorous results, cf.\ discussions in \cite{BPR,GorinBook,BorotGG}.

The main goal of this paper is to establish a direct limit transition between the two universal objects:
we prove that the centered height function of the Airy$_\beta$ line ensemble, as space
and time scales grow, converges towards the Gaussian free field for every $\beta>0$, as detailed in Theorem \ref{thm:surface} and Theorem \ref{Theorem_main_detailed}. In particular, the Airy line ensemble, for which the spatial and temporal directions look very different, has hidden conformal invariance on large scales; though, in order to see it, one has to change coordinates by two-dimensional map $(\ttt,x)\mapsto \frac{\ttt}{2}+\ii\sqrt{-x}$. The dependence on $\beta$ also becomes especially simple in the limit, as it turns into the deterministic prefactor $\sqrt{\frac{2}{\beta\pi}}$ in front of the GFF.

The technical difficulty in developing the convergence result is to find an effective description of the Airy$_\beta$ line ensemble. Here we rely on the recent constructions in \cite{GXZ}: that article presents explicit, but intricate formulas for the expectations of the products of Laplace transforms of the particle positions $\{\cA_i^\beta(\ttt)\}_{i=1}^\infty$. We analyze these formulas and find asymptotic cancellations and simplifications, which eventually lead to the Gaussian limits and the GFF covariance structure.

\subsection{The Airy surface and the main result}

\begin{figure}[t]
 \centering
 \includegraphics[width=0.99\linewidth]{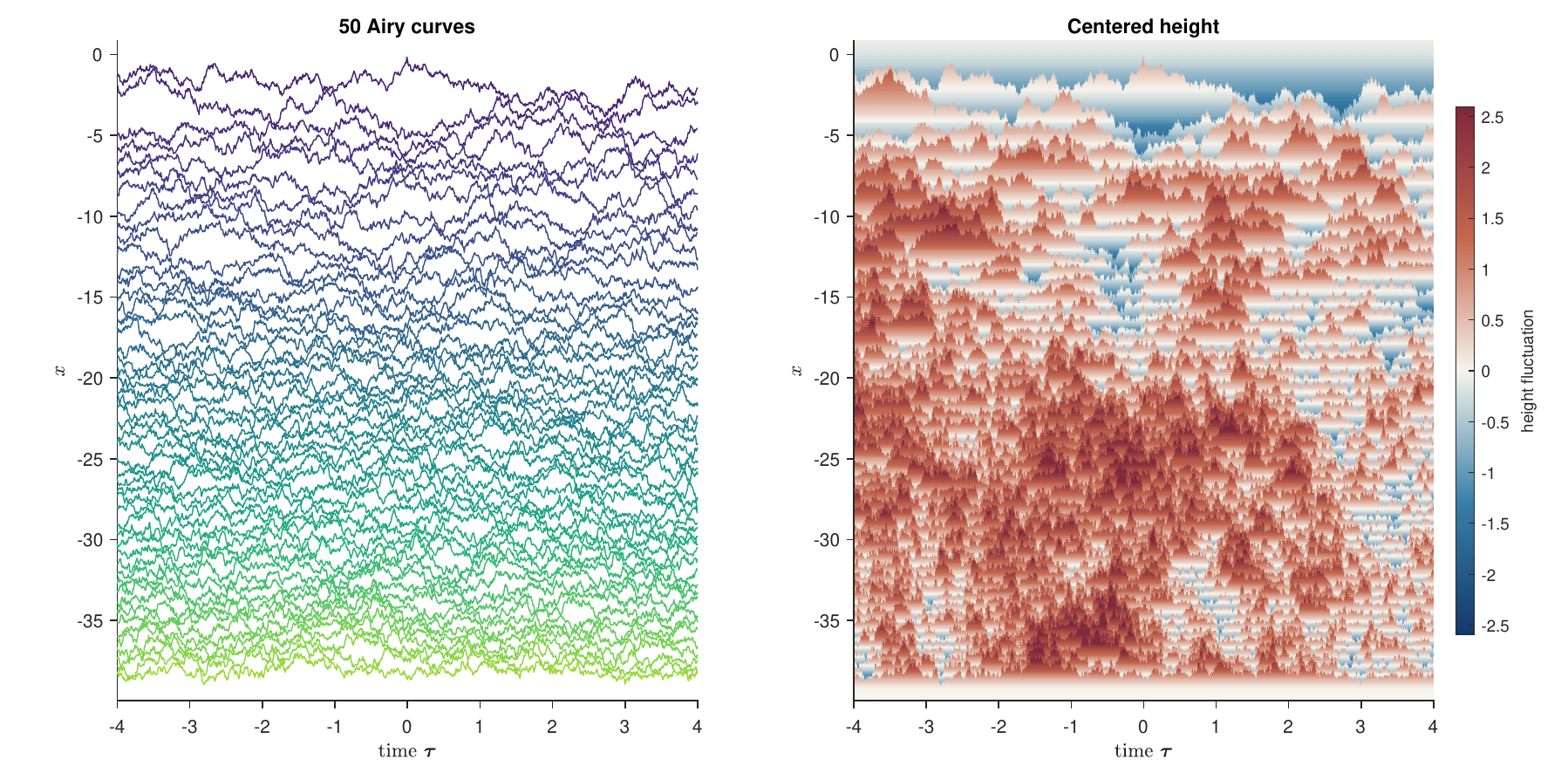}
 \caption{Airy$_1$ line ensemble and corresponding Airy surface
 fluctuation field, approximated by Dyson Brownian motion with
 $N=1000$ particles started from equilibrium. The left panel shows the
 top $50$ curves; the right panel shows the centered height function. The simulations in Figures \ref{Figure_Airy_surface_1}, \ref{Figure_Airy_surface_2}, \ref{Figure_Airy_surface_8} use the same random seed.}
 \label{Figure_Airy_surface_1}
\end{figure}

\begin{figure}[t!]
 \centering
 \includegraphics[width=0.99\linewidth]{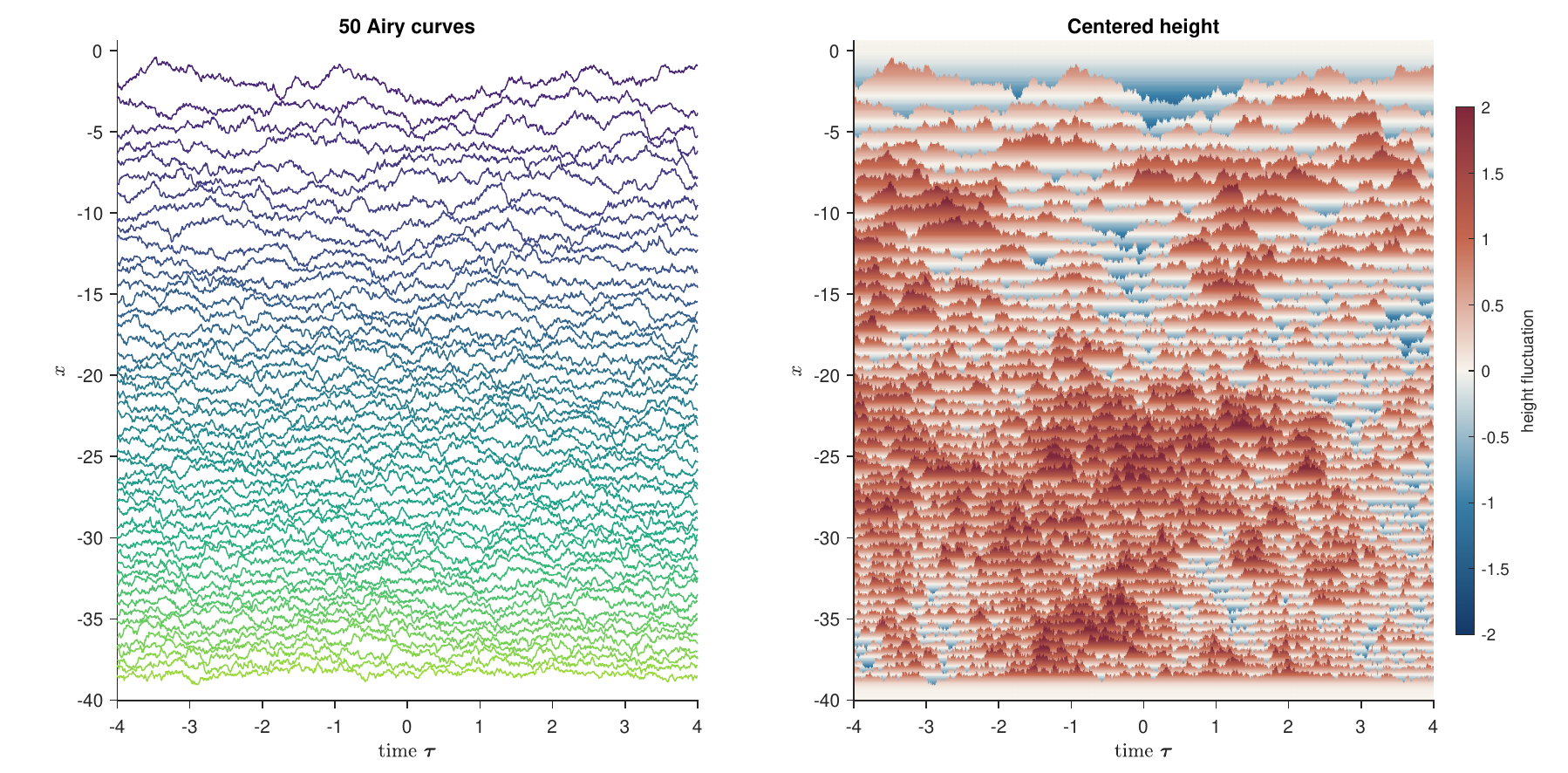}
 \caption{Airy$_2$ line ensemble and corresponding centered Airy surface
 fluctuation field.}
 \label{Figure_Airy_surface_2}
\end{figure}

\begin{figure}[t!]
 \centering
 \includegraphics[width=0.99\linewidth]{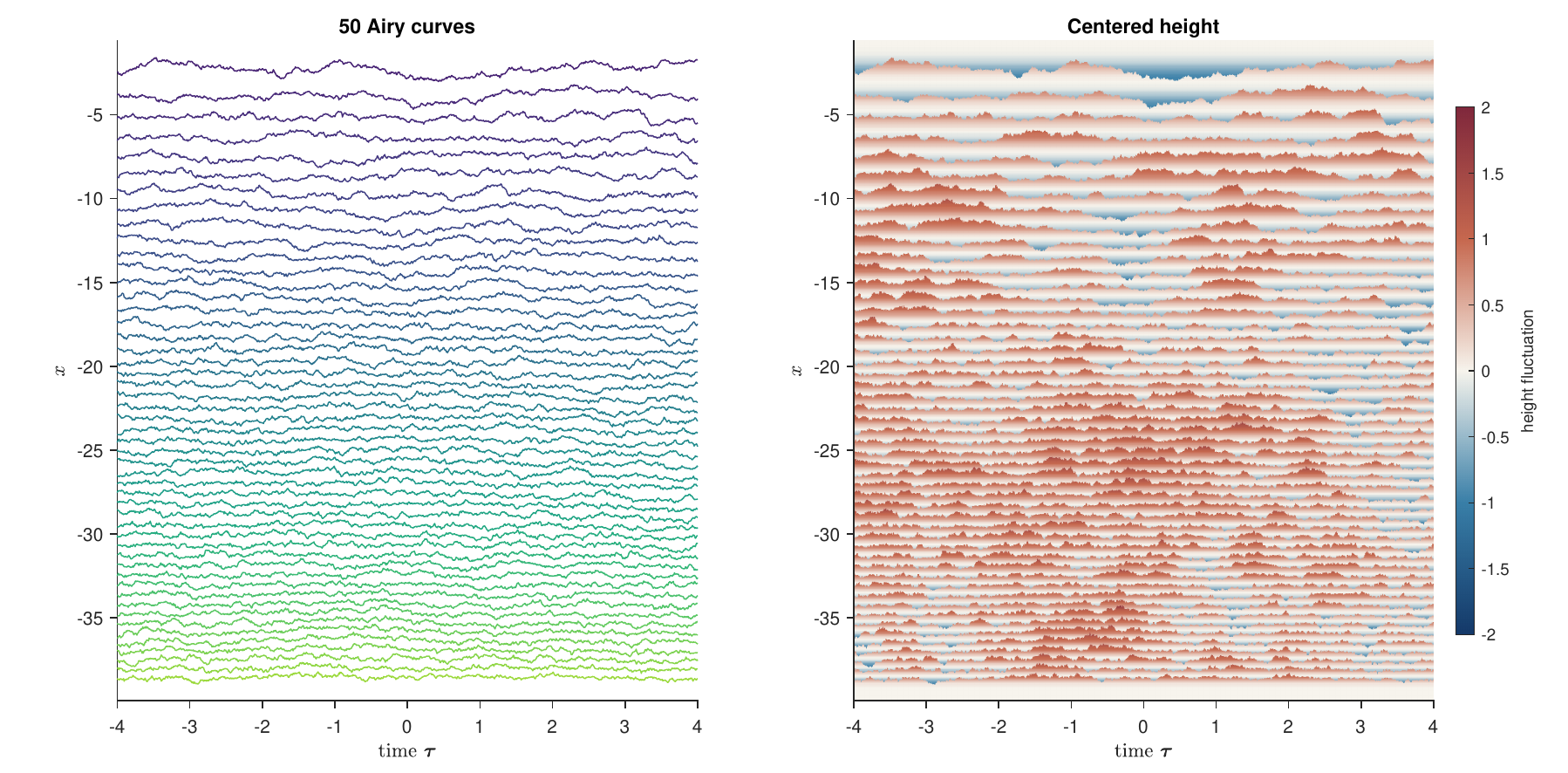}
 \caption{Airy$_8$ line ensemble and corresponding centered Airy surface
 fluctuation field. }
 \label{Figure_Airy_surface_8}
\end{figure}

We view the curves of the Airy$_\beta$ line ensemble as level lines of a
random surface. The \emph{Airy$_\beta$ surface}, or \emph{height function} of the Airy$_\beta$ line ensemble,
is defined by
\begin{equation}\label{eq:intro-airy-surface}
 \cA^\beta(\ttt,x)
 =\#\{i\in\N:\cA_i^\beta(\ttt)>x\},
 \qquad (\ttt,x)\in\R^2.
\end{equation}
For $\eta>0$ we introduce the centered height function in rescaled coordinates 
\begin{equation}\label{eq:intro-centered-surface}
 H_\eta(\ttt,x)
 =\cA^\beta(\eta^{-1/2}\ttt,\eta^{-1}x)
 -\E[\cA^\beta(\eta^{-1/2}\ttt,\eta^{-1}x)],  \qquad (\ttt,x)\in\R^2.
\end{equation}
Figures \ref{Figure_Airy_surface_1}--\ref{Figure_Airy_surface_8}
illustrate the line ensemble and the corresponding centered height
function for $\beta=1,2,8$ and $\eta=1$. We are interested in the asymptotics as $\eta\to0+$, so that the spatial scale grows at speed $\eta^{-1}$
and the time scale grows at speed $\eta^{-1/2}$. For each fixed $x<0$,
we expect the one-point asymptotics\footnote{We only prove the leading
term of the expectation formula \eqref{eq:EAconst} in
Corollary~\ref{cor:height-mean}. We expect that full proofs of these
asymptotics can be obtained using the stochastic Airy operator of
\cite{ES,RRV}. For $\beta=1,2$, these formulas are known; see
\cite{CCAiry,BGS}.}
\begin{equation}\label{eq:EAconst}
 \E[\cA^\beta(\eta^{-1/2}\ttt,\eta^{-1}x)]
 =\frac{2}{3\pi}\eta^{-3/2}(-x)^{3/2}
  +\frac{1}{2\beta}-\frac14+o(1),
\end{equation}
\begin{equation}\label{eq:intro-one-point-variance}
 \operatorname{Var}[\cA^\beta(\eta^{-1/2}\ttt,\eta^{-1}x)]= \operatorname{Var}[H_\eta(\ttt,x)]
 =\frac{3}{2\beta\pi^2}\log(\eta^{-1})
  +o\bigl(\log(\eta^{-1})\bigr).
\end{equation}
Note that the variance in \eqref{eq:intro-one-point-variance} diverges as $\eta\to 0+$, reflecting roughness of the height fluctuations. Instead of looking at one-point fluctuations, which have
to be rescaled, our asymptotic theorem involves integrals of the height function against sufficiently
smooth test-measures. Such integrals do not need any rescaling for distributional limits as $\eta\to 0+$ and give rise to the Gaussian free field, which we now define.

\medskip

Let $\U=\{z\in\C:\operatorname{Im}z>0\}$ be the upper half-plane. The Gaussian free field on $\U$ with Dirichlet boundary conditions, $\GFF_\U(z)$, $z\in\U$, is a centered generalized Gaussian field  with covariance
kernel
\begin{equation}\label{eq:intro-green}
 G_\U(z,w)
 =-\frac{1}{2\pi}\log\left|\frac{z-w}{z-\overline w}\right|, \qquad z,w\in\U.
\end{equation}
In parallel to \eqref{eq:intro-one-point-variance} growing as $\eta\to 0+$, the values of $\GFF_\U$ at individual points are not defined, while integrals against suitable test-measures are bona-fide mean $0$ Gaussian random variables with
covariances obtained by integrating \eqref{eq:intro-green}, see Section \ref{sec:pullback-proof} for further details.

The fluctuations of the Airy surface are governed not by $\GFF_\U(z)$ per se, but rather its version in a modified coordinate system.
The relevant change of coordinates is
\begin{equation}\label{eq:intro-omega}
 \Omega:\R\times\R_{<0}\longrightarrow\U,
 \qquad \Omega(\ttt,x)=\frac{\ttt}{2}+\ii\sqrt{-x}.
\end{equation}
We denote by $\overline{\GFF_\U\circ\Omega}$ the pullback of $\GFF_\U$ with respect to $\Omega$,
extended by zero to $x\ge0$. Its covariance on $x,y<0$ is
$G_\U(\Omega(\ttt,x),\Omega(\ttt',y))$, and is zero if either spatial
coordinate is nonnegative.

\begin{theorem}\label{thm:surface}
For each fixed $\beta>0$, we have
\begin{equation}\label{eq:intro-main-limit}
 \lim_{\eta\to 0+} \sqrt{\frac{\beta\pi}{2}}\,H_\eta
 \stackrel{d}= \overline{\GFF_\U\circ\Omega},
\end{equation}
where the convergence holds in joint moments (and in joint distribution) for integrals against all sufficiently smooth signed test-measures.
\end{theorem}
We refer to Theorem \ref{Theorem_main_detailed} for the description of the precise class of signed test-measures which we allow: they should have continuous density in $x$ direction, exponential decay as $x\to-\infty$, and integrability in $\ttt$ direction (no continuity in $\ttt$ direction is required and the measures are allowed to be supported on fixed $\ttt$ slices). For instance, we prove that for any continuous compactly supported $f(\ttt,x)$, the distributional limit
$$
\lim_{\eta\to 0+}
\sqrt{\frac{\beta\pi}{2}}
       \int_{\R^2} f(\ttt,x) H_\eta(\ttt,x)\, \dd \ttt \dd x
$$
is a mean $0$ Gaussian random variables with variance
$$
 \int_{\R\times \R_{<0}\times\R\times\R_{<0}}
   f(\ttt_1,x_1) f(\ttt_2,x_2)
   G_\U\bigl(\Omega(\ttt_1,x_1),\Omega(\ttt_2,x_2)\bigr)\, \dd \ttt_1 \dd x_1\,\dd \ttt_2 \dd x_2.
$$

\smallskip

For $\beta=2$, \eqref{eq:intro-main-limit} is also\footnote{The projects were carried out independently and in parallel.} established in \cite{DFZ} for a smaller class of signed test-measures using the determinantal point process structure of the Airy$_2$ line ensemble. The determinantal structures do not extend to other values of $\beta$ and we follow a different path.

\subsection{Methodology}

We first explain why one should expect Theorem \ref{thm:surface} to hold. For random domino and
lozenge tilings, a formal quadratic expansion of the variational principle of \cite{CKP} predicts GFF
fluctuations of the height function; see \cite[Lecture~12]{GorinBook}.  Although the general prediction
remains open, it has been proved for broad classes of domains, see, e.g.\ \cite{Kenyon,KenyonHoneycomb,Petrov,
BufetovGorin,HuangHeight,BorotGG}. On the other hand, the Airy$_2$ line ensemble is known to govern the edge limits for the liquid regions in such models, see \cite{JohanssonAztec,PetrovEdge,AH}. Hence, the $\beta=2$ version of Theorem \ref{thm:surface} can be predicted through a limit transition from the tilings$\to$GFF link.

For other values of $\beta$, the Airy$_\beta$ line ensemble is an object of random matrix theory, but the connection survives in several ways. It was first noticed in \cite{Spohn1998} that the Gaussian free field governs the asymptotic fluctuations of the Dyson Brownian motion (DBM) on the circle for all values of $\beta>0$. While the circular geometry does not allow an edge limit to the Airy$_\beta$ line ensemble, the DBM on the line is a close relative, which has a direct scaling limit to Airy$_\beta$, as shown in \cite{Landon,GXZ,HZ}. A complementary connection comes from matrix corners. \cite{Borodin,BorodinDynamic} realized that since joint eigenvalue distributions for corners of GUE random matrices can be obtained as a limit of random tilings, the GFF connection for tilings should lead to GFF asymptotics for corners, which can then be checked directly for real and complex Wigner matrices, corresponding to $\beta=1,2$. Similar phenomenology was subsequently extended to general $\beta$ corners processes in \cite{BG}. From here, \cite{GXZ} explains that the Airy$_\beta$ line ensemble is an edge scaling limit of the corners processes, hence predicting the connection to the Gaussian free field once again.

\bigskip

For the actual proof, we do not directly use any of the connections of the previous two paragraphs: the edge limit of random matrix models is quite delicate and rigorously transferring the GFF limit through the edge scaling limit is difficult. Instead, we rely on the construction of the Airy$_\beta$ line ensemble from \cite{GXZ}. The approach of that paper suggests working with centered
Laplace transforms, which are related to integrals of the height function by integration by parts:
\begin{equation}\label{eq:intro-laplace}
 \begin{split}
 X_\eta(\ttt,\bk)
 &=\sum_{i=1}^\infty
   \exp\left(\frac{\eta\bk}{2}
             \cA_i^\beta(\eta^{-1/2}\ttt)\right)
   -\E\left[\sum_{i=1}^\infty
   \exp\left(\frac{\eta\bk}{2}
             \cA_i^\beta(\eta^{-1/2}\ttt)\right)\right] \\
 &=\frac{\bk}{2}\int_\R e^{\bk x/2}H_\eta(\ttt,x)\dd x,
 \qquad \bk>0.
 \end{split}
\end{equation}
The formulas of \cite{GXZ}, recalled in
Section \ref{Section_moments_from_GXZ}, express joint moments of the
uncentered Laplace transforms through integrals over objects called blocks. Our main technical task
is to analyze these integrals as $\eta\to 0$. We pass from moments to cumulants and associate a graph
with each block configuration. A cancellation in Section \ref{sec:basic} expresses the cumulants in
terms of connected graphs. The naive asymptotic expansions of integrals for cumulants still contain divergent terms, which
are removed by further cancellations in Section \ref{sec:threecan}. Eventually, cumulants of order three and higher for \eqref{eq:intro-laplace} tend to $0$, while the limit of the second cumulant (i.e., covariance) is computed in Section \ref{sec:covariance} to be a fourfold integral,
three of whose integrations can be evaluated explicitly, leading to the exact limit expression:
\begin{equation}\label{eq:intro-laplace-covariance}
 \lim_{\eta\to0+}
 \Cov\bigl(X_\eta(\ttt_1,\bk_1),X_\eta(\ttt_2,\bk_2)\bigr)
 =\frac{1}{\beta\pi\sqrt{\bk_1\bk_2}}
 \int_0^\infty g\exp\left(
 -\frac{|\ttt_1-\ttt_2|g}{2}
 -\frac{g^2}{2\bk_1}-\frac{g^2}{2\bk_2}\right)\dd g.
\end{equation}
Using the Fourier representation of the logarithmic kernel \eqref{eq:intro-green}, we convert in Lemma \ref{lem:pullback-energy} the right-hand side of \eqref{eq:intro-laplace-covariance} into covariance of pairings of the GFF against exponential test-functions on fixed $\ttt$ slices.

The remaining work in Section \ref{sec:pullback-proof} is to extend convergence from exponential test-functions to general admissible test-measures. This is done using moment comparisons of independent interest, which we develop in Lemma \ref{lem:surface-association}, as a corollary from the positive association of the
Airy$_\beta$ particles at fixed time proved in \cite{BKM}.

\medskip
\noindent\textbf{Acknowledgements.} The work of V.G.\ was partially supported by NSF grant DMS-2553163. The work of L.Z. was partially supported by NSF grant DMS-2505625 and a Sloan fellowship.
AI tools were used when preparing this manuscript.

\section{Moments of the Laplace transforms}
\label{Section_moments_from_GXZ}

In this section we recall the main construction of \cite{GXZ}. The definitions are pulled directly from \cite[Section 2]{GXZ} and we refer to that article for the much more detailed discussions; we also synchronized the numbering of definitions with that paper.

\begin{customthm}{GXZ}[Theorem 1.4 in \cite{GXZ}]\label{Theorem_moments}
For any $m\in \N$, $\vec \bk \in \R^m_+ $, and $\vec \ttt\in \R^m$ such that $\ttt_1\le \cdots \le \ttt_m$, the joint $m$th moment for the Laplace transforms of $\{\AB_i\}_{i\ge 1}$ is
\begin{multline} \label{eq_expected_Laplace}
\E\left[ \prod_{\ell=1}^m \left(\sum_{i=1}^\infty
\exp(\bk_\ell \AB_i(\ttt_\ell)/2) \right) \right]=\bL_\beta(\vec\bk, \vec\ttt)\\
= \mathrm{p.v.}\int_{\sK[\vec\bk]} \exp\left(\sum_{\ell=1}^{m-1}
(\ttt_\ell-\ttt_{\ell+1}) \bH(\bQ_\ell)/2 \right)\bI_{\beta,\vec\bk}[\vec \bp, \vec\bb, \bH] \dd(\vec\bp,  \vec\bb,\bH),
\end{multline}
where $\bQ_\ell = \sum_{\ell'=1}^\ell \bk_{\ell'}$ and the principal value integral is taken in the following sense.
For any small enough $\epsilon>0$, we let  $\sK_\epsilon=\sK_\epsilon[\vec\bk]$ denote the space of all blocks $(\vec\bp, \vec\bb, \bH)\in \sK[\vec \bk]$, such that for each $\ell\in\llbracket 1,m\rrbracket$, $\bp^0$ is constant in $(\bQ_{\ell-1}, \bQ_{\ell-1}+\epsilon)$, and either $\bb_\ell>\epsilon$, or $\bb_\ell=\varnothing$.
We integrate over $\sK_\epsilon$, then send $\epsilon\to0+$.
\end{customthm}

\subsection{Recap of definitions from \cite{GXZ}}

The definitions of the combinatorial structure \emph{blocks}, which is the integration space in the formula of Theorem~\ref{Theorem_moments}, the measure on this space, and the function being integrated are given below.

\begin{definition}  \label{defn:dbfdef}
For an interval $[a, b]$, we call a function $f:[a, b]\to \R_{\ge 0}$ a \emph{block function}, if it is right-continuous on $[a,b)$, has finitely many left-discontinuities, and $f(a)=f(b)=f(b-)=0$, and the support of $f$ is connected, and the image $f([a,b])$ is a finite set.
We denote by $\fB([a, b])$ the space of all block functions $f$ on $[a,b]$.
\end{definition}
\begin{definition} \label{Definition_block_process}
Take $m\in \N$, $\vec \bk=(\bk_1,\ldots, \bk_m)\in \R_+^m$, and denote $\bQ_\ell=\sum_{\ell'=1}^\ell \bk_{\ell'}$ for each\footnote{Here and below, $\llbracket a, b\rrbracket=\{a,a+1,\dots,b\}$.} $\ell \in \llbracket 0, m\rrbracket$.
A \emph{block process with times $\vec \bk$}
is a family $\vec \bp=\{\bp_j\}_{j=1}^{m+\bu}$, where  $\bu\in \Z_{\ge 0}$, and each $\bp_j:[0, \bQ_m]\to \R_{\ge 0}$ is a function satisfying the following:
\begin{itemize}
    \item for each $\ell \in \llbracket 1, m\rrbracket$, $\bp_\ell\in \fB([0,\bQ_{m}])$, and $\bp_\ell= 0$ on $[\bQ_{\ell-1}, \bQ_m]$;
    \item for each $j\in \llbracket m+1, m+\bu\rrbracket$, $\bp_j \in \fB([0,\bQ_m])$ and $\bp_j$ is not identically zero.
    \item for each $j\in \llbracket 1, m+\bu\rrbracket$, $\ell \in \llbracket 1, m\rrbracket$, $\bp_j$ is left-continuous at $x=\bQ_{\ell-1}$ unless $j=\ell$.
    \item the regularity condition spelled in Definition \ref{defn:blp} below.
\end{itemize}
\end{definition}
In particular, the definition implies that $\bp_1$ is always identically zero. Given a block process and $\ell\in\llbracket 0, m\rrbracket$, we also denote
\begin{equation}\label{eq_p_ell}
\bp^\ell=\sum_{j=\ell+1}^{m+\bu} \bp_j.
\end{equation}

The block functions are encoded by the jumps at all the discontinuity points, for which we also set up notations. For a block process $\vec \bp$ as defined above, and  $j\in\llbracket 1, m+\bu\rrbracket$, $\ell\in\llbracket 1,m\rrbracket$, we let
\begin{equation}
\label{eq_discont_def_1}
 \bDel_{j,\ell}=\{x\in (\bQ_{\ell-1}, \bQ_\ell): \bp_j(x-)\neq \bp_j(x)\}, \qquad \bdel_{j,\ell}=|\bDel_{j,\ell}|\qquad  \text{ and}
\end{equation}
\begin{equation}
\label{eq_discont_def_2}
\bdel=\sum_{j=1}^{m+\bu}\sum_{\ell=1}^m \bdel_{j,\ell}, \quad  \bDel=\bigcup_{j=1}^{m+\bu}\bigcup_{\ell=1}^m \bDel_{j,\ell}.
\end{equation}

\begin{definition}  \label{defn:blp}
For a block process $\vec \bp$ as defined above, we say it is \emph{regular}, if all $\bDel_{j,\ell}$ are mutually disjoint and  $\min \bigcup_{\ell=1}^m\bDel_{j,\ell} < \min \bigcup_{\ell=1}^m\bDel_{j',\ell}$ for all $m+1\le j<j'\le m+\bu$. We consider only regular block processes.
\end{definition}
\begin{definition} \label{Definition_virtual_block}
For a block process $\vec \bp$, a compatible \emph{virtual block process} is a vector $\vec\bb=(\bb_1,\ldots, \bb_m)\in (\R_+\cup \{\varnothing\})^m$ satisfying the following conditions: for each $\ell \in \llbracket 1,m\rrbracket$,
\begin{itemize}
    \item if $\bb_\ell\neq \varnothing$, then $\bQ_{\ell-1}+\bb_\ell < \min\bigl( \bigcup_{j=1}^{m+\bu}\bDel_{j,\ell}\cup \{\bQ_\ell\}\bigr)$;
    \item if $\ell=1$, or $\ell>1$ and $\bp_\ell(\bQ_{\ell-1}-)=0$, then $\bb_\ell=\varnothing$.
\end{itemize}
\end{definition}
\begin{definition} \label{Definition_block_height}
For a block process $\vec \bp$ and compatible virtual block process $\vec\bb$, a compatible \emph{block height} is a function $\bH:\bDel\cup\{\bQ_\ell\mid \ell\in  \llbracket 0, m\rrbracket\}\cup\{\bQ_{\ell-1}+\bb_\ell \mid \ell\in  \llbracket 1, m\rrbracket,\, \bb_\ell\ne \varnothing\} \to \R_{\ge 0}$, such that
\begin{itemize}
    \item $\bH(\bQ_\ell)=\bp^0(\bQ_\ell-)$ for each $\ell \in \llbracket 1, m\rrbracket$ and $\bH(0)=0$;
    \item $\bH(\bQ_{\ell-1}+\bb_\ell)=\bp^0(\bQ_{\ell-1}-)$ for each $\ell$ with $\bb_\ell\neq\varnothing$;
    \item for all $x\in \bDel$, we have $\bH(x)\ge  \max\bigl(\bp^0(x),\bp^0(x-)\bigr)$;
    \item if $j\in \llbracket 1, m+\bu\rrbracket$ and $x\in \bDel$ is a discontinuity point for $\bp_j$, then the jump $\bp_j(x)-\bp_j(x-)$ and the difference $\bH(x)-\bp^0(x-)-\bp_j(x)$ have the same sign;
    \item if $\bb_\ell=\varnothing$ for $\ell\in \llbracket 1, m\rrbracket$, then  $\bH\bigl(\min\bigl(\bigcup_{j=1}^{m+\bu}\bDel_{j,\ell}\cup \{\bQ_\ell\}\bigr)\bigr)\ge \bH(\bQ_{\ell-1})$.
\end{itemize}
\end{definition}
In particular, if $\bp_j(x-)>0$ and $\bp_j(x)=0$ for $x\in\bDel$, then we must have $\bH(x)=\bp^0(x-)$, because $\bH(x)\ge \max\bigl(\bp^0(x), \bp^0(x-)\bigr)$ on one hand, and $\bH(x)-\bp^0(x-)\le \bp_j(x)$ on the other hand.

\begin{definition} ``Blocks'' is a triplet $(\vec\bp,  \vec\bb, \bH)$: block process, virtual block process, and compatible block height.
For each $\vec \bk\in \R_+^m$, we let $\sK=\sK[\vec\bk]$ be the space of all blocks with times $\vec\bk$.
\end{definition}

\begin{definition}  \label{defn:measurecb}
 The measure $\d(\vec\bp,  \vec\bb,\bH)$ on $\sK[\vec\bk]$ views it as a disjoint union of subspaces with fixed values of $\bu=0,1,\dots$ and $\bdel=0,1,\dots$. For the fixed $(\bu,\bdel)$ the space is further partitioned according to the values of $\bdel_{j,\ell}$.
 We equip each part with the product of Lebesgue measures over all parameters of blocks, except that whenever only a single value for a parameter is involved, then we assign weight $1$ for this value and for $\R_+\cup \{\varnothing\}$ in the definition of virtual blocks, $\varnothing$ comes with the Dirac delta measure of weight $1$, while $\R_+$ comes with the Lebesgue measure.
\end{definition}

The value of $\bL_\beta(\vec \bk, \vec\ttt)$ is an integral of a function $\bI_{\beta,\vec\bk}:\sK[\vec\bk]\to\R$. This function
involves expectations with respect to the following three stochastic processes. For $W:\R_{\ge 0}\to \R$ being a standard Wiener process,
\begin{itemize}
    \item for an interval $[a, b]$, and $g, h\in \R$, $x\mapsto B(x)= W(x-a)+g+\frac{x-a}{b-a}(-W(b-a) +h-g)$ is a Brownian bridge on $[a, b]$, with $B(a)=g$ and $B(b)=h$;
    \item for $W$ conditional on $W>-\epsilon$ on $[0, \epsilon^{-1}]$, its $\epsilon\to0+$ limit is the Bessel$_3$ process;
    \item for $a>0$, the Brownian excursion on $[0, a]$ is the $\epsilon\to0+$ limit of $W$ on $[0, a]$, conditional on $W>-\epsilon$ on $[0, a]$ and $W(a)<\epsilon$.
\end{itemize}

\stepcounter{theorem}

\begin{definition}   \label{defn:bIs}
For $x>0$, and $h, g\ge 0$, we define
\[
\bF(x;h,g) = \frac{1}{\sqrt{2\pi x}} \left( \exp\left(-\frac{(g-h)^2}{2x}\right)-\exp\left(-\frac{(g+h)^2}{2x}\right)\right),
\]
and
\[
\bI(x;h,g) = \bF(x;h,g)\E\left[ \exp\left( \beta^{-1}\int_0^x B(y) \d y \right) \right],
\]
where $B:[0, x]\to \R$ is a Brownian bridge with $B(0)=g$ and $B(x)=h$, conditional on $B\ge 0$. We also define
\begin{equation}
\label{eq_I0}
\bF_0(x;h) = \frac{2h}{\sqrt{2\pi x^3}} \exp\left(-\frac{h^2}{2x}\right), \qquad
\bI_0(x;h) =\bF_0(x;h)\E\left[ \exp\left( \beta^{-1}\int_0^x B_3(y) \d y \right) \right],
\end{equation}
where $B_3:[0, x]\to \R_{\ge 0}$ is a Bessel$_3$ process, conditional on $B_3(x)=h$.
Finally,
\begin{equation}
\label{eq_I00}
\bF_{0,0}(x)=\frac{2}{\sqrt{2\pi x^3}}, \qquad \bI_{0,0}(x) = \frac{2}{\sqrt{2\pi x^3}} \E\left[ \exp\left( \beta^{-1}\int_0^x B_e(y) \d y \right) \right],
\end{equation}
where $B_e:[0, x]\to \R_{\ge 0}$ is a Brownian excursion.
\end{definition}

\smallskip

Let us fix the blocks structure $(\vec\bp,
\vec\bb, \bH)\in \sK[\vec\bk]$ and consider the set of points on which $\bH$ is defined --- $\bDel\cup\{\bQ_\ell\mid \ell\in  \llbracket 0, m\rrbracket\}\cup\{\bQ_{\ell-1}+v_\ell \mid \ell\in  \llbracket 2, m\rrbracket,\, v_\ell\ne \varnothing\}$; we denote them $(0=x_0<x_1<x_2<\dots<x_n)$. We further let $\Xi$ be the set of $n$ pairs of adjacent points: $\Xi=\{(x_0,x_1), (x_1,x_2),\dots,(x_{n-1},x_n)\}$. We subdivide $\Xi$ into four subsets $\Xi_1,\Xi_2,\Xi_3,\Xi_4$, depending on the types of the constraints Definition \ref{Definition_block_height} imposed on $\bH$ in these points. For that, we introduce $A=\bp^0(x-)$ if $x=Q_\ell$, and otherwise $A=\bp^0(x)=\bp^0(y-)$; $A$ has the meaning of a lower bound for the Brownian motion in $[x, y]$. We compare possible values of $\bH(x)$ and $\bH(y)$ (dictated by Definition \ref{Definition_block_height}) with $A$: $(x,y)\in \Xi_i$ if
\begin{enumerate}
    \item[$\Xi_1:$] $\bH(x)=\bH(y)=A$. This happens when $x=\bQ_{\ell-1}$ for some $\ell \in \llbracket 1,m\rrbracket$, and either $y=\bQ_\ell$,  or $y=\bQ_{\ell-1}+\bb_\ell$, or $y=\sup\{y': \bp_j(y')>0\}$ for some $j\in\llbracket m+1, m+\bu\rrbracket$;
    \item[$\Xi_2:$] $\bH(x)=A$, $\bH(y)>A$. This happens when $x=\bQ_{\ell-1}$ for some $\ell \in \llbracket 1,m\rrbracket$, and $y\in \bDel_{j,\ell}$ with either $j\in\llbracket 1, m\rrbracket$ or  $j\in\llbracket m+1, m+\bu\rrbracket$ and $y\ne \sup\{y': \bp_j(y')>0\}$ in the latter case;
    \item[$\Xi_3:$] $\bH(x)>A$, $\bH(y)=A$. This happens when $x
    \not\in \{\bQ_\ell\}_{\ell=0}^{m-1}$, and either $y=\bQ_\ell$ for some $\ell \in \llbracket 1,m\rrbracket$, or $y=\sup\{y': \bp_j(y')>0\}$ for some $j\in\llbracket m+1, m+\bu\rrbracket$;
    \item[$\Xi_4:$] $\bH(x)>A$, $\bH(y)>A$. This happens when $x
    \not\in \{\bQ_\ell\}_{\ell=0}^{m-1}$, and $y\in \bDel_{j,\ell}$ with either $j\in\llbracket 1, m\rrbracket$ or  $j\in\llbracket m+1, m+\bu\rrbracket$ and $y\ne \sup\{y': \bp_j(y')>0\}$ in the latter case.
\end{enumerate}

\stepcounter{theorem}

\begin{definition}  \label{def:xibi} Using the above $\Xi$, $\Xi_1$, $\Xi_2$, $\Xi_3$, $\Xi_4$ notations, we define
\begin{equation} \label{eq:defnibk}
\bI_{\beta,\vec\bk}[\vec\bp,\vec\bb, \bH]=2^{-\bdel-|\{\ell\in\llbracket 1,m\rrbracket: \bb_\ell\neq\varnothing\}|}\prod_{(x,y)\in\Xi} \bI_\beta[x,y], \qquad \text{where}
\end{equation}
\begin{equation*}
\bI_\beta[x,y]=
\begin{cases}
\exp\big(\beta^{-1}(y-x)(\bH(x)-\bp^0(x))\big)\bI_{0,0}(y-x), & (x,y) \in \Xi_1, \\
\exp\big(\beta^{-1}(y-x)(\bH(x)-\bp^0(x))\big) \bI_0(y-x; \bH(y)-\bH(x)), & (x,y) \in \Xi_2, \\
(-1)^{\don[\bp^0(x)<\bp^0(x-)]}\bI_0(y-x;\bH(x)-\bH(y)) , & (x,y) \in \Xi_3, \\
(-1)^{\don[\bp^0(x)<\bp^0(x-)]} \bI(y-x;\bH(x)-\bp^0(x), \bH(y)-\bp^0(y_-) ), & (x,y) \in \Xi_4.
\end{cases}
\end{equation*}
\end{definition}

\subsection{The graph of connectivity}

In our arguments, we pass from moments described in Theorem \ref{Theorem_moments} to cumulants. This leads to cancellations and only \emph{connected} combinatorial structures survive. The connectivity of the blocks structure is described in terms of the graph which we now define. This graph does not depend on virtual block process or block height and is fully determined by a block process $ \vec\bp$.

\begin{definition} \label{Definition_connectivity_graph} Given a block process $\vec\bp$, its connectivity graph $\Gamma(\vec\bp)$ is a graph (undirected; with possible multi-edges and loops) with $m$ vertices and $\bdel-\bu=e$ edges. The vertex set is $\llbracket1,m\rrbracket$. For each $1\le \ell_1\le \ell_2\le m$ and each block  function $\bp_j$ of $\vec{\bp}$, such that two consecutive discontinuities of $\bp_j$ belong to intervals $[\bQ_{\ell_1-1},\bQ_{\ell_1})$ and $[\bQ_{\ell_2-1},\bQ_{\ell_2})$, we create an edge $\ell_1\to\ell_2$.
\end{definition}
Note that we included the left end-points in the intervals, which is important for the situation when $2\le j\le m$ and $\bp_j$ is not identically zero (no additional edges are created when $\bp_j\equiv 0$): in this case $\bp_j$ necessarily has a discontinuity at $\bQ_{j-1}$.

We call the blocks \((\vec\bp,\vec\bb,\bH)\) \emph{connected} when the graph $\Gamma(\vec\bp)$ is connected. In this situation, $e\ge m-1$.

\section{Cumulants via connected blocks}
\label{sec:basic}

If $(Y_i)_{i = 1}^n$ is an $n$-tuple of random variables, we define their joint cumulant through
\begin{equation}
\label{eq_cumulant_def}
\E^{(\textnormal{c})}\big[Y_1,\ldots,Y_n\big] = \ii^{-n} \partial_{s_1} \cdots \partial_{s_n} \log \E\Bigg[\exp\bigg(\ii \sum_{i= 1}^n s_i Y_i\bigg)\Bigg]\Bigg|_{s_1 = \ldots = s_n = 0}.
\end{equation}
Joint cumulants are related to joint moments through
\begin{equation}
\label{eq_moments_cumulants}
\E \bigg[\prod_{i = 1}^n Y_i \bigg] = \sum_{\mathbb{J} \vdash \llbracket 1, n\rrbracket} \prod_{J\, \in\, \mathbb{J}} \E^{(\textnormal{c})}\big[(Y_j)_{j \in J}\big], \qquad  \E^{(\textnormal{c})}\big[Y_1,\ldots,Y_n\big]
 =
  \sum_{\mathbb{J} \vdash \llbracket 1, n\rrbracket}
 (-1)^{|\mathbb{J}|-1}(|\mathbb{J}|-1)!
  \prod_{J\, \in\, \mathbb{J}} \E\left[\prod_{j\in J}Y_j\right],
\end{equation}
where the sums range over all set partitions of $\llbracket 1, n\rrbracket$, \textit{i.e.} tuples of nonempty pairwise disjoint subsets of $\llbracket 1, n\rrbracket$, and $|\mathbb{J}|$ is the number of sets in $\mathbb{J}$.  For more details on moments and cumulants see, \textit{e.g.}, \cite[Chapter 3]{PT}.

The goal of this section is to establish a version of Theorem \ref{Theorem_moments} for the cumulants.

\begin{definition}  \label{defn:Lbeta} For any $m\in \N$, $\vec \bk \in \R^m_+ $, $\vec \ttt\in \R^m$ such that $\ttt_1\le \cdots \le \ttt_m$, a graph $\Gamma$ on vertices $\llbracket1,m\rrbracket$, and reals $\epsilon,\beta>0$, we define
\begin{equation}
\label{eq_L_graph_cutoff}
 \bL^\Gamma_{\beta;\, (\epsilon)}(\vec\bk, \vec\ttt)
 =
 \int_{\substack{\sK_\epsilon[\vec\bk]\\
                   \Gamma(\vec\bp)=\Gamma}}
 \exp\left(\sum_{\ell=1}^{m-1}
(\ttt_\ell-\ttt_{\ell+1}) \bH(\bQ_\ell)/2 \right)\bI_{\beta,\vec\bk}[\vec \bp, \vec\bb, \bH] \d(\vec\bp,  \vec\bb,\bH),
\end{equation}
and
\begin{equation}\label{eq_L_graph}
  \bL^\Gamma_{\beta}(\vec\bk, \vec\ttt)=\lim_{\epsilon\to0+}\bL^\Gamma_{\beta;\, (\epsilon)}(\vec\bk, \vec\ttt).
\end{equation}
\end{definition}
The integral in Theorem~\ref{Theorem_moments} and its finite $\epsilon$ version are rewritten in the new notation as:
\begin{equation}
\label{eq_Graph_expansion}
 \bL_{\beta}(\vec\bk, \vec\ttt)=\sum_{\Gamma\text{ is a graph on }\llbracket1,m\rrbracket} \bL^\Gamma_{\beta}(\vec\bk, \vec\ttt), \qquad \qquad  \bL_{\beta;\, (\epsilon)}(\vec\bk, \vec\ttt)=\sum_{\Gamma\text{ is a graph on }\llbracket1,m\rrbracket} \bL^\Gamma_{\beta;\, (\epsilon)}(\vec\bk, \vec\ttt).
\end{equation}
For the cumulants, we should only sum over \emph{connected} graphs on $\llbracket1,m\rrbracket$. We define:
\begin{equation}
\label{eq_Graph_expansion_connected}
 \bLC_{\beta}(\vec\bk, \vec\ttt)=\sum_{\Gamma\text{ is connected}} \bL^\Gamma_{\beta}(\vec\bk, \vec\ttt), \qquad \qquad  \bLC_{\beta;\, (\epsilon)}(\vec\bk, \vec\ttt)=\sum_{\Gamma\text{ is connected}} \bL^\Gamma_{\beta;\, (\epsilon)}(\vec\bk, \vec\ttt).
\end{equation}
We note that some justifications are needed for these definitions: the integral in \eqref{eq_L_graph_cutoff} is actually absolutely convergent, the limit in \eqref{eq_L_graph} is well-defined, and the summations in \eqref{eq_Graph_expansion} are also absolutely convergent.
Such arguments are essentially present in \cite{GXZ}. The only issue is that the notion of graph $\Gamma(\vec\bp)$ is not present in \cite{GXZ}, but all the cancellations exploited in proving that the principal value in Theorem~\ref{Theorem_moments} exists act inside the classes of blocks with fixed $\Gamma$. For more details on justifying these definitions, see Appendix \ref{app:graph-beta-expansion}.

\begin{theorem}\label{Theorem_cumulants}
For any $m\in \N$, $\vec \bk \in \R^m_+ $, and $\vec \ttt\in \R^m$ such that $\ttt_1\le \cdots \le \ttt_m$, the joint $m$th cumulant for the Laplace transforms of $\{\AB_i\}_{i\ge 1}$ is
\begin{equation} \label{eq_expected_Laplace_cum}
\E^{(\textnormal{c})}\left[  \left(\sum_{i=1}^\infty
\exp(\bk_\ell \AB_i(\ttt_\ell)/2) \right)_{\ell=1}^m \right]=  \bLC_{\beta}(\vec\bk, \vec\ttt).
\end{equation}
\end{theorem}

 For the future asymptotic analysis it is convenient to expand the cumulant formula in Theorem~\ref{Theorem_cumulants} in a series in $\beta^{-1}$. For that we
note that $\bI_{\beta,\vec\bk}[\vec\bp,\vec\bb, \bH]$ of Definition \ref{def:xibi} depends on $\beta$ only through $\beta^{-1}$ prefactors in various exponents and therefore has a convergent large $\beta$ Taylor series expansion:
\begin{equation}
\label{eq_I_r}
 \bI_{\beta,\vec\bk}[\vec\bp,\vec\bb, \bH]=\sum_{r\ge 0} \beta^{-r} \cdot \bI_{\vec\bk}^r[\vec\bp,\vec\bb, \bH].
\end{equation}

\begin{definition} For any $m\in \N$, $\vec \bk \in \R^m_+ $, $\vec \ttt\in \R^m$ such that $\ttt_1\le \cdots \le \ttt_m$, a graph $\Gamma$ on vertices $\llbracket1,m\rrbracket$, a non-negative integer $r$, and a real $\epsilon>0$, we define
\begin{equation}
\label{eq_L_graph_cutoff_r}
 \bL^\Gamma_{[r];\, (\epsilon)}(\vec\bk, \vec\ttt)
 =
 \int_{\substack{\sK_\epsilon[\vec\bk]\\
                   \Gamma(\vec\bp)=\Gamma}}
 \exp\left(\sum_{\ell=1}^{m-1}
(\ttt_\ell-\ttt_{\ell+1}) \bH(\bQ_\ell)/2 \right)\bI^r_{\vec\bk}[\vec \bp, \vec\bb, \bH] \d(\vec\bp,  \vec\bb,\bH),
\end{equation}
and
\begin{equation}\label{eq_L_graph_r}
  \bL^\Gamma_{[r]}(\vec\bk, \vec\ttt)=\lim_{\epsilon\to0+}\bL^\Gamma_{[r];\, (\epsilon)}(\vec\bk, \vec\ttt).
\end{equation}
\end{definition}
As in Definition \ref{defn:Lbeta}, justifications are also needed for \eqref{eq_L_graph_cutoff_r} and \eqref{eq_L_graph_r}. These can be achieved by extending \eqref{eq_L_graph_cutoff} and \eqref{eq_L_graph} to complex $\beta$, and using Cauchy's integral formula to take contour integrals. See Appendix \ref{app:graph-beta-expansion} for the details.

\begin{lemma} \label{Lemma_Graph_expansion_beta}
There is $\epsilon_0>0$, such that for each graph $\Gamma$, the expansions
 \begin{equation}
 \label{eq_beta_expansion}
 \bL^\Gamma_{\beta}(\vec\bk, \vec\ttt)= \sum_{r=0}^\infty \beta^{-r}\cdot \bL^\Gamma_{[r]}(\vec\bk, \vec\ttt), \qquad
   \bL^\Gamma_{\beta;\, (\epsilon)}(\vec\bk, \vec\ttt)= \sum_{r=0}^{\infty} \beta^{-r}\cdot \bL^\Gamma_{[r];\, (\epsilon)}(\vec\bk, \vec\ttt)
 \end{equation}
 are absolutely convergent series, uniformly in $0<\epsilon<\epsilon_0$. Moreover, the total summations
  \begin{equation}
 \label{eq_beta_expansion_total}
 \begin{split}
     \bL_\beta(\vec\bk, \vec\ttt)&= \sum_{\Gamma\text{ is a graph on }\llbracket1,m\rrbracket}\sum_{r=0}^\infty \beta^{-r}\cdot \bL^\Gamma_{[r]}(\vec\bk, \vec\ttt), \\
   \bL_{\beta;\, (\epsilon)}(\vec\bk, \vec\ttt)&=\sum_{\Gamma\text{ is a graph on }\llbracket1,m\rrbracket} \sum_{r=0}^{\infty} \beta^{-r}\cdot \bL^\Gamma_{[r];\, (\epsilon)}(\vec\bk, \vec\ttt),
 \end{split}
 \end{equation}
  are also absolutely convergent series.
\end{lemma}
The proof of Lemma \ref{Lemma_Graph_expansion_beta} is an adaptation of the tail bounds
developed in \cite{GXZ} and is provided in Appendix \ref{app:graph-beta-expansion}. Plugging \eqref{eq_beta_expansion} into \eqref{eq_Graph_expansion_connected} we get a $\beta^{-1}$ expansion for the cumulants.

In the rest of this section we prove Theorem \ref{Theorem_cumulants}. The proof is based on a factorization identity.

\begin{lemma}\label{Lemma_factorization}
 Take a graph $\Gamma$ with vertices $\llbracket1,m\rrbracket$, represented as a union of connected components:
 $$
  \Gamma=\Gamma_1\sqcup \Gamma_2\sqcup\dots \sqcup \Gamma_d.
 $$
 For each $0<\epsilon<\min_{1\le i\le m} \bk_i$, we have
 \begin{equation} \label{eq_integral_factorization}
   \bL^\Gamma_{\beta;\, (\epsilon)}(\vec\bk, \vec\ttt)=\prod_{i=1}^d  \bL^{\Gamma_i}_{\beta;\, (\epsilon)}\left(\vec\bk\bigr|_{\Gamma_i}, \vec\ttt\bigr|_{\Gamma_i}\right)\qquad\text{ and }\qquad \bL^\Gamma_{\beta}(\vec\bk, \vec\ttt)=\prod_{i=1}^d  \bL^{\Gamma_i}_{\beta}\left(\vec\bk\bigr|_{\Gamma_i}, \vec\ttt\bigr|_{\Gamma_i}\right),
 \end{equation}
 where for $\Gamma_i$ with $m_i$ vertices, in $\bL^{\Gamma_i}$ we rename the vertices of $\Gamma_i$ into $\llbracket1,m_i\rrbracket$, and let $\vec\bk\bigr|_{\Gamma_i}$ and $\vec\ttt\bigr|_{\Gamma_i}$ be dimension $m_i$ vectors corresponding to the coordinates with indices of $\Gamma_i$ vertices.
\end{lemma}
\begin{proof} We prove the cutoff version of Lemma~\ref{Lemma_factorization}; the second version is obtained by \eqref{eq_L_graph}.

Let us take a blocks structure $(\vec\bp,\vec\bb,\bH)$ with $\Gamma(\vec\bp)=\Gamma$, and construct the corresponding $d$--tuple of blocks corresponding to $\Gamma_i$, $1\le i \le d$.

We start from the block process of Definition \ref{Definition_block_process}. Each $\bp_j$ with $1\le j\le m$ is assigned to the blocks of $\Gamma_i$ which has the vertex $j$. For $j>m$, $\bp_j$ is assigned to the blocks of $\Gamma_i$ which has the vertices of the intervals of discontinuities of $\bp_j$ (if $\bp_j$ has a discontinuity in $(\bQ_{\ell-1},\bQ_{\ell})$, then it is assigned go $\Gamma_i$ which has vertex $\ell$). Note that $\bp_j$ has discontinuities only in the segments $[\bQ_{\ell-1},\bQ_{\ell})$ of vertices of $\Gamma_i$ to which it was assigned, and is constant on all other segments. Hence, we can simply ignore these other segments, and get a block function entering  $\bL^{\Gamma_i}_{\beta;\, (\epsilon)}\left(\vec\bk\bigr|_{\Gamma_i}, \vec\ttt\bigr|_{\Gamma_i}\right)$. The $\bp^0$ for $\Gamma_i$ is reconstructed as the sum of all $\bp_j$ assigned to $\Gamma_i$.

Next, the virtual block process (see Definition \ref{Definition_virtual_block}) for $\Gamma_i$ has the coordinates equal to the coordinates of the virtual block process for $\Gamma$ corresponding to the vertices of $\Gamma_i$.

For the block height of Definition \ref{Definition_block_height}, we preserve not heights themselves, but rather the differences $\bH(x)-\bp^0(x-)$, with the convention $\bp^0(0-)=0$. At the end-points $\bQ_\ell$ of the intervals, the differences are set to $0$ both in $ \bL^\Gamma_{\beta;\, (\epsilon)}(\vec\bk, \vec\ttt)$ and in each  $\bL^{\Gamma_i}_{\beta;\, (\epsilon)}(\vec\bk, \vec\ttt)$. Other points where we define $\bH$ are discontinuities of various $\bp_j$ or virtual-block endpoints: these points are in bijection with the corresponding points for $\Gamma_i$, according to the assignments of the block and virtual block processes. At these interior points, preserving the differences is equivalent to setting
\begin{equation}
\label{eq_Height_shift}
 \bH_{\Gamma_i}=\bH_{\Gamma}-\bp^0_\Gamma+\bp^0_{\Gamma_i}.
\end{equation}
The construction can be also run in the opposite direction: given blocks for each $\Gamma_i$, we extend all $\bp_j$ by constants outside the segments of $\Gamma_i$, and then unite all $\bp_j$ into a single list. The virtual blocks for $\Gamma_i$ directly give rise to virtual blocks for $\Gamma$ and the correspondence for the heights is the same \eqref{eq_Height_shift} read in the opposite direction.

Let us check that various required properties of the blocks are preserved in the correspondence. For the conditions on the block process of Definition \ref{Definition_block_process}, the check is not immediate only for the first one. Equivalence of this condition for $\Gamma_i$ and for $\Gamma$ becomes clear once we notice that the $\bp_j$, $1\le j \le m$, for $\Gamma$ turns into $\bp_{j'}$ for $\Gamma_i$, where $j$ is $(j')$th smallest vertex of $\Gamma_i$. The conditions on virtual blocks in Definition \ref{Definition_virtual_block} are isolated inside a single interval $[\bQ_{\ell-1},\bQ_\ell)$ and therefore are the same between $\Gamma_i$ and $\Gamma$. The conditions on the block height of Definition \ref{Definition_block_height} are unchanged under simultaneous shift of $\bH$ and $\bp^0$ by a constant, and therefore are the same between $\Gamma_i$ and $\Gamma$.

We conclude that the integration space in the definition of $\bL^\Gamma_{\beta;\, (\epsilon)}(\vec\bk, \vec\ttt)$ splits into a direct product
of those for $\Gamma_i$, $1\le i \le d$. Since the measure of integration is a combination of Lebesgue and counting measures, it also factorizes. It remains to analyze the integrand in \eqref{eq_L_graph_cutoff}. Using $\bH(0)=\bH(\bQ_m)=0$, we rewrite
$$
\exp\left(\sum_{\ell=1}^{m-1}
(\ttt_\ell-\ttt_{\ell+1}) \bH(\bQ_\ell)/2 \right)= \exp\left(\sum_{\ell=1}^{m}
\ttt_\ell \bigl(\bH(\bQ_\ell)-\bH(\bQ_{\ell-1})\bigr)/2 \right),
$$
from which the factorization of this factor over $\Gamma_i$ is clear, because the increment $\bH(\bQ_\ell)-\bH(\bQ_{\ell-1})$ does not change under the shift \eqref{eq_Height_shift}. For the $\bI_{\beta,\vec\bk}[\vec\bp,\vec\bb, \bH]$ factor we note that classification into
the four cases in Definition \ref{def:xibi} remains exactly the same between $\Gamma_i$ and $\Gamma$, because $x$ and $y$ are always inside the same $[\bQ_{\ell-1},\bQ_{\ell}]$ segment. The factors corresponding to the four cases are also preserved, because $\bH$ and $\bp^0$ are shifted by the same constants, while in the $\Xi_3$ and $\Xi_4$ cases $x$ can not be one of $\bQ_\ell$ endpoints (and therefore the power of $(-1)$ is unchanged when passing from $\Gamma_i$ to $\Gamma$). Finally, the power of $2$ in \eqref{eq:defnibk} cleanly factorizes over $\Gamma_i$.

Applying Fubini's theorem proves Lemma~\ref{Lemma_factorization}.
\end{proof}

\begin{proof}[Proof of Theorem \ref{Theorem_cumulants}] We compute the cumulant in Theorem~\ref{Theorem_cumulants} by expressing it through moments using \eqref{eq_moments_cumulants} and Theorem \ref{Theorem_moments}, and then applying Lemma \ref{Lemma_factorization}. We get a sum over $d=1,2,\dots,m$, and disjoint unions of $d$ connected graphs $\Gamma_i$; the sum over graphs is infinite, but absolutely convergent, because so is \eqref{eq_Graph_expansion}. The $d=1$ term comes from $\mathbb{J}$ being a partition into a single set and matches the right-hand side of the identity in Theorem~\ref{Theorem_cumulants}. It remains to show that other terms cancel out.

For $d>1$, take a $d$--tuple of connected graphs $\Gamma_i$ such that their vertex sets are disjoint and have the union $\llbracket 1,m\rrbracket$. The product corresponding to these graphs,
$$
\prod_{i=1}^d  \bL^{\Gamma_i}_{\beta}\left(\vec\bk\bigr|_{\Gamma_i}, \vec\ttt\bigr|_{\Gamma_i}\right),
$$
arises with the coefficient
$$
 \sum_{\mathbb{J} \vdash (\Gamma_1,\dots,\Gamma_d)}
 (-1)^{|\mathbb{J}|-1}(|\mathbb{J}|-1)!=\sum_{\mathbb{J'} \vdash \llbracket 1,d\rrbracket}
 (-1)^{|\mathbb{J}'|-1}(|\mathbb{J}'|-1)!,
$$
where the first summation goes over all set partitions of $\llbracket 1,m\rrbracket$ whose parts are unions of the vertex sets of $\Gamma_1,\dots,\Gamma_d$, and the second summation goes over all set partitions of $\llbracket 1,d\rrbracket$. The second sum is equal to zero, because it computes the $d$--th cumulant of the (deterministic) random variable $1$ by \eqref{eq_moments_cumulants}.
\end{proof}

\section{Rescaling of the graph expansion}
\label{Section_scaling}

For the GFF limit, we need the following rescaling of the parameters:
\begin{equation}
 (\vec\bk,\vec\ttt)
 \longmapsto
 (\eta\vec\bk,\eta^{-1/2}\vec\ttt), \qquad \eta\to0+.
\end{equation}
The way the cumulant expansion obtained from Theorem~\ref{Theorem_cumulants}, \eqref{eq_Graph_expansion_connected}, and \eqref{eq_beta_expansion} changes under this rescaling is described in the following proposition.

\begin{proposition} \label{Proposition_graph_rescaling} For any $m\in \N$, $\vec \bk \in \R^m_+ $, $\vec \ttt\in \R^m$ such that $\ttt_1\le \cdots \le \ttt_m$, $\beta>0$, $r=0,1,\dots$, and a graph $\Gamma$ on vertices $\llbracket1,m\rrbracket$ with $e$ edges, we have
\begin{equation}
\label{eq_Graph_scaling}
 \bL^\Gamma_{[r]}\bigl(\eta \vec\bk, \eta^{-1/2}\vec\ttt\bigr)=  \eta^{\frac32(e+r-m)} \cdot \bL^\Gamma_{[r]}\bigl(\vec\bk, \vec\ttt\bigr).
\end{equation}
\end{proposition}
\begin{proof} We make the following change of variables in \eqref{eq_L_graph_cutoff}:
\begin{equation}\label{eq_rescaling_integral}
  \vec\bk
 \mapsto
 \eta\vec\bk, \qquad \vec\ttt\mapsto \eta^{-1/2}\vec\ttt, \qquad x\mapsto\eta x,\qquad
 \bp_j\mapsto \eta^{1/2}\,\bp_j,\qquad
 \bH\mapsto\eta^{1/2}\,\bH,\qquad
 \bb_\ell\mapsto\eta\bb_\ell,
\end{equation}
where $x$ is the horizontal coordinate for the discontinuities of block process $\vec\bp$. The Jacobian of the transformation \eqref{eq_rescaling_integral}
leads to the factor:
\begin{equation}
\label{eq_x1}
 \eta^{\bdel} \cdot \eta^{|\{\ell\in\llbracket1,m\rrbracket:\bb_\ell\neq\varnothing\}|} \cdot \eta^{(\bdel-\bu)/2} \cdot \eta^{(\bdel-\bu)/2}=\eta^{2\bdel-\bu+|\{\ell\in\llbracket1,m\rrbracket:\bb_\ell\neq\varnothing\}|},
\end{equation}
where the factors correspond to positions of discontinuities of block process, the coordinates of the virtual block process, the sizes
of jumps of the block process, and the (free) coordinates of the block height $\bH$.

The exponent in the integrand in \eqref{eq_L_graph_cutoff} is unchanged under $\eta$--rescaling. For the $\bI^r_{\eta\vec\bk}[\eta^{1/2} \vec \bp, \eta \vec\bb, \eta^{1/2}\bH]$ factor in the integrand, we note that the
Brownian scaling multiplies \(\bI_{0,0},\bI_0,\bI\) by
\(\eta^{-3/2},\eta^{-1},\eta^{-1/2}\), respectively, and rescales the $\beta^{-1}$ prefactor in the exponents by $\eta^{3/2}$. In order to compute the total change of the product in Definition \ref{def:xibi} and \eqref{eq_I_r}, we need to count the numbers of members in $\Xi_1,\Xi_2,\Xi_3,\Xi_4$ cases. We have:
\begin{enumerate}
\item The first two cases correspond to $x=\bQ_{\ell-1}$. Hence, $|\Xi_1|+|\Xi_2|=m$.
\item The second and fourth cases correspond to $y$ being a discontinuity which is not the last one for the corresponding $\bp_j$. Hence, $|\Xi_2|+|\Xi_4|=\bdel-\bu$.
\item Overall, $\bH$ is defined at $n+1=\bdel+m+1+|\{\ell\in\llbracket1,m\rrbracket:\bb_\ell\neq\varnothing\}|$ points. Hence, $|\Xi_1|+|\Xi_2|+|\Xi_3|+|\Xi_4|=\bdel+m+|\{\ell\in\llbracket1,m\rrbracket:\bb_\ell\neq\varnothing\}|$.
\end{enumerate}
We conclude that
\begin{multline}
\label{eq_x2}
  \bI^r_{\eta\vec\bk}[\eta^{1/2} \vec \bp, \eta \vec\bb, \eta^{1/2}\bH]= \eta^{-\frac{1}{2}\left(3|\Xi_1|+2|\Xi_2|+2|\Xi_3|+|\Xi_4|\right)} \cdot \eta^{\frac{3}{2}r} \cdot \bI^r_{\vec\bk}[\vec \bp, \vec\bb, \bH]\\
  = \eta^{-\frac{1}{2}\bdel-\frac{1}{2}\bu-\frac{3}{2}m-|\{\ell\in\llbracket1,m\rrbracket:\bb_\ell\neq\varnothing\}| + \frac{3}{2}r} \cdot \bI^r_{\vec\bk}[\vec \bp, \vec\bb, \bH].
\end{multline}
Multiplying \eqref{eq_x1} by \eqref{eq_x2}, recalling that $\bdel-\bu=e$ in Definition \ref{Definition_connectivity_graph}, and sending $\epsilon\to0+$ we obtain
Proposition~\ref{Proposition_graph_rescaling}.
\end{proof}

\begin{corollary} \label{Corollary_cumulant_expansion}
 For any $m\in \N$, $\vec \bk \in \R^m_+ $, and $\vec \ttt\in \R^m$ such that $\ttt_1\le \cdots \le \ttt_m$, and $\eta>0$, the rescaled joint $m$th cumulant for the Laplace transforms of $\{\AB_i\}_{i\ge 1}$ is
\begin{multline} \label{eq_expected_Laplace_cum_rescaled}
\E^{(\textnormal{c})}\left[  \left(\sum_{i=1}^\infty
\exp\left(\eta \bk_\ell \AB_i(\eta^{-1/2}\ttt_\ell)/2\right) \right)_{\ell=1}^m \right]\\ =
 \sum_{e,r=0}^\infty  \beta^{-r}  \eta^{\frac32(e+r-m)}  \cdot \sum_{\begin{smallmatrix}\Gamma\text{ is a connected graph on } \llbracket 1,m\rrbracket \\ \text{with }e\text{ edges}\end{smallmatrix}} \bL^\Gamma_{[r]}(\vec\bk, \vec\ttt),
\end{multline}
where the series is absolutely convergent.
\end{corollary}
\begin{proof} Proposition \ref{Proposition_graph_rescaling} gives  the $\eta^{\frac32(e+r-m)}$ rescaling of the individual terms of the series delivered by combining Theorem \ref{Theorem_cumulants} with expansions \eqref{eq_beta_expansion}. The absolute convergence of the series for fixed $\eta$ follows from \eqref{eq_beta_expansion_total}.
\end{proof}

\section{Limit of expectation}

In this section, as a warm-up for the later analysis of higher cumulants, we use Corollary \ref{Corollary_cumulant_expansion} with $m=1$ to get the first term of the asymptotic expansion of the mean of the Airy surface as $\eta\to0+$. Since the ensemble is stationary in time $\ttt$, it is sufficient to only deal with the $\ttt=0$ case.

\begin{theorem} \label{Theorem_mean}
 For any real $\bk>0$, as $\eta\to0+$
\begin{equation}
\label{eq_mean_expansion}
 \E\left[  \sum_{i=1}^\infty
\exp\left( \frac{\eta \bk}{2} \AB_i(0)\right)  \right]= \sqrt{\frac{2}{\pi}} \bk^{-3/2} \eta^{-3/2} + \frac{1}{2}\left(\frac{1}{\beta}-\frac{1}{2}\right) + O\left(\eta^{3/2}\right).
\end{equation}
\end{theorem}
There are other ways to get this result, e.g., it follows from \cite[Propositon 2.14, Corollary 2.15]{GS} and the fact (see \cite{Janson}) that the expected area under Brownian excursion is $\sqrt{\frac{\pi}{8}}$.

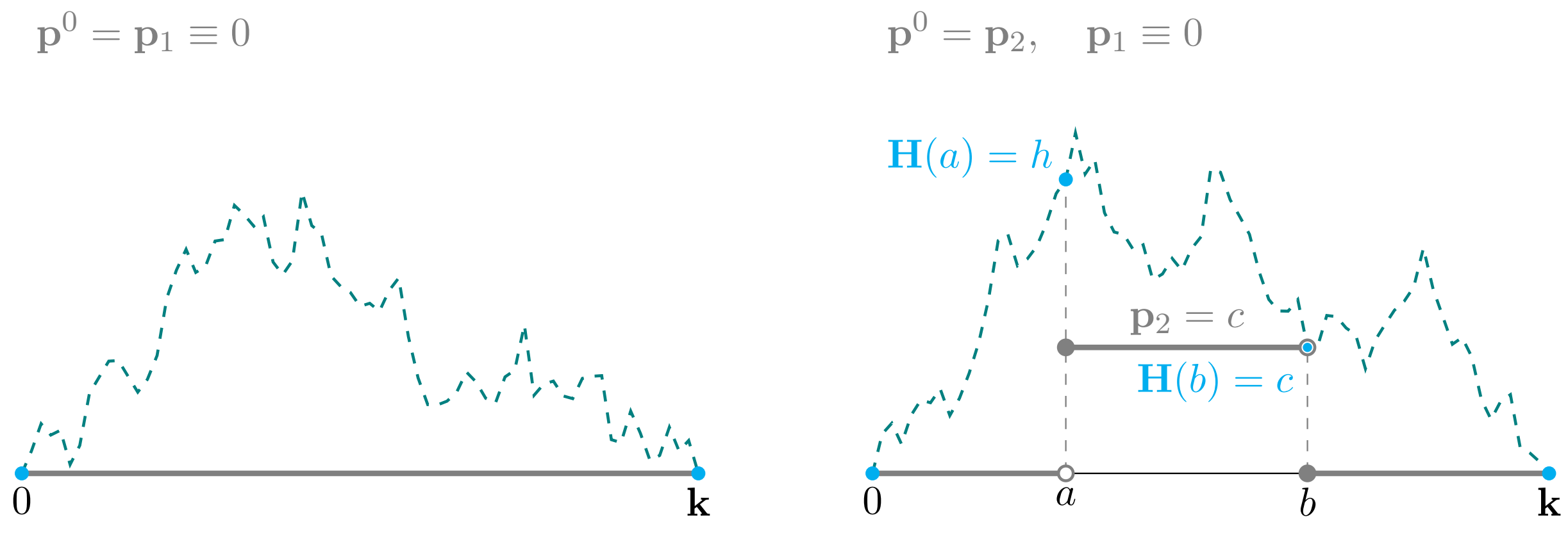
\begin{figure}[!t]
  \begin{subfigure}[b]{0.46\textwidth}
    \centering
    \scalebox{1.084}{%
    \begin{tikzpicture}[x=0.46cm,y=1cm]
    \draw[gray,line width=1.5pt] (0,0)--(14,0);
    \draw[teal,thick,dashed] plot coordinates {
      (0,0) (0.2,0.213) (0.4,0.473) (0.6,0.365) (0.8,0.411) (1,0.084)
      (1.2,0.266) (1.4,0.762) (1.6,0.916) (1.8,1.069) (2,1.075) (2.2,0.924)
      (2.4,0.774) (2.6,0.896) (2.8,1.123) (3,1.668) (3.2,1.936) (3.4,2.134)
      (3.6,1.914) (3.8,1.97) (4,2.213) (4.2,2.227) (4.4,2.553) (4.6,2.465)
      (4.8,2.355) (5,2.444) (5.2,2.013) (5.4,1.892) (5.6,2.026) (5.8,2.68)
      (6,2.36) (6.2,2.286) (6.4,1.88) (6.6,1.78) (6.8,1.723) (7,1.591)
      (7.2,1.62) (7.4,1.544) (7.6,1.741) (7.8,1.862) (8,1.302) (8.2,0.91)
      (8.4,0.655) (8.6,0.652) (8.8,0.688) (9,0.774) (9.2,0.966) (9.4,0.871)
      (9.6,0.715) (9.8,0.668) (10,0.919) (10.2,0.983) (10.4,1.415) (10.6,0.744)
      (10.8,0.851) (11,0.879) (11.2,0.738) (11.4,0.712) (11.6,0.905) (11.8,0.925)
      (12,0.931) (12.2,0.321) (12.4,0.293) (12.6,0.593) (12.8,0.383) (13,0.122)
      (13.2,0.171) (13.4,0.44) (13.6,0.214) (13.8,0.313) (14,0)
    };

    \foreach \x in {0,14}
      \draw[cyan,fill=cyan,thick] (\x,0) circle (1.5pt);
    \node[below] at (0,0) {$0$};
    \node[below] at (14,0) {$\bk$};

    \node[anchor=west,text=gray] at (0,4.2) {$\bp^0=\bp_1\equiv0$};
  \end{tikzpicture}%
    }
    \caption{No jumps, $\bu=0$.}
  \end{subfigure}%
  \hfill
  \begin{subfigure}[b]{0.46\textwidth}
    \centering
    \scalebox{1.084}{%
    \begin{tikzpicture}[x=0.46cm,y=1cm]
    \draw (0,0)--(14,0);
    \draw[gray,dashed] (4,0)--(4,2.8);
    \draw[gray,dashed] (9,0)--(9,1.2);

    \draw[gray,line width=1.5pt] (0,0)--(4,0);
    \draw[gray,line width=1.5pt] (4,1.2)--(9,1.2);
    \draw[gray,line width=1.5pt] (9,0)--(14,0);

    \draw[teal,thick,dashed] plot coordinates {
      (0,0) (0.2,0.375) (0.4,0.478) (0.6,0.279) (0.8,0.554) (1,0.699)
      (1.2,0.673) (1.4,0.806) (1.6,0.551) (1.8,0.716) (2,0.954) (2.2,1.249)
      (2.4,1.656) (2.6,2.217) (2.8,2.278) (3,1.977) (3.2,2.041) (3.4,2.163)
      (3.6,2.382) (3.8,2.665) (4,2.8) (4.2,3.243) (4.4,2.851) (4.6,2.988)
      (4.8,2.483) (5,2.299) (5.2,2.275) (5.4,2.128) (5.6,2.177) (5.8,1.844)
      (6,1.899) (6.2,2.048) (6.4,1.934) (6.6,2.139) (6.8,2.255) (7,2.9)
      (7.2,2.868) (7.4,2.607) (7.6,2.445) (7.8,2.277) (8,1.922) (8.2,1.657)
      (8.4,1.547) (8.6,1.545) (8.8,1.655) (9,1.2) (9.2,1.193) (9.4,1.504)
      (9.6,1.493) (9.8,1.385) (10,1.319) (10.2,0.996) (10.4,1.267) (10.6,1.414)
      (10.8,1.554) (11,1.638) (11.2,1.776) (11.4,2.146) (11.6,1.756) (11.8,1.495)
      (12,1.23) (12.2,1.294) (12.4,1.11) (12.6,0.702) (12.8,0.513) (13,0.688)
      (13.2,0.75) (13.4,0.251) (13.6,0.195) (13.8,0.112) (14,0)
    };

    \draw[gray,fill=white,thick] (4,0) circle (2pt);
    \draw[gray,fill=gray,thick] (4,1.2) circle (2pt);
    \draw[gray,fill=white,thick] (9,1.2) circle (2pt);
    \draw[gray,fill=gray,thick] (9,0) circle (2pt);

    \foreach \x/\y in {0/0,4/2.8,14/0}
      \draw[cyan,fill=cyan,thick] (\x,\y) circle (1.5pt);
    \fill[cyan] (9,1.2) circle (1.25pt);
    \node[above right,text=cyan,yshift=8pt] at (0,2.4) {$\bH(a)=h$};
    \node[below left,text=cyan] at (9,1.2) {$\bH(b)=c$};

    \node[below] at (0,0) {$0$};
    \node[below] at (4,0) {$a$};
    \node[below] at (9,0) {$b$};
    \node[below] at (14,0) {$\bk$};
    \node[below,text=gray] at (6.5,1.70) {$\bp_2=c$};
    \node[anchor=west,text=gray] at (0,4.2) {$\bp^0=\bp_2,\quad\bp_1\equiv0$};
  \end{tikzpicture}%
    }
    \caption{Two jumps, $\bu=1$.}
  \end{subfigure}
  \caption{The blocks for the $\Gamma_0$ and $\Gamma_1$ cases.}
  \label{Figure_one-interval-block-cases}
\end{figure}

Before proving Theorem \ref{Theorem_mean} we state a helpful property for the function $\bF_0$ of \eqref{eq_I0}.

\begin{lemma}\label{lem:FP_convolution}
The function $\bF_0$ satisfies a convolution property:
\begin{equation}\label{eq_FP_convolution}
 \int_0^T \bF_0(t;u)\bF_0(T-t;v)\dd t
 =
 2\bF_0(T;u+v),
 \qquad u,v,T>0.
\end{equation}
\end{lemma}
\begin{proof}
This can be proven by direct integration. Alternatively, we can interpret $\tfrac{1}{2}\bF_0(t;u)\dd t$ as the probability density for the random variable $\tau_u=\inf\{t\ge 0:\, B(t)=u\}$, where $B(t)$ is the standard Brownian motion, and then use the identity $\tau_{u+v}=\tau_u+\tilde \tau_v$, which follows from the strong Markov property for the Brownian motion and where $\tau_u$ and $\tilde \tau_v$ are independent.
\end{proof}

\begin{proof}[Proof of Theorem \ref{Theorem_mean}]
We send $\eta\to0+$ in the $m=1$ case of Corollary~\ref{Corollary_cumulant_expansion}. Since the series is absolutely convergent, we can send $\eta\to0+$ term by term. The $e=r=0$ term gives the first term of the expansion in Theorem~\ref{Theorem_mean}, as soon as we prove:
\begin{equation}
\label{eq_x3}
\bL^{\Gamma_0}_{[0]}(\bk, 0) \stackrel{?}{=}\sqrt{\frac{2}{\pi}} \bk^{-3/2}, \qquad \Gamma_0= \text{graph with one vertex and no edges}
\end{equation}
Indeed, for this graph $\bdel=\bu=0$, there are no non-zero block functions or virtual blocks. A single pair $(x,y)$ in Definition \ref{def:xibi} is $(x,y)=(0,\bk)$ of type $\Xi_1$, as in the left panel of Figure \ref{Figure_one-interval-block-cases}, and $\bI^0_{\bk}$ evaluates using \eqref{eq_I00} to $\bF_{0,0}(\bk)=\frac{2}{\sqrt{2\pi \bk^3}}$, matching \eqref{eq_x3}.

The second term of the expansion in Theorem~\ref{Theorem_mean} is a combination of $(e,r)=(0,1)$ and $(e,r)=(1,0)$ terms in the expansion of Corollary~\ref{Corollary_cumulant_expansion}, and we should prove:
\begin{equation}
\label{eq_x4}
\bL^{\Gamma_0}_{[1]}(\bk, 0) \stackrel{?}{=}\frac{1}{2}, \qquad \bL^{\Gamma_1}_{[0]}(\bk, 0) \stackrel{?}{=}-\frac{1}{4}, \qquad \Gamma_1= \text{graph with one vertex and one edge (loop)}.
\end{equation}
For the first identity, compared to \eqref{eq_x3} we replace the $0$--th term of $\frac{1}{\beta}$ expansion of $\Xi_1$  by the first one. Looking back at \eqref{eq_I00}, we therefore need to multiply \eqref{eq_x3} by the expected area under Brownian excursion going from $0$ to $0$ in time $\bk$. This expectation is $\bk^{3/2} \sqrt{\frac{\pi}{8}}$, see, e.g., \cite{Janson}, and we get the first equality of \eqref{eq_x4}.

For the second identity, we deal with $\Gamma_1$. The corresponding blocks structure has $\bdel=2$, $\bu=1$, as in the right panel of Figure \ref{Figure_one-interval-block-cases}, i.e., $\bp_1\equiv 0$, while the function $\bp_2$ has the form
$$
 \bp_2(x)=\begin{cases} c, & a\le x<b,\\ 0,& 0\le x<a\quad \text{ or }\quad b\le x \le \bk, \end{cases}
$$
where $0<a<b<\bk$ and $c>0$ are three parameters over which we will be integrating. There is no virtual block. The heights $\bH$ are prescribed at points $0,b,\bk$ to be $0,c,0$, respectively; one more value $\bH(a)=h>c$ is the fourth variable of integration. Altogether, we need to compute a four-dimensional integral. Since we deal with $r=0$ term, all the exponents involving $\beta$ in Definitions \ref{defn:bIs} and \ref{def:xibi} can be ignored. The three intervals in Definition \ref{def:xibi} are $(0,a)$, $(a,b)$, $(b,\bk)$ of types $\Xi_2$, $\Xi_3$, $\Xi_3$, respectively. Hence, we need to compute the integral:
\begin{equation}\label{eq_four_tuple_integral}
 \iint_{0<a<b<\bk} \iint_{0<c<h} 2^{-2} \bF_0(a;h) \cdot \bF_0(b-a;h-c) \cdot (-1) \bF_0(\bk-b;c) \dd a\dd b \dd c \dd h,
\end{equation}
where $\bF_0$ is given by \eqref{eq_I0}. Applying Lemma~\ref{lem:FP_convolution} twice to integrate over $a$ and $b$ in \eqref{eq_four_tuple_integral}, we transform it into
$$
-\iint_{0<c<h} \bF_0(\bk;2h)\dd c \dd h=-\int_0^{\infty}  \frac{4h^2}{\sqrt{2\pi \bk^3}} \exp\left(-2 \frac{h^2}{\bk}\right)\dd h=
-\frac{1}{2\sqrt{\pi}}\int_0^{\infty}  z^{1/2} \exp\left(-z\right) \dd z=-\frac{1}{4},
$$
where we changed the variable $z=2h^2/\bk$ and used $\Gamma(\tfrac{3}{2})=\tfrac{\sqrt{\pi}}{2}$. We obtain the desired \eqref{eq_x4}.
\end{proof}
\begin{remark}
  Probabilistically, the integral \eqref{eq_four_tuple_integral}  represents expected area under the Brownian excursion (up to multiplication with an explicit constant). This excursion is obtained by concatenating three parts on $[0,a]$, $[a,b]$, $[b,\bk]$ intervals, as in the right panel of Figure \ref{Figure_one-interval-block-cases}. We start by fixing a trajectory of an excursion. Integration over $a$ and $h$ specifies that this excursion passes through the point $(a,h)$; $c$ then can be integrated over the segment $[0,h]$, creating an $h$ factor; once $(a,h,c)$ are fixed, $b$ is uniquely determined as the first time excursion reaches level $c$ after time $a$. Hence, we are integrating $h$ over all possible values of $(a,h)$ on the excursion trajectory, which gives the area.
\end{remark}
We can translate Theorem~\ref{Theorem_mean} to asymptotics of the Airy surface expectation. We recall the notation  $\cA^\beta(\ttt,x) =\#\{i\in\N:\cA_i^\beta(\ttt)>x\}$.
\begin{corollary}\label{cor:height-mean}
For every fixed $\beta>0$, as $x\to+\infty$,
\begin{equation}
\label{eq_x13}
 \E[\cA^\beta(0,-x)]
 =\frac{2}{3\pi}x^{3/2}+o(x^{3/2}).
\end{equation}
\end{corollary}
\begin{proof}
For every $q>0$, $x\in\mathbb R$, the definition of $\cA^\beta(\ttt,x)$ gives
\begin{equation}   \label{eq:AbqA}
 0\le\E[\cA^\beta(0,-x)]
 \le e^{qx} \cdot \E\sum_{i=1}^\infty \exp\bigl(q\cA_i^\beta(0)\bigr).
\end{equation}
The expectation in the right-hand side is finite by
Theorem~\ref{Theorem_moments} in Section \ref{Section_moments_from_GXZ}. For any $\eta>0$,
taking $q=\eta/2$ as $x\to+\infty$ and $q=2\eta$ as
$x\to-\infty$ shows that, whenever $|x|\to \infty$,
\begin{equation}   \label{eq:boundaryterm}
 e^{-\eta x}\E[\cA^\beta(0,-x)]\to 0.
\end{equation}
Thus, integrating by parts and using \eqref{eq:boundaryterm} for the boundary terms, we get
\[
 \eta\int_{-\infty}^{+\infty} e^{-\eta x}\E[\cA^\beta(0,-x)]\dd x
 =\E\sum_{i=1}^\infty \exp\bigl(\eta\cA_i^\beta(0)\bigr).
\]
Next, we use \eqref{eq:AbqA} with $q=1$ and get that for $0<\eta<1$,
\[
 0\le\eta\int_{-\infty}^0e^{-\eta x}\E[\cA^\beta(0,-x)]\dd x
 \le\frac{\eta}{1-\eta}
       \E\sum_{i=1}^\infty \exp\bigl(\cA_i^\beta(0)\bigr)
 =O(\eta), \qquad \eta\to 0+.
\]
The last two equations and Theorem~\ref{Theorem_mean} with $\bk=2$ give that
\begin{equation}
\label{eq_x10}
 \int_0^\infty e^{-\eta x}\E[\cA^\beta(0,-x)]\dd x
 =\frac{1}{2\sqrt\pi}\eta^{-5/2} + O(\eta^{-1}), \qquad \eta\to 0+.
\end{equation}
Now write
\[
 m^{(\eta)}_\beta(x)=\eta^{3/2}\E[\cA^\beta(0,-\eta^{-1}x)],
 \qquad x\ge0.
\]
For every $s>0$, the estimate \eqref{eq_x10} gives
\begin{equation}   \label{eq:intsLmeta}
 \lim_{\eta\to 0+}\int_0^\infty e^{-sx}m^{(\eta)}_\beta(x)\dd x
 =\lim_{\eta\to 0+} \eta^{5/2}\int_0^\infty e^{-s\eta x}\E[\cA^\beta(0,-x)]\dd x
 =
 \frac{s^{-5/2}}{2\sqrt\pi}
 =\frac{2}{3\pi} \int_0^\infty e^{-sx}x^{3/2}\dd x.
\end{equation}
Let us show that we can replace $e^{-sx}$ in the above formula by any continuous compactly supported
function $\phi:[0,\infty)\to\R$. In other words, we claim that
\begin{equation}
\label{eq_x11}
 \lim_{\eta\to 0+}\int_0^\infty\phi(x)m^{(\eta)}_\beta(x)\dd x
 \stackrel{?}{=} \frac{2}{3\pi}\int_0^\infty\phi(x)x^{3/2}\dd x.
\end{equation}
Indeed, the function
$$
q\mapsto \begin{cases} q^{-1}\phi(-\log q), &q>0,\\ 0, &q=0 \end{cases}
$$
is continuous on $[0,1]$. Hence, it can be approximated by a polynomial. Writing $q=e^{-x}$, we conclude that for any $\varepsilon>0$ there exists a polynomial $P_\eps$ such that
\[
 |\phi(x)-e^{-x}P_\eps(e^{-x})|\le\varepsilon e^{-x},\qquad \forall\,x\ge0.
\]
A finite linear combination of  \eqref{eq:intsLmeta} gives
\[
 \lim_{\eta\to 0+} \int_0^\infty e^{-x}P_\eps(e^{-x})m^{(\eta)}_\beta(x)\dd x
 = \frac{2}{3\pi}\int_0^\infty e^{-x}P_\eps(e^{-x})x^{3/2}\dd x.
\]
Sending $\eps\to 0$, we get \eqref{eq_x11}.

For $0<\delta<1$, take $\phi$ to be two nonnegative continuous functions of integral one
supported in $(1-\delta,1)$ and $(1,1+\delta)$, respectively. Using \eqref{eq_x11} and monotonicity of $m^{(\eta)}_\beta$ we obtain
\[
  \frac{2}{3\pi}(1-\delta)^{3/2}
 \le\liminf_{\eta\to0+}m^{(\eta)}_\beta(1)
 \le\limsup_{\eta\to0+}m^{(\eta)}_\beta(1)
 \le  \frac{2}{3\pi}(1+\delta)^{3/2}.
\]
Sending $\delta\to0+$ proves the claimed asymptotic \eqref{eq_x13} for $x=\eta^{-1}$.
\end{proof}
\begin{remark}
The above argument uses nonnegativity and monotonicity in $x$ of $\E[\cA^\beta(0,-x)]$. For the remainder term
$\E[\cA^\beta(0,-\eta^{-1}x)]-\frac{2}{3\pi}\eta^{-3/2}x^{3/2}$, we do not know such properties \emph{a priori}, which is an obstacle for recovering the constant order term $\frac{1}{2\beta}-\frac{1}{4}$ in \eqref{eq:EAconst} from
Laplace asymptotics in Theorem~\ref{Theorem_mean}. However, if we knew that an $\eta\to 0$ asymptotic expansion of the form \eqref{eq:EAconst} exists, then arguments similar to those in Corollary \ref{cor:height-mean} can be used to obtain the constant order term.
\end{remark}

\section{Three cancellation identities}   \label{sec:threecan}

Our next goal is to analyze the expansion in Corollary~\ref{Corollary_cumulant_expansion} for $m\ge 2$. We expect that the second cumulant is of constant order as $\eta\to0+$, while higher cumulants go to $0$. Hence, we need to show that the relevant leading coefficients are equal to $0$. Since the graphs in the sum of Corollary \ref{Corollary_cumulant_expansion} are connected, we have
\begin{equation}
 e\ge m-1.
\end{equation}
The following three propositions address the minimal  $e=m-1$ and $e=m$ cases.

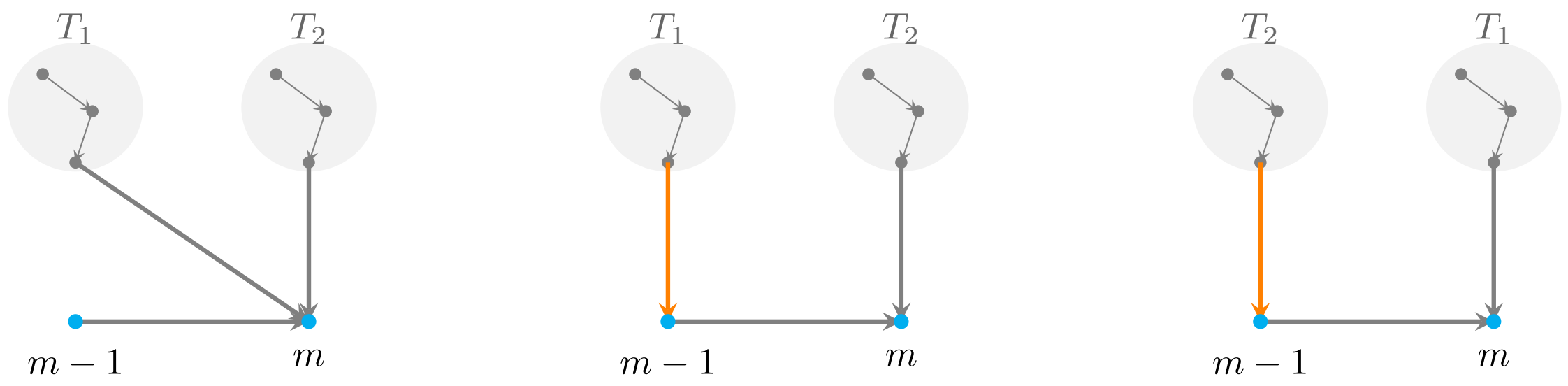
\begin{figure}[t]
\centering
\begin{subfigure}[t]{.31\textwidth}
\centering
\begin{tikzpicture}[x=.94cm,y=.78cm,font=\small,>=stealth]
\path[use as bounding box] (-.25,-.3) rectangle (4.4,4.65);
\fill[gray!10] (0.85,3.25) ellipse (.68 and .78);
\node[text=gray!80!black] at (0.85,4.2) {$T_1$};
\draw[gray,->] (0.52,3.65)--(1.02,3.2);
\draw[gray,->] (1.02,3.2)--(0.85,2.58);
\fill[gray] (0.52,3.65) circle (1.65pt);
\fill[gray] (1.02,3.2) circle (1.65pt);
\fill[gray] (0.85,2.58) circle (1.65pt);
\draw[gray,->,line width=1.2pt] (0.85,2.58)--(3.2,.65);
\fill[gray!10] (3.2,3.25) ellipse (.68 and .78);
\node[text=gray!80!black] at (3.2,4.2) {$T_2$};
\draw[gray,->] (2.87,3.65)--(3.37,3.2);
\draw[gray,->] (3.37,3.2)--(3.2,2.58);
\fill[gray] (2.87,3.65) circle (1.65pt);
\fill[gray] (3.37,3.2) circle (1.65pt);
\fill[gray] (3.2,2.58) circle (1.65pt);
\draw[gray,->,line width=1.2pt] (3.2,2.58)--(3.2,.65);
\draw[gray,->,line width=1.2pt] (.85,.65)--(3.2,.65);
\fill[cyan] (.85,.65) circle (2pt);
\fill[cyan] (3.2,.65) circle (2pt);
\node[below] at (.85,.48) {$m-1$};
\node[below] at (3.2,.48) {$m$};
\end{tikzpicture}
\caption{$\Gamma_0$}
\end{subfigure}
\hfill
\begin{subfigure}[t]{.31\textwidth}
\centering
\begin{tikzpicture}[x=.94cm,y=.78cm,font=\small,>=stealth]
\path[use as bounding box] (-.25,-.3) rectangle (4.4,4.65);
\fill[gray!10] (0.85,3.25) ellipse (.68 and .78);
\node[text=gray!80!black] at (0.85,4.2) {$T_1$};
\draw[gray,->] (0.52,3.65)--(1.02,3.2);
\draw[gray,->] (1.02,3.2)--(0.85,2.58);
\fill[gray] (0.52,3.65) circle (1.65pt);
\fill[gray] (1.02,3.2) circle (1.65pt);
\fill[gray] (0.85,2.58) circle (1.65pt);
\draw[orange,->,line width=1.2pt] (0.85,2.58)--(0.85,.65);
\fill[gray!10] (3.2,3.25) ellipse (.68 and .78);
\node[text=gray!80!black] at (3.2,4.2) {$T_2$};
\draw[gray,->] (2.87,3.65)--(3.37,3.2);
\draw[gray,->] (3.37,3.2)--(3.2,2.58);
\fill[gray] (2.87,3.65) circle (1.65pt);
\fill[gray] (3.37,3.2) circle (1.65pt);
\fill[gray] (3.2,2.58) circle (1.65pt);
\draw[gray,->,line width=1.2pt] (3.2,2.58)--(3.2,.65);
\draw[gray,->,line width=1.2pt] (.85,.65)--(3.2,.65);
\fill[cyan] (.85,.65) circle (2pt);
\fill[cyan] (3.2,.65) circle (2pt);
\node[below] at (.85,.48) {$m-1$};
\node[below] at (3.2,.48) {$m$};
\end{tikzpicture}
\caption{$\Gamma_1$}
\end{subfigure}
\hfill
\begin{subfigure}[t]{.31\textwidth}
\centering
\begin{tikzpicture}[x=.94cm,y=.78cm,font=\small,>=stealth]
\path[use as bounding box] (-.25,-.3) rectangle (4.4,4.65);
\fill[gray!10] (3.2,3.25) ellipse (.68 and .78);
\node[text=gray!80!black] at (3.2,4.2) {$T_1$};
\draw[gray,->] (2.87,3.65)--(3.37,3.2);
\draw[gray,->] (3.37,3.2)--(3.2,2.58);
\fill[gray] (2.87,3.65) circle (1.65pt);
\fill[gray] (3.37,3.2) circle (1.65pt);
\fill[gray] (3.2,2.58) circle (1.65pt);
\draw[gray,->,line width=1.2pt] (3.2,2.58)--(3.2,.65);
\fill[gray!10] (0.85,3.25) ellipse (.68 and .78);
\node[text=gray!80!black] at (0.85,4.2) {$T_2$};
\draw[gray,->] (0.52,3.65)--(1.02,3.2);
\draw[gray,->] (1.02,3.2)--(0.85,2.58);
\fill[gray] (0.52,3.65) circle (1.65pt);
\fill[gray] (1.02,3.2) circle (1.65pt);
\fill[gray] (0.85,2.58) circle (1.65pt);
\draw[orange,->,line width=1.2pt] (0.85,2.58)--(0.85,.65);
\draw[gray,->,line width=1.2pt] (.85,.65)--(3.2,.65);
\fill[cyan] (.85,.65) circle (2pt);
\fill[cyan] (3.2,.65) circle (2pt);
\node[below] at (.85,.48) {$m-1$};
\node[below] at (3.2,.48) {$m$};
\end{tikzpicture}
\caption{$\Gamma_2$}
\end{subfigure}
\caption{The third case of
Definition~\ref{Definition_group_of_graphs} for $n=2$.
The shaded trees $T_1,T_2$ are unchanged in all three
graphs. Starting from $\Gamma_0$, reattaching $T_i$
to $m-1$ instead of $m$ gives $\Gamma_i$.
The orange edge marks the changed attachment.}
\label{Figure_tree_group}
\end{figure}

\begin{proposition}
\label{Proposition_tree_zero_area}
For every $m\ge 2$, $\vec\bk\in\R_+^m$,
$\ttt_1\leq\cdots\leq\ttt_m$, and every connected graph $\Gamma$ on  $\llbracket 1,m\rrbracket$  with no $m\to m$ loop, we have $\bL^\Gamma_{[0]}(\vec\bk, \vec\ttt)=0$.
\end{proposition}
\begin{proposition}
\label{Proposition_unicycle_zero_area} For every $m\ge 2$, $\vec\bk\in\R_+^m$,
$\ttt_1\leq\cdots\leq\ttt_m$, and every connected graph $\Gamma$ on  $\llbracket 1,m\rrbracket$ with $m$ edges  we have $\bL^\Gamma_{[0]}(\vec\bk, \vec\ttt)=0$.
\end{proposition}
\begin{definition}
\label{Definition_group_of_graphs}
 Fix $m\ge 3$. A \emph{group of trees} $G$ is a collection of connected graphs on $\llbracket 1,m\rrbracket$ with $m-1$ edges each and such that either:
 \begin{itemize}
  \item $|G|=1$ and the single tree in $G$ has a vertex $\ell<m$, with no edges connecting it to $\llbracket \ell+1,m\rrbracket$.
  \item $|G|=1$ and the single tree in $G$ has at least two edges connecting $m-1$ with $\llbracket 1,m-2\rrbracket$.
  \item $|G|=n+1>1$ and trees $\Gamma_0,\dots,\Gamma_n\in G$ are such that all of them have $m-1\to m$ edge. $\Gamma_0$ has no edges connecting $m-1$ with $\llbracket 1,m-2\rrbracket$, but has $n$ edges connecting $m$ with $\llbracket 1,m-2\rrbracket$. Graphs $\Gamma_1,\dots,\Gamma_n$ differ from $\Gamma_0$ by changing a single edge $i\to m$ into $i\to (m-1)$, $i\le m-2$.
 \end{itemize}
\end{definition}
\begin{remark}
 The groups allowed by this definition do not have to be disjoint. However, the set of all connected graphs on $\llbracket 1,m\rrbracket$ with $m-1$ edges can be partitioned into a disjoint union of subsets, such that each one is a group of trees.
\end{remark}
Figure \ref{Figure_tree_group} illustrates the third case in Definition \ref{Definition_group_of_graphs}.

\begin{proposition}
 \label{Proposition_tree_one_area} For every $m\ge 3$, $\vec\bk\in\R_+^m$,
$\ttt_1\leq\cdots\leq\ttt_m$, and every group of trees $G$, we have $\sum_{\Gamma\in G} \bL^\Gamma_{[1]}(\vec\bk, \vec\ttt)=0$.
\end{proposition}

Note that in the case $e=m-1$, Proposition \ref{Proposition_tree_zero_area} applies. In the rest of this section we prove these propositions.

\subsection{Proof of Proposition \ref{Proposition_tree_zero_area}}
\label{Section_tree_constant}

The proof  is based on a cancellation between the terms with $\bb_m=\varnothing$ and terms with $\bb_m>0$.

\begin{lemma} \label{Lemma_no_loop_vanishing}
 Take blocks $(\vec\bp, \vec\bb, \bH)$ with a connected graph $\Gamma(\vec \bp)$ on $\llbracket 1, m \rrbracket$, which has no loop $m\to m$.
 For almost every such blocks, if the last virtual block is empty, $\bb_m=\varnothing$, then $\bp_m$ is identically zero on the entire segment $[0,\bQ_m]$.
\end{lemma}
\begin{proof}
We argue by contradiction. If $\bp_m$ is not identically zero, then $\bp_m(\bQ_{m-1}-)=h>0$ (the situation when $\bp_m$ jumps to $0$ before time $\bQ_{m-1}$ has measure $0$ in Definition \ref{defn:measurecb} and is therefore ignored) and $\bp_m(\bQ_{m-1})=0$. In this situation, let us take a look at the first jump from $\bDel$ in the interval $(\bQ_{m-1},\bQ_m)$. If there is no such jump, then the third and fifth points in Definition \ref{Definition_block_height} imply $0=\bH(\bQ_m)\ge\bH(\bQ_{m-1})\ge h>0$, which is impossible. Hence, such a jump $x$ exists and corresponds to some function $\bp_j$, $j\ge 1$. We claim that this is a jump to $0$, i.e.\ $\bp_j(x)=0$. Otherwise, $\bp_j$ would have another jump in the interval $(\bQ_{m-1},\bQ_m)$, which would create an edge $m\to m$ that is prohibited. Hence, using the fifth point in Definition \ref{Definition_block_height} and the remark after it, $\bp^0(x-)=\bH(x)\ge \bH(\bQ_{m-1}) > \bp^0(Q_{m-1})$, which is impossible because there is no jump in $(\bQ_{m-1},x)$.
\end{proof}

\begin{proof}[Proof of Proposition \ref{Proposition_tree_zero_area}]

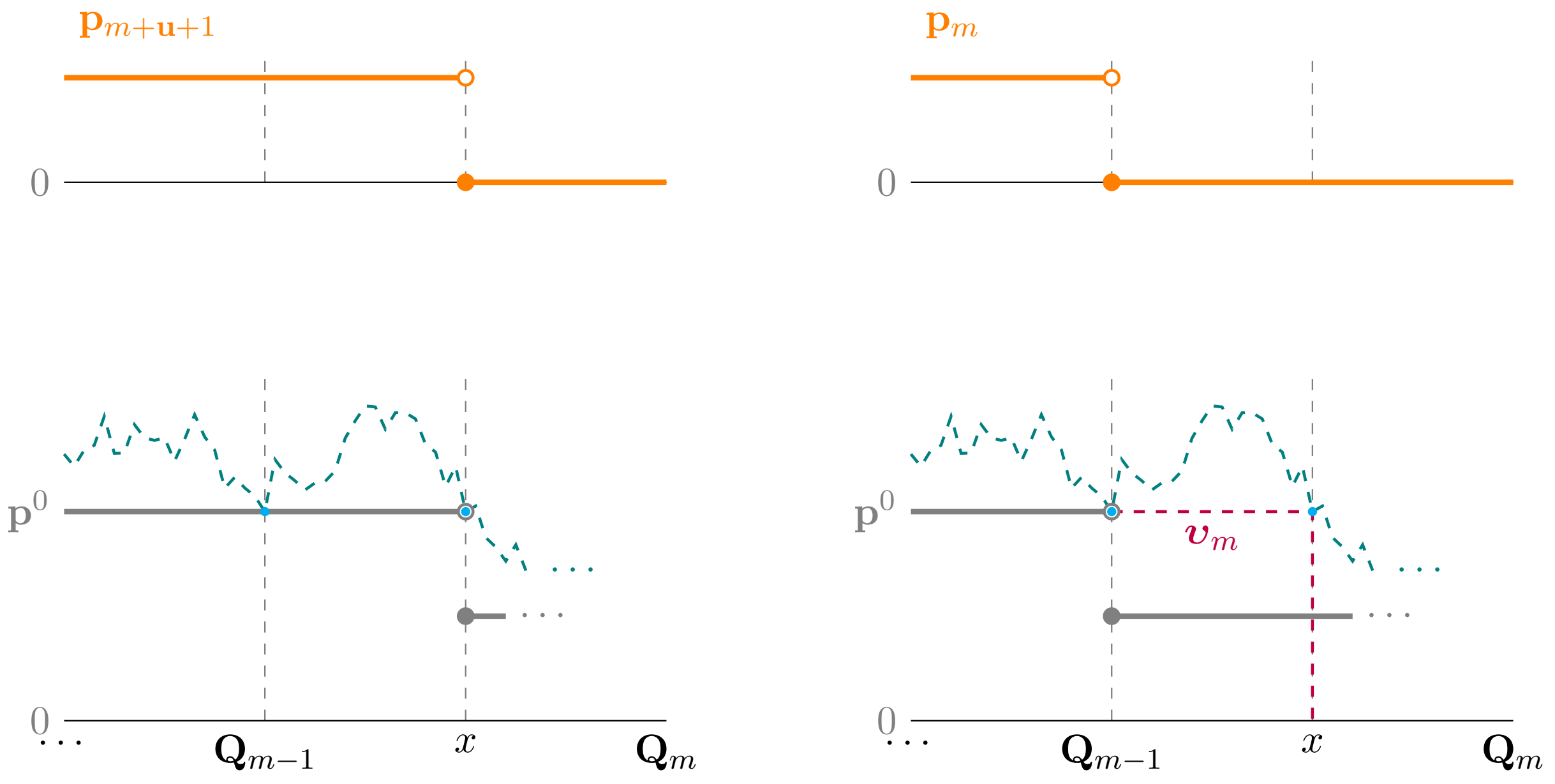
\begin{figure}[t]
  \centering
  \begin{subfigure}[b]{0.46\textwidth}
    \centering
    \scalebox{1.084}{%
    \begin{tikzpicture}[x=0.48cm,y=1cm]
      \path[use as bounding box] (-.8,-.48) rectangle (12.8,6.85);
      \draw[gray,dashed,thin] (4,0)--(4,3.35);
      \draw[gray,dashed,thin] (8,0)--(8,3.35);

      \begin{scope}[yshift=0.65cm]
      \draw[gray,dashed,thin] (4,4.5)--(4,5.75);
      \draw[gray,dashed,thin] (8,4.5)--(8,5.75);
      \draw (0,4.5)--(12,4.5);
      \node[left,text=gray] at (0,4.5) {$0$};
      \node[anchor=west,text=orange] at (0,6.0) {$\bp_{m+\bu+1}$};
      \draw[orange,line width=1.5pt] (0,5.5)--(8,5.5);
      \draw[orange,line width=1.5pt] (8,4.5)--(12,4.5);
      \draw[orange,fill=white,thick] (8,5.5) circle (2pt);
      \draw[orange,fill=orange,thick] (8,4.5) circle (2pt);
      \end{scope}

      \draw (0,0)--(12,0);
      \node[left,text=gray] at (0,0) {$0$};
      \draw[gray,line width=1.5pt] (0,2)--(8,2);
      \draw[gray,line width=1.5pt] (8,1)--(8.8,1);
      \node[text=gray] at (9.6,1) {$\cdots$};
      \draw[teal,thick,dashed] plot coordinates {
      (0,2.55) (0.2,2.433) (0.4,2.588) (0.6,2.634) (0.8,2.911) (1,2.558)
      (1.2,2.561) (1.4,2.833) (1.6,2.707) (1.8,2.681) (2,2.707) (2.2,2.482)
      (2.4,2.691) (2.6,2.928) (2.8,2.715) (3,2.589) (3.2,2.224) (3.4,2.329)
      (3.6,2.226) (3.8,2.147) (4,2) (4.2,2.499) (4.4,2.371) (4.6,2.295)
      (4.8,2.207) (5,2.275) (5.2,2.289) (5.4,2.39) (5.6,2.702) (5.8,2.869)
      (6,3.011) (6.2,3.001) (6.4,2.788) (6.6,2.947) (6.8,2.946) (7,2.887)
      (7.2,2.646) (7.4,2.569) (7.6,2.241) (7.8,2.433) (8,2) (8.2,2.054)
      (8.4,1.752) (8.6,1.661) (8.8,1.526) (9,1.683) (9.2,1.436)
      };
      \node[text=teal] at (10.2,1.436) {$\cdots$};
      \draw[gray,fill=white,thick] (8,2) circle (2pt);
      \draw[gray,fill=gray,thick] (8,1) circle (2pt);

      \foreach \xx in {4,8}
        \fill[cyan] (\xx,2) circle (1.25pt);
      \node[left,text=gray] at (0,2) {$\bp^0$};

      \node[below] at (0,0) {$\cdots$};
      \node[below] at (4,0) {$\bQ_{m-1}$};
      \node[below] at (8,0) {$x$};
      \node[below] at (12,0) {$\bQ_m$};
    \end{tikzpicture}%
    }
  \end{subfigure}
  \hfill
  \begin{subfigure}[b]{0.46\textwidth}
    \centering
    \scalebox{1.084}{%
    \begin{tikzpicture}[x=0.48cm,y=1cm]
      \path[use as bounding box] (-.8,-.48) rectangle (12.8,6.85);
      \draw[gray,dashed,thin] (4,0)--(4,3.35);
      \draw[gray,dashed,thin] (8,0)--(8,3.35);

      \begin{scope}[yshift=0.65cm]
      \draw[gray,dashed,thin] (4,4.5)--(4,5.75);
      \draw[gray,dashed,thin] (8,4.5)--(8,5.75);
      \draw (0,4.5)--(12,4.5);
      \node[left,text=gray] at (0,4.5) {$0$};
      \node[anchor=west,text=orange] at (0,6.0) {$\bp_m$};
      \draw[orange,line width=1.5pt] (0,5.5)--(4,5.5);
      \draw[orange,line width=1.5pt] (4,4.5)--(12,4.5);
      \draw[orange,fill=white,thick] (4,5.5) circle (2pt);
      \draw[orange,fill=orange,thick] (4,4.5) circle (2pt);
      \end{scope}

      \draw (0,0)--(12,0);
      \node[left,text=gray] at (0,0) {$0$};
      \draw[gray,line width=1.5pt] (0,2)--(4,2);
      \draw[gray,line width=1.5pt] (4,1)--(8.8,1);
      \node[text=gray] at (9.6,1) {$\cdots$};
      \draw[purple,thick,dashed] (4,2)--(8,2)--(8,0);
      \node[below,text=purple] at (6,2) {$\bb_m$};
      \draw[teal,thick,dashed] plot coordinates {
      (0,2.55) (0.2,2.433) (0.4,2.588) (0.6,2.634) (0.8,2.911) (1,2.558)
      (1.2,2.561) (1.4,2.833) (1.6,2.707) (1.8,2.681) (2,2.707) (2.2,2.482)
      (2.4,2.691) (2.6,2.928) (2.8,2.715) (3,2.589) (3.2,2.224) (3.4,2.329)
      (3.6,2.226) (3.8,2.147) (4,2) (4.2,2.499) (4.4,2.371) (4.6,2.295)
      (4.8,2.207) (5,2.275) (5.2,2.289) (5.4,2.39) (5.6,2.702) (5.8,2.869)
      (6,3.011) (6.2,3.001) (6.4,2.788) (6.6,2.947) (6.8,2.946) (7,2.887)
      (7.2,2.646) (7.4,2.569) (7.6,2.241) (7.8,2.433) (8,2) (8.2,2.054)
      (8.4,1.752) (8.6,1.661) (8.8,1.526) (9,1.683) (9.2,1.436)
      };
      \node[text=teal] at (10.2,1.436) {$\cdots$};
      \draw[gray,fill=white,thick] (4,2) circle (2pt);
      \draw[gray,fill=gray,thick] (4,1) circle (2pt);

      \foreach \xx in {4,8}
        \fill[cyan] (\xx,2) circle (1.25pt);
      \node[left,text=gray] at (0,2) {$\bp^0$};

      \node[below] at (0,0) {$\cdots$};
      \node[below] at (4,0) {$\bQ_{m-1}$};
      \node[below] at (8,0) {$x$};
      \node[below] at (12,0) {$\bQ_m$};
    \end{tikzpicture}%
    }
  \end{subfigure}%
  \caption{The matching between $\bp_m\equiv 0$, $\bb_m=\varnothing$ (on the left) and $\bp_m\ne 0$, $\bb_m>0$ (on the right) cases. The upper row shows the individual
    block function before and after replacement; the lower row shows $\bp^0$.}
  \label{fig:virtual-block-bijection}
\end{figure}

Using Lemma \ref{Lemma_no_loop_vanishing} and the second point in Definition \ref{Definition_virtual_block}, we split possible blocks $(\vec\bp, \vec\bb, \bH)$ with $\Gamma(\vec \bp)=\Gamma$ into two disjoint situations: either $\bp_m\equiv 0$, $\bb_m=\varnothing$, or $\bp_m\ne 0$, $\bb_m>0$.

We take blocks in the second situation and match them with blocks in the first situation by introducing a new function $\bp_{m+\bu+1}$, such that $\bp_{m+\bu+1}=\bp_m$ on $[0,\bQ_{m-1})$ and
$$
\bp_{m+\bu+1}(x)=\begin{cases}
 \bp_m(\bQ_{m-1}-), & \bQ_{m-1}\le x < \bQ_{m-1}+\bb_m,\\
 0, & \bQ_{m-1}+\bb_m \le x\le \bQ_m.
\end{cases}.
$$
Once $\bp_{m+\bu+1}$ is introduced, we set $\bp_m\equiv 0$, and set $\bb_m=\varnothing$. We would also need to permute the indices of all $\bp_j$ with $m+1\le j \le m+\bu+1$, to ensure the ordering condition in Definition \ref{defn:blp}.
 This procedure is illustrated in Figure \ref{fig:virtual-block-bijection}.

This matching is bijective, because in the $\bp_m\equiv 0$ situation, the new function is uniquely determined as the function $\bp_j$ which makes the first jump in $(\bQ_{m-1},\bQ_m)$. The connectivity of $\Gamma$ implies the existence of such a function; $j>m$, because $\bp_1,\dots,\bp_m$ are all identically zero on $(\bQ_{m-1},\bQ_m)$ by Definition \ref{Definition_block_process}. This first jump in $(\bQ_{m-1},\bQ_m)$ has to be a jump to $0$, as otherwise the graph $\Gamma$ would be forced to have a loop $m\to m$.

Fix $\epsilon>0$.
It remains to compare the contributions of the matched functions to $\bL^\Gamma_{[0];\, (\epsilon)}(\vec\bk, \vec\ttt)$. The power of $2$ in Definition \ref{def:xibi}  is the same, because the change in $\bdel$ is compensated by the change in $|\{\ell\in\llbracket 1,m\rrbracket: \bb_\ell\neq\varnothing\}|$. The sign changes because of the $\Xi_3$ factor for $x$ being the jump of $\bp_j$ in $(\bQ_{m-1},\bQ_m)$: the inequality $\bp^0(x)<\bp^0(x-)$ is no longer valid when $x$ becomes the virtual block endpoint, $x=\bQ_{m-1}+\bb_m$.

The definition of $\bp^0$ is different for the matched functions on the interval $(\bQ_{m-1},\bQ_m)$, and therefore the exponentiated area factor changes in the $\Xi_1$ case of Definition \ref{def:xibi} for the first interval of $(\bQ_{m-1},\bQ_m)$. If we denote $z=\bp_m(\bQ_{m-1}-)$ in the case when $\bp_m$ is not identically zero, then we collect an additional factor $\exp(\beta^{-1} z \bb_m)$ in this situation. Therefore, adding the signed contributions of the matched cases produces the factor $\exp(\beta^{-1} z \bb_m)-1$. It has no constant term and means that the matched cases cancel in $\bL^\Gamma_{[0];\, (\epsilon)}(\vec\bk, \vec\ttt)$. Hence, $\bL^\Gamma_{[0];\, (\epsilon)}(\vec\bk, \vec\ttt)=0$ for each $0<\epsilon<\bk_m$. It remains to send $\epsilon\to0+$. \end{proof}

Our proof of Proposition \ref{Proposition_tree_zero_area} equivalently says that the contribution of $[\bQ_{m-1},\bQ_m]$ to the integral $\bL^\Gamma_{\beta;\, (\epsilon)}(\vec\bk, \vec\ttt)$ has a prefactor $\beta^{-1}$, coming from $\exp(\beta^{-1} z \bb_m)-1$ in the argument. We can be precise in computing this contribution, which will be helpful in later arguments. Note that since $\Gamma$ has no $m\to m$ loop, the only possible jumps in this segment are jumps of some $\bp_j$ from a positive value to $0$. Hence, the part of the integral corresponding to $[\bQ_{m-1},\bQ_m]$ simply integrates over positions of these jumps and over the $\bb_m$.

\begin{lemma} \label{Lemma_first_matching_integral}
 For $n\ge 1$, let $0<z_1<z_2<\dots<z_n$
be all non-zero values\footnote{The situation when two block functions $\bp_j$ have the same positive value has measure zero in the integral \eqref{eq_L_graph_cutoff_r} and therefore can be ignored.} of $\bigl\{\bp_j(\bQ_{m-1}-)\bigr\}_{j\ge 1}$, and let $\bQ_{m-1}=x_0<x_1<x_2<\dots<x_{n+1}=\bQ_m$ be all the points at which $\bH$ is defined in $[\bQ_{m-1},\bQ_m]$. Then
\begin{equation}
\label{eq_first_matching_result}
 \sum_{\text{choices in } [\bQ_{m-1},\bQ_m]}\int 2^{-n}\prod_{i=1}^{n+1} \bI_\beta[x_{i-1},x_i] \dd x_1\cdots\dd x_n= (-1)^{n-1}(n-1)! \beta^{-1} J_{\epsilon}\left(\sum_{i=1}^n z_i\right) + O(\beta^{-2}),
\end{equation}
where the integration goes over ${\bQ_{m-1}+\epsilon<x_1<\dots<x_{n}<\bQ_m}$ and summation is over all combinatorial choices allowed for a block structure: whether $\bb_{m}=\varnothing$ or $\bb_{m}=x_1-\bQ_{m-1}$, and correspondence between $z_i$ and $x_j$. Finally,
\begin{equation} \label{eq_x12}
 J_{\epsilon}\left(Z\right)=\frac{Z}{2} \int_{\epsilon}^{\bk_m} x\bF_{0,0}(x) \bF_0(\bk_m-x; Z)\dd x.
\end{equation}
\end{lemma}
\begin{proof}
 If $\bb_m\ne\emptyset$, then necessarily $x_1=\bQ_{m-1}+\bb_m$. Let $z_i=\bp_m(\bQ_{m-1}-)$ in this situation. All other $z_{i'}$, $i'\ne i$ correspond to some block functions, which stay constant until time $x_{j}$, $j>1$, and then jump to $0$. The correspondence between $i'$ and $j$ indices is parameterized by permutations, leading to $(n-1)!$ cases and corresponding prefactor in \eqref{eq_first_matching_result}. Next, the product in \eqref{eq_first_matching_result} has one factor of type $\Xi_1$ and $n$ factors of type $\Xi_3$. Each of the latter contributes a factor $(-1)$, except for the situation when the left end-point comes from the virtual block, leading to $(-1)^{n-1}$ in \eqref{eq_first_matching_result}. When we match two situations as in the proof of Proposition \ref{Proposition_tree_zero_area}, we get a factor $\exp(\beta^{-1} z_i x)-1$ because of the area change, where $x=\bb_m=x_1-\bQ_{m-1}$. This factor is $\beta^{-1} z_i x +O(\beta^{-2})$. The $x$ here is precisely the $x$ factor in \eqref{eq_x12}. After extracting the linear in $\beta^{-1}$ factor, we replace the remaining Brownian kernels $\bI_{0,0}$ and $\bI_0$ by their constant terms $\bF_{0,0}$ and $\bF_0$, respectively. Therefore, we can use Lemma~\ref{lem:FP_convolution} $(n-1)$ times to integrate out the variables $x_2,\dots,x_{n}$. We obtain a factor $2^{n-1}$, which combines with $2^{-n}$ in \eqref{eq_first_matching_result} to give $\frac{1}{2}$ in \eqref{eq_x12}. The remaining integral over $x_1$ then gives
 $$
  (-1)^{n-1}(n-1)! \beta^{-1} \frac{z_i}{2} \int_{\epsilon}^{\bk_m} x\bF_{0,0}(x) \bF_0(\bk_m-x; \sum_{a=1}^n z_a)\dd x + O(\beta^{-2}).
 $$
 It remains to sum this expression over $i=1,2,\dots,n$.
\end{proof}

\subsection{Proof of Proposition \ref{Proposition_unicycle_zero_area}}
\label{Section_unicycle_zero_area}

If a graph $\Gamma$ has no loop $m\to m$, then $\bL^\Gamma_{[0]}(\vec\bk, \vec\ttt)=0$ by the results of the previous section, and therefore we assume throughout the present argument that such a loop exists.

\begin{lemma}
\label{Lemma_m0_loop_condition}
 Choose a connected graph $\Gamma$ with $m$ edges and a loop $m\to m$. The total contribution to  $\bL^\Gamma_{[0];\, (\epsilon)}(\vec\bk, \vec\ttt)$ of all blocks $(\vec\bp, \vec\bb, \bH)$ with $\Gamma(\vec{\bp})=\Gamma$, except those satisfying the following conditions, is equal to zero:
 \begin{itemize}
  \item $\bb_m=\varnothing$; and
  \item There is a unique jump $x\in (\bQ_{m-1},\bQ_m)\cap\bDel$, such that it is not a jump to $0$, i.e.\ for some $j>m$, $\bp_j(x-)\ne \bp_j(x)\ne 0$; and
  \item $x$ is the first jump in $(\bQ_{m-1},\bQ_m)$, i.e.\ $(\bQ_{m-1},x)\cap\bDel=\emptyset$.
 \end{itemize}
\end{lemma}
\begin{proof}
Compared to the arguments of Section \ref{Section_tree_constant}, when there is an $m\to m$ loop, Lemma \ref{Lemma_no_loop_vanishing} no longer works and we get a new case: $\bb_m=\varnothing$, but $\bp_m$ is not identically $0$. From the proof of the lemma, in order for that to be possible, the first jump in $(\bQ_{m-1},\bQ_m)\cap\bDel$ should not be a jump to $0$ --- hence this is the desired $x$. The jump-not-to-$0$ in $(\bQ_{m-1},\bQ_m)\cap\bDel$ is unique, as an additional such jump would lead to the second $m\to m$ loop in $\Gamma$, which (together with connectivity) requires $e\ge m+1$. We got the three claimed conditions, and additionally $\bp_m$ is not identically zero.

Outside this new case, we can still exploit the matching and cancellation of Section \ref{Section_tree_constant}, however it is no longer true that the matching is a bijection. It is still an injection, yet the inverse does not produce a preimage when (in the $\bp_m\equiv 0$ situation) the first jump in $(\bQ_{m-1},\bQ_m)$ is not a jump to $0$. But this is again exactly the situation with three claimed conditions, but this time $\bp_m\equiv 0$.
\end{proof}

It remains to study the blocks $(\vec\bp, \vec\bb, \bH)$ satisfying the three conditions of Lemma \ref{Lemma_m0_loop_condition}. The cancellation will be isolated in the contribution of the $(\bQ_{m-1},\bQ_{m})$ segment. Let
\begin{equation}
\label{eq_x5}
0<z_1<z_2<\dots<z_n
\end{equation}
be all non-zero values of $\{\bp_j(\bQ_{m-1}-)\}_{j\ge 1}$. Note that the situation when two block functions $\bp_j$ have the same positive value has measure zero in the integral \eqref{eq_L_graph_cutoff_r} and therefore can be ignored. We are going to show that for each fixed $n$ and fixed $\{z_i\}$, the total contribution of the blocks corresponding to this data is equal to $0$. We fix the block data before $\bQ_{m-1}$, leaving unspecified which incoming function is assigned to $\bp_m$. These earlier
data contribute a common factor, which is irrelevant for the cancellations. The freedom we have is in assigning $\{z_i\}$ to different functions $\bp_j$. There are four cases in the integrand for $\bL^\Gamma_{[0]}(\vec\bk, \vec\ttt)$, as shown in Figure \ref{Figure_four_cases}:

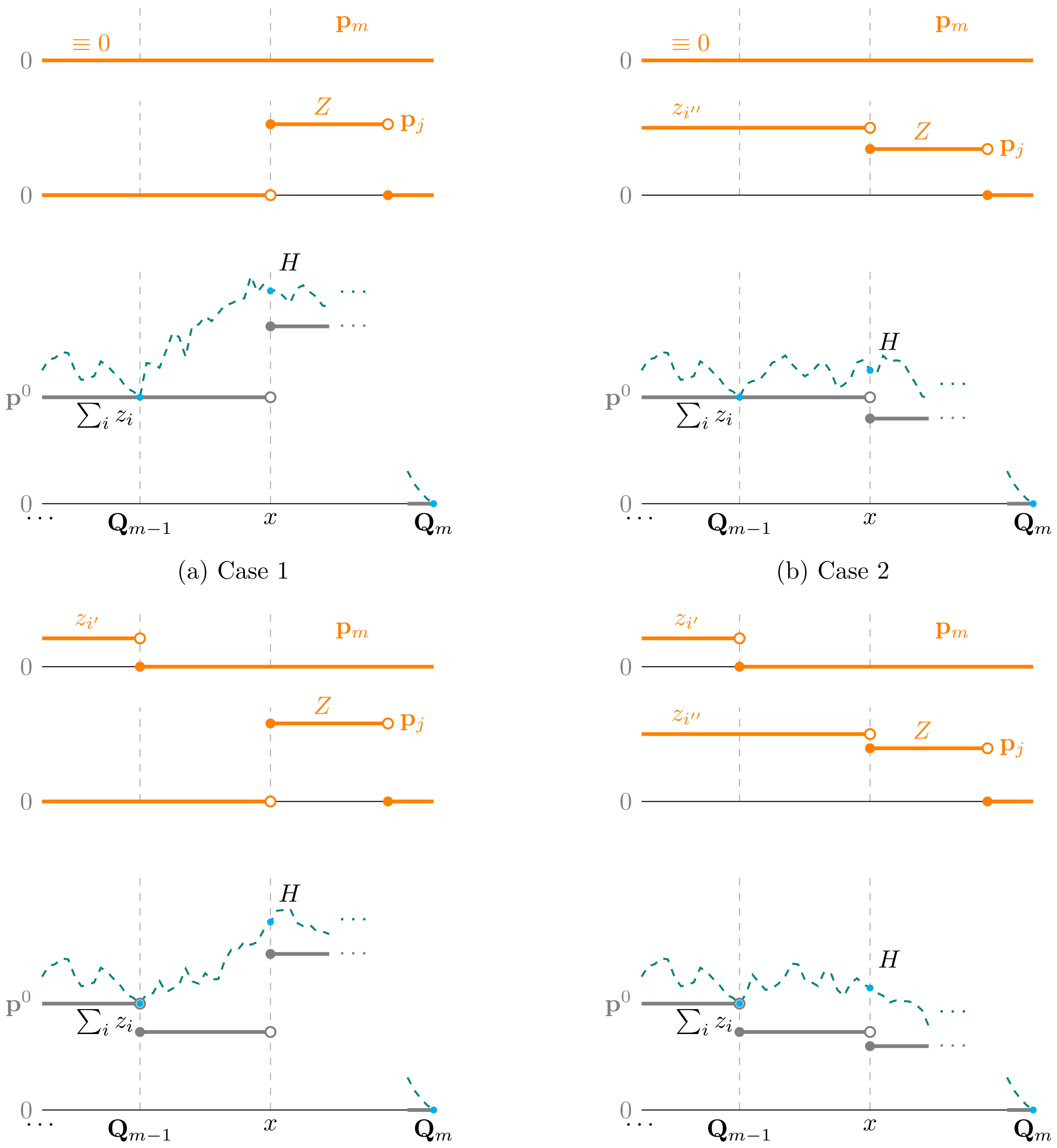
\begin{figure}[t!]
\centering
\begin{subfigure}[t]{0.48\textwidth}
\centering
\begin{tikzpicture}[x=0.46cm,y=1cm,font=\small,
  block/.style={orange,line width=1.5pt},
  barrier/.style={gray,line width=1.5pt},
  brownian/.style={teal,thick,dashed},
  guide/.style={gray!65,dashed,thin}]
\path[use as bounding box] (-1.0,-.55) rectangle (12.8,7.12);
\foreach \xx in {3,7} {
  \draw[guide] (\xx,0)--(\xx,3.35);
  \draw[guide] (\xx,4.35)--(\xx,5.68);
  \draw[guide] (\xx,6.25)--(\xx,7.02);}
\foreach \yy in {0,4.35,6.25} {
  \draw[thin] (0,\yy)--(12,\yy);
  \node[left,text=gray] at (0,\yy) {$0$};}
\node[anchor=west,text=orange] at (8.7,6.76) {$\bp_m$};
\node[anchor=west,text=orange,inner sep=2pt] at (10.8,5.35) {$\bp_j$};
\draw[block] (0,6.25)--(12,6.25);
\node[above,text=orange] at (1.5,6.25) {$\equiv0$};
\draw[block] (0,4.35)--(7,4.35);
\draw[block] (7,5.35)--(10.6,5.35);
\draw[block] (10.6,4.35)--(12,4.35);
\draw[orange,fill=white,thick] (7,4.35) circle (2pt);
\fill[orange] (7,5.35) circle (2pt);
\draw[orange,fill=white,thick] (10.6,5.35) circle (2pt);
\fill[orange] (10.6,4.35) circle (2pt);
\node[above,text=orange] at (8.6,5.35) {$Z$};
\draw[barrier] (0,1.5)--(3,1.5);
\draw[barrier] (3,1.50)--(7,1.50);
\draw[barrier] (7,2.50)--(8.8,2.50);
\node[text=gray] at (9.6,2.50) {$\cdots$};
\draw[barrier] (11.2,0)--(12,0);
\node[left,text=gray] at (0,1.5) {$\bp^0$};
\draw[brownian] plot coordinates {
(0.00,1.880) (0.20,2.028) (0.40,2.056) (0.60,2.136) (0.80,2.125) (1.00,1.886)
        (1.20,1.746) (1.40,1.760) (1.60,1.801) (1.80,2.006) (2.00,1.944) (2.20,1.846)
        (2.40,1.754) (2.60,1.625) (2.80,1.581) (3.00,1.500) (3.20,1.982) (3.40,1.971)
        (3.60,1.917) (3.80,2.154) (4.00,2.408) (4.20,2.350) (4.40,2.070) (4.60,2.507)
        (4.80,2.527) (5.00,2.640) (5.20,2.574) (5.40,2.685) (5.60,2.798) (5.80,2.824)
        (6.00,2.869) (6.20,2.896) (6.40,3.208) (6.60,2.983) (6.80,3.101) (7.00,3.000)
        (7.20,3.010) (7.40,2.922) (7.60,2.829) (7.80,3.032) (8.00,3.080) (8.20,2.988)
        (8.40,2.916) (8.60,2.792) (8.80,2.780)
};
\node[text=teal] at (9.6,2.98) {$\cdots$};
\draw[brownian] plot coordinates {
(11.20,0.460) (11.40,0.284) (11.60,0.165) (11.80,0.068) (12.00,0.000)
};
\draw[gray,fill=white,thick] (7,1.50) circle (2pt);
\fill[gray] (7,2.50) circle (2pt);
\fill[cyan] (3,1.5) circle (1.4pt);
\fill[cyan] (7,3.0) circle (1.4pt);
\fill[cyan] (12,0) circle (1.4pt);
\node[below left,inner sep=2pt] at (3,1.5) {$\sum_i z_i$};
\node[above right,inner sep=2pt,xshift=1pt,yshift=6pt] at (7,3.0) {$H$};
\node[below] at (0,0) {$\cdots$};
\node[below] at (3,0) {$\bQ_{m-1}$};
\node[below] at (7,0) {$x$};
\node[below] at (12,0) {$\bQ_m$};
\end{tikzpicture}
\caption{Case 1}
\end{subfigure}
\hfill
\begin{subfigure}[t]{0.48\textwidth}
\centering
\begin{tikzpicture}[x=0.46cm,y=1cm,font=\small,
  block/.style={orange,line width=1.5pt},
  barrier/.style={gray,line width=1.5pt},
  brownian/.style={teal,thick,dashed},
  guide/.style={gray!65,dashed,thin}]
\path[use as bounding box] (-1.0,-.55) rectangle (12.8,7.12);
\foreach \xx in {3,7} {
  \draw[guide] (\xx,0)--(\xx,3.35);
  \draw[guide] (\xx,4.35)--(\xx,5.68);
  \draw[guide] (\xx,6.25)--(\xx,7.02);}
\foreach \yy in {0,4.35,6.25} {
  \draw[thin] (0,\yy)--(12,\yy);
  \node[left,text=gray] at (0,\yy) {$0$};}
\node[anchor=west,text=orange] at (8.7,6.76) {$\bp_m$};
\node[anchor=west,text=orange,inner sep=2pt] at (10.8,5.00) {$\bp_j$};
\draw[block] (0,6.25)--(12,6.25);
\node[above,text=orange] at (1.5,6.25) {$\equiv0$};
\draw[block] (0,5.30)--(7,5.30);
\draw[block] (7,5.00)--(10.6,5.00);
\draw[block] (10.6,4.35)--(12,4.35);
\draw[orange,fill=white,thick] (7,5.30) circle (2pt);
\fill[orange] (7,5.00) circle (2pt);
\draw[orange,fill=white,thick] (10.6,5.00) circle (2pt);
\fill[orange] (10.6,4.35) circle (2pt);
\node[above,text=orange] at (8.6,5.00) {$Z$};
\node[above,text=orange] at (1.4,5.30) {$z_{i''}$};
\draw[barrier] (0,1.5)--(3,1.5);
\draw[barrier] (3,1.50)--(7,1.50);
\draw[barrier] (7,1.20)--(8.8,1.20);
\node[text=gray] at (9.6,1.20) {$\cdots$};
\draw[barrier] (11.2,0)--(12,0);
\node[left,text=gray] at (0,1.5) {$\bp^0$};
\draw[brownian] plot coordinates {
(0.00,1.880) (0.20,2.028) (0.40,2.056) (0.60,2.136) (0.80,2.125) (1.00,1.886)
        (1.20,1.746) (1.40,1.760) (1.60,1.801) (1.80,2.006) (2.00,1.944) (2.20,1.846)
        (2.40,1.754) (2.60,1.625) (2.80,1.581) (3.00,1.500) (3.20,1.686) (3.40,1.729)
        (3.60,1.737) (3.80,1.845) (4.00,1.988) (4.20,2.019) (4.40,2.090) (4.60,1.970)
        (4.80,1.876) (5.00,1.795) (5.20,1.851) (5.40,1.957) (5.60,1.989) (5.80,1.860)
        (6.00,1.628) (6.20,1.677) (6.40,1.761) (6.60,1.991) (6.80,2.028) (7.00,1.880)
        (7.20,1.825) (7.40,2.086) (7.60,2.000) (7.80,2.022) (8.00,2.004) (8.20,1.842)
        (8.40,1.691) (8.60,1.515) (8.80,1.480)
};
\node[text=teal] at (9.6,1.68) {$\cdots$};
\draw[brownian] plot coordinates {
(11.20,0.460) (11.40,0.284) (11.60,0.165) (11.80,0.068) (12.00,0.000)
};
\draw[gray,fill=white,thick] (7,1.50) circle (2pt);
\fill[gray] (7,1.20) circle (2pt);
\fill[cyan] (3,1.5) circle (1.4pt);
\fill[cyan] (7,1.88) circle (1.4pt);
\fill[cyan] (12,0) circle (1.4pt);
\node[below left,inner sep=2pt] at (3,1.5) {$\sum_i z_i$};
\node[above right,inner sep=2pt,xshift=1pt,yshift=6pt] at (7,1.88) {$H$};
\node[below] at (0,0) {$\cdots$};
\node[below] at (3,0) {$\bQ_{m-1}$};
\node[below] at (7,0) {$x$};
\node[below] at (12,0) {$\bQ_m$};
\end{tikzpicture}
\caption{Case 2}
\end{subfigure}

\medskip

\begin{subfigure}[t]{0.48\textwidth}
\centering
\begin{tikzpicture}[x=0.46cm,y=1cm,font=\small,
  block/.style={orange,line width=1.5pt},
  barrier/.style={gray,line width=1.5pt},
  brownian/.style={teal,thick,dashed},
  guide/.style={gray!65,dashed,thin}]
\path[use as bounding box] (-1.0,-.55) rectangle (12.8,7.12);
\foreach \xx in {3,7} {
  \draw[guide] (\xx,0)--(\xx,3.35);
  \draw[guide] (\xx,4.35)--(\xx,5.68);
  \draw[guide] (\xx,6.25)--(\xx,7.02);}
\foreach \yy in {0,4.35,6.25} {
  \draw[thin] (0,\yy)--(12,\yy);
  \node[left,text=gray] at (0,\yy) {$0$};}
\node[anchor=west,text=orange] at (8.7,6.76) {$\bp_m$};
\node[anchor=west,text=orange,inner sep=2pt] at (10.8,5.45) {$\bp_j$};
\draw[block] (0,6.65)--(3,6.65);
\draw[block] (3,6.25)--(12,6.25);
\draw[orange,fill=white,thick] (3,6.65) circle (2pt);
\fill[orange] (3,6.25) circle (2pt);
\node[above,text=orange] at (1.4,6.65) {$z_{i'}$};
\draw[block] (0,4.35)--(7,4.35);
\draw[block] (7,5.45)--(10.6,5.45);
\draw[block] (10.6,4.35)--(12,4.35);
\draw[orange,fill=white,thick] (7,4.35) circle (2pt);
\fill[orange] (7,5.45) circle (2pt);
\draw[orange,fill=white,thick] (10.6,5.45) circle (2pt);
\fill[orange] (10.6,4.35) circle (2pt);
\node[above,text=orange] at (8.6,5.45) {$Z$};
\draw[barrier] (0,1.5)--(3,1.5);
\draw[barrier] (3,1.10)--(7,1.10);
\draw[barrier] (7,2.20)--(8.8,2.20);
\node[text=gray] at (9.6,2.20) {$\cdots$};
\draw[barrier] (11.2,0)--(12,0);
\node[left,text=gray] at (0,1.5) {$\bp^0$};
\draw[brownian] plot coordinates {
(0.00,1.880) (0.20,2.028) (0.40,2.056) (0.60,2.136) (0.80,2.125) (1.00,1.886)
        (1.20,1.746) (1.40,1.760) (1.60,1.801) (1.80,2.006) (2.00,1.944) (2.20,1.846)
        (2.40,1.754) (2.60,1.625) (2.80,1.581) (3.00,1.500) (3.20,1.633) (3.40,1.646)
        (3.60,1.818) (3.80,1.661) (4.00,1.712) (4.20,1.769) (4.40,1.993) (4.60,1.814)
        (4.80,1.781) (5.00,1.928) (5.20,1.841) (5.40,1.846) (5.60,2.109) (5.80,2.261)
        (6.00,2.258) (6.20,2.373) (6.40,2.329) (6.60,2.368) (6.80,2.546) (7.00,2.650)
        (7.20,2.809) (7.40,2.821) (7.60,2.840) (7.80,2.641) (8.00,2.594) (8.20,2.615)
        (8.40,2.526) (8.60,2.517) (8.80,2.480)
};
\node[text=teal] at (9.6,2.68) {$\cdots$};
\draw[brownian] plot coordinates {
(11.20,0.460) (11.40,0.284) (11.60,0.165) (11.80,0.068) (12.00,0.000)
};
\draw[gray,fill=white,thick] (3,1.5) circle (2pt);
\fill[gray] (3,1.10) circle (2pt);
\draw[gray,fill=white,thick] (7,1.10) circle (2pt);
\fill[gray] (7,2.20) circle (2pt);
\fill[cyan] (3,1.5) circle (1.4pt);
\fill[cyan] (7,2.65) circle (1.4pt);
\fill[cyan] (12,0) circle (1.4pt);
\node[below left,inner sep=2pt] at (3,1.5) {$\sum_i z_i$};
\node[above right,inner sep=2pt,xshift=1pt,yshift=6pt] at (7,2.65) {$H$};
\node[below] at (0,0) {$\cdots$};
\node[below] at (3,0) {$\bQ_{m-1}$};
\node[below] at (7,0) {$x$};
\node[below] at (12,0) {$\bQ_m$};
\end{tikzpicture}
\caption{Case 3}
\end{subfigure}
\hfill
\begin{subfigure}[t]{0.48\textwidth}
\centering
\begin{tikzpicture}[x=0.46cm,y=1cm,font=\small,
  block/.style={orange,line width=1.5pt},
  barrier/.style={gray,line width=1.5pt},
  brownian/.style={teal,thick,dashed},
  guide/.style={gray!65,dashed,thin}]
\path[use as bounding box] (-1.0,-.55) rectangle (12.8,7.12);
\foreach \xx in {3,7} {
  \draw[guide] (\xx,0)--(\xx,3.35);
  \draw[guide] (\xx,4.35)--(\xx,5.68);
  \draw[guide] (\xx,6.25)--(\xx,7.02);}
\foreach \yy in {0,4.35,6.25} {
  \draw[thin] (0,\yy)--(12,\yy);
  \node[left,text=gray] at (0,\yy) {$0$};}
\node[anchor=west,text=orange] at (8.7,6.76) {$\bp_m$};
\node[anchor=west,text=orange,inner sep=2pt] at (10.8,5.10) {$\bp_j$};
\draw[block] (0,6.65)--(3,6.65);
\draw[block] (3,6.25)--(12,6.25);
\draw[orange,fill=white,thick] (3,6.65) circle (2pt);
\fill[orange] (3,6.25) circle (2pt);
\node[above,text=orange] at (1.4,6.65) {$z_{i'}$};
\draw[block] (0,5.30)--(7,5.30);
\draw[block] (7,5.10)--(10.6,5.10);
\draw[block] (10.6,4.35)--(12,4.35);
\draw[orange,fill=white,thick] (7,5.30) circle (2pt);
\fill[orange] (7,5.10) circle (2pt);
\draw[orange,fill=white,thick] (10.6,5.10) circle (2pt);
\fill[orange] (10.6,4.35) circle (2pt);
\node[above,text=orange] at (8.6,5.10) {$Z$};
\node[above,text=orange] at (1.4,5.30) {$z_{i''}$};
\draw[barrier] (0,1.5)--(3,1.5);
\draw[barrier] (3,1.10)--(7,1.10);
\draw[barrier] (7,0.90)--(8.8,0.90);
\node[text=gray] at (9.6,0.90) {$\cdots$};
\draw[barrier] (11.2,0)--(12,0);
\node[left,text=gray] at (0,1.5) {$\bp^0$};
\draw[brownian] plot coordinates {
(0.00,1.880) (0.20,2.028) (0.40,2.056) (0.60,2.136) (0.80,2.125) (1.00,1.886)
        (1.20,1.746) (1.40,1.760) (1.60,1.801) (1.80,2.006) (2.00,1.944) (2.20,1.846)
        (2.40,1.754) (2.60,1.625) (2.80,1.581) (3.00,1.500) (3.20,1.638) (3.40,1.907)
        (3.60,1.796) (3.80,1.693) (4.00,1.715) (4.20,1.784) (4.40,1.888) (4.60,2.065)
        (4.80,2.051) (5.00,1.838) (5.20,1.803) (5.40,1.749) (5.60,1.951) (5.80,1.929)
        (6.00,1.694) (6.20,1.607) (6.40,1.785) (6.60,1.857) (6.80,1.771) (7.00,1.720)
        (7.20,1.656) (7.40,1.610) (7.60,1.517) (7.80,1.546) (8.00,1.536) (8.20,1.530)
        (8.40,1.489) (8.60,1.404) (8.80,1.180)
};
\node[text=teal] at (9.6,1.38) {$\cdots$};
\draw[brownian] plot coordinates {
(11.20,0.460) (11.40,0.284) (11.60,0.165) (11.80,0.068) (12.00,0.000)
};
\draw[gray,fill=white,thick] (3,1.5) circle (2pt);
\fill[gray] (3,1.10) circle (2pt);
\draw[gray,fill=white,thick] (7,1.10) circle (2pt);
\fill[gray] (7,0.90) circle (2pt);
\fill[cyan] (3,1.5) circle (1.4pt);
\fill[cyan] (7,1.72) circle (1.4pt);
\fill[cyan] (12,0) circle (1.4pt);
\node[below left,inner sep=2pt] at (3,1.5) {$\sum_i z_i$};
\node[above right,inner sep=2pt,xshift=1pt,yshift=6pt] at (7,1.72) {$H$};
\node[below] at (0,0) {$\cdots$};
\node[below] at (3,0) {$\bQ_{m-1}$};
\node[below] at (7,0) {$x$};
\node[below] at (12,0) {$\bQ_m$};
\end{tikzpicture}
\caption{Case 4}
\end{subfigure}
\caption{The four cases in Section \ref{Section_unicycle_zero_area}.}
\label{Figure_four_cases}
\end{figure}

\begin{enumerate}
 \item[Case 1:] $\bp_m\equiv 0$. Some block function $\bp_j$ has value zero at $\bQ_{m-1}-$, jumps from $0$ to $Z>0$ at time $x$ and then jumps back to $0$ at some later time. $z_1,\dots,z_n$ correspond to functions $\bp_j$ which all jump from level $z_i$ to $0$ at some times in $(x,\bQ_m)$ interval.
 \item[Case 2:] $\bp_m\equiv 0$. One value $z_{i''}$ corresponds to a block function $\bp_j$ which jumps from level $z_{i''}$ to $Z>0$ at time $x$, and then jumps to $0$ at a later time. Remaining $n-1$ values $z_i$ correspond  to block functions which all jump from level $z_i$ to $0$ at some times in $(x,\bQ_m)$ interval.
 \item[Case 3:] $\bp_m(\bQ_{m-1}-)=z_{i'}$ for some $1\le i' \le n$.  Some block function $\bp_j$ has value zero at $\bQ_{m-1}-$, jumps from $0$ to $Z>0$ at time $x$ and then jumps back to $0$ at some later time. The $n-1$ values  $z_i$, $i\ne i'$ correspond to functions $\bp_j$ which all jump from level $z_i$ to $0$ at some times in $(x,\bQ_m)$ interval.
 \item[Case 4:] $\bp_m(\bQ_{m-1}-)=z_{i'}$ for some $1\le i' \le n$.  One value $z_{i''}$, $i''\ne i'$ corresponds to a block function $\bp_j$ which jumps from level $z_{i''}$ to $Z>0$ at time $x$, and then jumps to $0$ at a later time. Remaining $n-2$ values $z_i$ correspond  to block functions which all jump from level $z_i$ to $0$ at some times in $(x,\bQ_m)$ interval.
\end{enumerate}

The integral itself goes over the following variables: $\bQ_{m-1}<x<\bQ_{m}$, $Z>0$, the value $H$ of the height process at point $x$; the positions of remaining jumps to $0$ in $(x,\bQ_m)$ interval (the heights at these positions are all fixed to be the corresponding values of $\bp^0$).

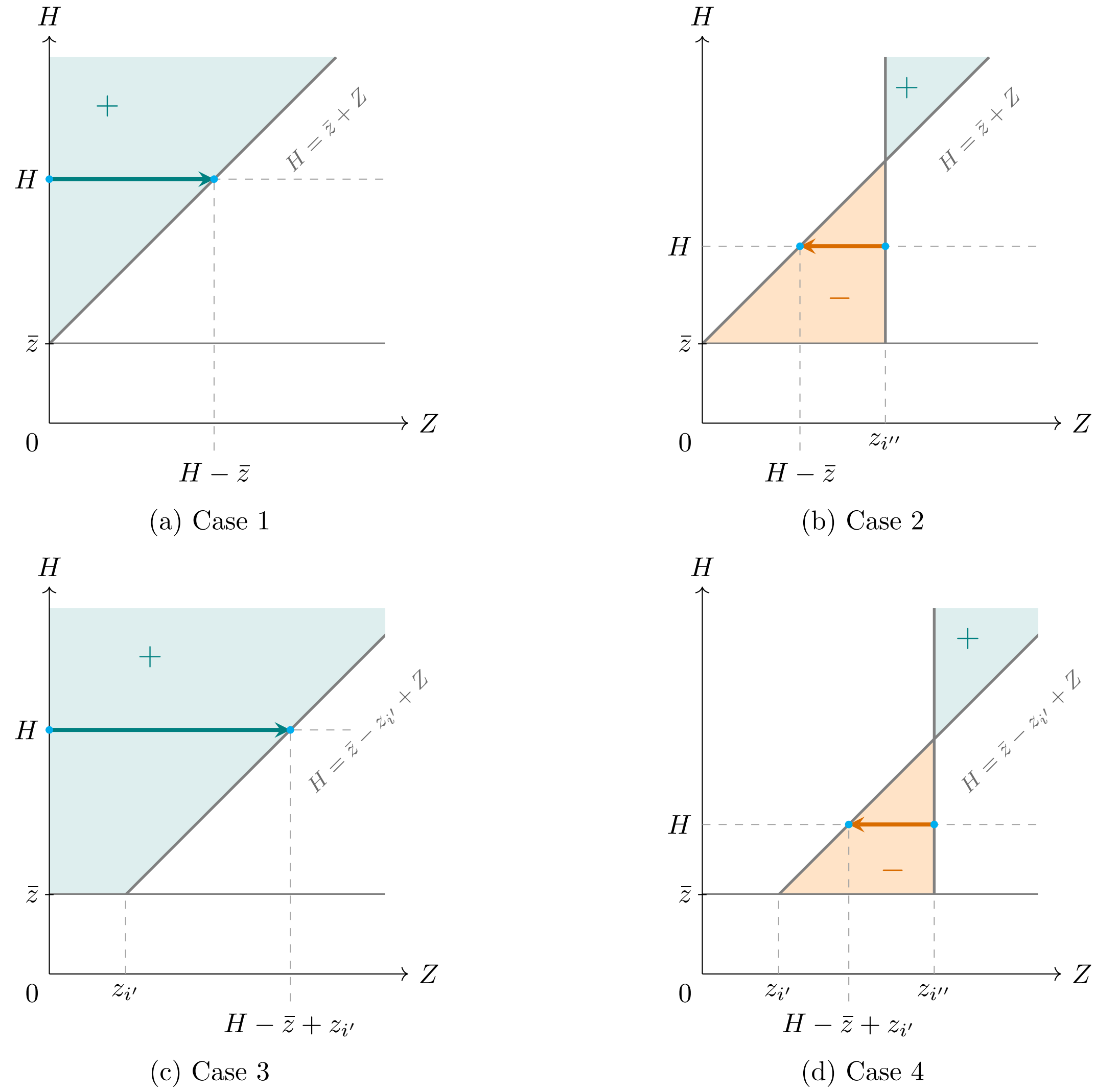
\begin{figure}[t!]
\centering
\begin{subfigure}[t]{0.48\textwidth}
\centering
\begin{tikzpicture}[x=0.79cm,y=0.79cm,font=\small,
  boundary/.style={gray,line width=1pt},
  slice/.style={gray!70,dashed,thin},
  positive/.style={teal,line width=1.5pt,->,>=stealth},
  negative/.style={orange!85!black,line width=1.5pt,->,>=stealth}]
\path[use as bounding box] (-1.3,-1.1) rectangle (6.6,6.55);
\begin{scope}
\clip (0,1.3) rectangle (5.5,6.0);
\fill[teal!13] (0.000,1.300)--(4.700,6.000)--(0.000,6.000)--cycle;
\draw[boundary] (0,1.3)--(4.7,6.0);
\draw[boundary] (0,1.3)--(5.5,1.3);
\end{scope}
\draw[->,thin] (0,0)--(5.88,0) node[right] {$Z$};
\draw[->,thin] (0,0)--(0,6.35) node[above] {$H$};
\node[below left] at (0,0) {$0$};
\draw (-.07,1.3)--(.07,1.3);
\node[left] at (0,1.3) {$\bar z$};
\draw[slice] (0,4.0)--(5.5,4.0);
\node[left] at (0,4.0) {$H$};
\draw[slice] (2.7,-.45)--(2.7,4.0);
\node[below] at (2.7,-.48) {$H-\bar z$};
\draw[positive] (0.0,4.0)--(2.7,4.0);
\fill[cyan] (0.0,4.0) circle (1.4pt);
\fill[cyan] (2.7,4.0) circle (1.4pt);
\node[text=teal,font=\large] at (0.95,5.2) {$+$};
\node[rotate=45,text=gray!85!black,font=\scriptsize,fill=white,inner sep=1pt] at (4.55,4.8) {$H=\bar z+Z$};
\end{tikzpicture}
\caption{Case 1}
\end{subfigure}
\hfill
\begin{subfigure}[t]{0.48\textwidth}
\centering
\begin{tikzpicture}[x=0.79cm,y=0.79cm,font=\small,
  boundary/.style={gray,line width=1pt},
  slice/.style={gray!70,dashed,thin},
  positive/.style={teal,line width=1.5pt,->,>=stealth},
  negative/.style={orange!85!black,line width=1.5pt,->,>=stealth}]
\path[use as bounding box] (-1.3,-1.1) rectangle (6.6,6.55);
\begin{scope}
\clip (0,1.3) rectangle (5.5,6.0);
\fill[teal!13] (3.000,4.300)--(4.700,6.000)--(3.000,6.000)--cycle;
\fill[orange!22] (0.0,1.3)--(3.0,1.3)--(3.0,4.3)--cycle;
\draw[boundary] (0,1.3)--(4.7,6.0);
\draw[boundary] (0,1.3)--(5.5,1.3);
\draw[boundary] (3.0,1.3)--(3.0,6.0);
\end{scope}
\draw[->,thin] (0,0)--(5.88,0) node[right] {$Z$};
\draw[->,thin] (0,0)--(0,6.35) node[above] {$H$};
\node[below left] at (0,0) {$0$};
\draw (-.07,1.3)--(.07,1.3);
\node[left] at (0,1.3) {$\bar z$};
\draw[slice] (0,2.9)--(5.5,2.9);
\node[left] at (0,2.9) {$H$};
\draw[slice] (1.5999999999999999,-.45)--(1.5999999999999999,2.9);
\draw[slice] (3.0,0)--(3.0,1.3);
\node[below] at (3.0,0) {$z_{i''}$};
\node[below] at (1.5999999999999999,-.48) {$H-\bar z$};
\draw[negative] (3.0,2.9)--(1.5999999999999999,2.9);
\fill[cyan] (3.0,2.9) circle (1.4pt);
\fill[cyan] (1.5999999999999999,2.9) circle (1.4pt);
\node[text=teal,font=\large] at (3.35,5.5) {$+$};
\node[text=orange!85!black,font=\large] at (2.25,2.05) {$-$};
\node[rotate=45,text=gray!85!black,font=\scriptsize,fill=white,inner sep=1pt] at (4.55,4.8) {$H=\bar z+Z$};
\end{tikzpicture}
\caption{Case 2}
\end{subfigure}

\bigskip

\begin{subfigure}[t]{0.48\textwidth}
\centering
\begin{tikzpicture}[x=0.79cm,y=0.79cm,font=\small,
  boundary/.style={gray,line width=1pt},
  slice/.style={gray!70,dashed,thin},
  positive/.style={teal,line width=1.5pt,->,>=stealth},
  negative/.style={orange!85!black,line width=1.5pt,->,>=stealth}]
\path[use as bounding box] (-1.3,-1.1) rectangle (6.6,6.55);
\begin{scope}
\clip (0,1.3) rectangle (5.5,6.0);
\fill[teal!13] (0.000,1.300)--(1.250,1.300)--(5.950,6.000)--(0.000,6.000)--cycle;
\draw[boundary] (0,0.050000000000000044)--(5.95,6.0);
\draw[boundary] (0,1.3)--(5.5,1.3);
\end{scope}
\draw[->,thin] (0,0)--(5.88,0) node[right] {$Z$};
\draw[->,thin] (0,0)--(0,6.35) node[above] {$H$};
\node[below left] at (0,0) {$0$};
\draw (-.07,1.3)--(.07,1.3);
\node[left] at (0,1.3) {$\bar z$};
\draw[slice] (0,4.0)--(5.5,4.0);
\node[left] at (0,4.0) {$H$};
\draw[slice] (3.95,-.45)--(3.95,4.0);
\draw[slice] (1.25,0)--(1.25,1.3);
\node[below] at (1.25,0) {$z_{i'}$};
\node[below] at (3.95,-.48) {$H-\bar z+z_{i'}$};
\draw[positive] (0.0,4.0)--(3.95,4.0);
\fill[cyan] (0.0,4.0) circle (1.4pt);
\fill[cyan] (3.95,4.0) circle (1.4pt);
\node[text=teal,font=\large] at (1.65,5.2) {$+$};
\node[rotate=45,text=gray!85!black,font=\scriptsize,fill=white,inner sep=1pt] at (5.25,3.95) {$H=\bar z-z_{i'}+Z$};
\end{tikzpicture}
\caption{Case 3}
\end{subfigure}
\hfill
\begin{subfigure}[t]{0.48\textwidth}
\centering
\begin{tikzpicture}[x=0.79cm,y=0.79cm,font=\small,
  boundary/.style={gray,line width=1pt},
  slice/.style={gray!70,dashed,thin},
  positive/.style={teal,line width=1.5pt,->,>=stealth},
  negative/.style={orange!85!black,line width=1.5pt,->,>=stealth}]
\path[use as bounding box] (-1.3,-1.1) rectangle (6.6,6.55);
\begin{scope}
\clip (0,1.3) rectangle (5.5,6.0);
\fill[teal!13] (3.800,3.850)--(5.950,6.000)--(3.800,6.000)--cycle;
\fill[orange!22] (1.25,1.3)--(3.8,1.3)--(3.8,3.8499999999999996)--cycle;
\draw[boundary] (0,0.050000000000000044)--(5.95,6.0);
\draw[boundary] (0,1.3)--(5.5,1.3);
\draw[boundary] (3.8,1.3)--(3.8,6.0);
\end{scope}
\draw[->,thin] (0,0)--(5.88,0) node[right] {$Z$};
\draw[->,thin] (0,0)--(0,6.35) node[above] {$H$};
\node[below left] at (0,0) {$0$};
\draw (-.07,1.3)--(.07,1.3);
\node[left] at (0,1.3) {$\bar z$};
\draw[slice] (0,2.45)--(5.5,2.45);
\node[left] at (0,2.45) {$H$};
\draw[slice] (2.4000000000000004,-.45)--(2.4000000000000004,2.45);
\draw[slice] (1.25,0)--(1.25,1.3);
\node[below] at (1.25,0) {$z_{i'}$};
\draw[slice] (3.8,0)--(3.8,1.3);
\node[below] at (3.8,0) {$z_{i''}$};
\node[below] at (2.4000000000000004,-.48) {$H-\bar z+z_{i'}$};
\draw[negative] (3.8,2.45)--(2.4000000000000004,2.45);
\fill[cyan] (3.8,2.45) circle (1.4pt);
\fill[cyan] (2.4000000000000004,2.45) circle (1.4pt);
\node[text=teal,font=\large] at (4.35,5.5) {$+$};
\node[text=orange!85!black,font=\large] at (3.12,1.7) {$-$};
\node[rotate=45,text=gray!85!black,font=\scriptsize,fill=white,inner sep=1pt] at (5.25,3.95) {$H=\bar z-z_{i'}+Z$};
\end{tikzpicture}
\caption{Case 4}
\end{subfigure}
\caption{The signed indicator $\sigma$ in the $(Z,H)$ plane.
For $Z>0$ and $H>\bar z$, teal shading marks $\sigma=+1$,
orange shading marks $\sigma=-1$, and unshaded regions have
$\sigma=0$.}
\label{Figure_integration_domains}
\end{figure}

As in the proof of Lemma \ref{Lemma_first_matching_integral}, we will
integrate the times of jumps to $0$ using Lemma~\ref{lem:FP_convolution}.
For the cancellation, it is useful to keep one of these times
fixed. Denote by $y\in(x,\bQ_m)$ the jump to $0$ of the
function whose value after $x$ is $Z$. We introduce two subsets of the indices $\llbracket 1,n\rrbracket$, $\mathcal A$ and $\mathcal B$. An index $i$ belongs to $\mathcal A$ whenever the block function $\bp_j$ which entered at $\bp_j(\bQ_{m-1}-)=z_i$, makes a jump to $0$ after time $y$.  $i$ belongs to $\mathcal B$ whenever the jump to $0$ is in $(x,y)$ interval.
 Also denote:
\begin{equation}\label{eq_uniform_cut_height}
 \bar z=\sum_{i=1}^n z_i,\qquad z^{\mathcal A}=\sum_{i\in \mathcal A}z_i.
\end{equation}
Note that $\bH(y)=\bp^0(y-)=Z+z^{\mathcal A}$. Let us record the function of four variables $(x,y,Z,H)$ which is obtained by integrating the rest; this function is slightly different in the four cases. Define:
\begin{equation}\label{eq_uniform_indicators}
 u_0=\don[Z<H-\bar z],\qquad u_i=\don[Z<H-\bar z+z_i],\qquad v_i=\don[Z<z_i].
\end{equation}
A function $\sigma$ is built out of these three in a case-dependent fashion:
\begin{equation}\label{eq_sigma_four_cases}
\sigma=\begin{cases} u_0, &\text{Case 1},\\
            u_0-v_{i''}, &\text{Case 2},\\
            u_{i'}, &\text{Case 3},\\
            u_{i'}-v_{i''}, &\text{Case 4}
        \end{cases}
\end{equation}
This function takes values $0,\pm 1$ and encodes the integration domain in $(Z,H)$ plane and the sign at $x$ for the integrand. It is shown in Figure \ref{Figure_integration_domains}.

\begin{lemma}\label{Lemma_cut_integral}
Up to sets of measure zero, admissible blocks satisfy
$H>\bar z$, $0<Z<H-z^{\mathcal A}$. Choose one of the four cases, the indices $i'$, $i''$, the sets $\mathcal A$, $\mathcal B$, and two permutations: $\rho$ of $\mathcal A$ and $\pi$ of $\mathcal B$, which encode the order of jumps to $0$ for the block functions in $\mathcal A$ and $\mathcal B$. For the integral $\bL^\Gamma_{[0]}(\vec\bk, \vec\ttt)$, consider the constant term in $\beta^{-1}$ of the factors
in Definition~\ref{def:xibi} belonging to
$[\bQ_{m-1},\bQ_m]$, including a factor $1/2$ for each jump
in this interval. Keeping $x,y,Z,H$ fixed and integrating
over all other jump times gives
\begin{equation}\label{eq_uniform_unicycle_factor}
 (-1)^{|\mathcal A|+|\mathcal B|+1} \frac{\sigma}{4}\bF_0(x-\bQ_{m-1};H-\bar z)
 \bF_0(y-x;H-Z-z^{\mathcal A})\bF_0(\bQ_m-y;Z+z^{\mathcal A}).
\end{equation}
\end{lemma}

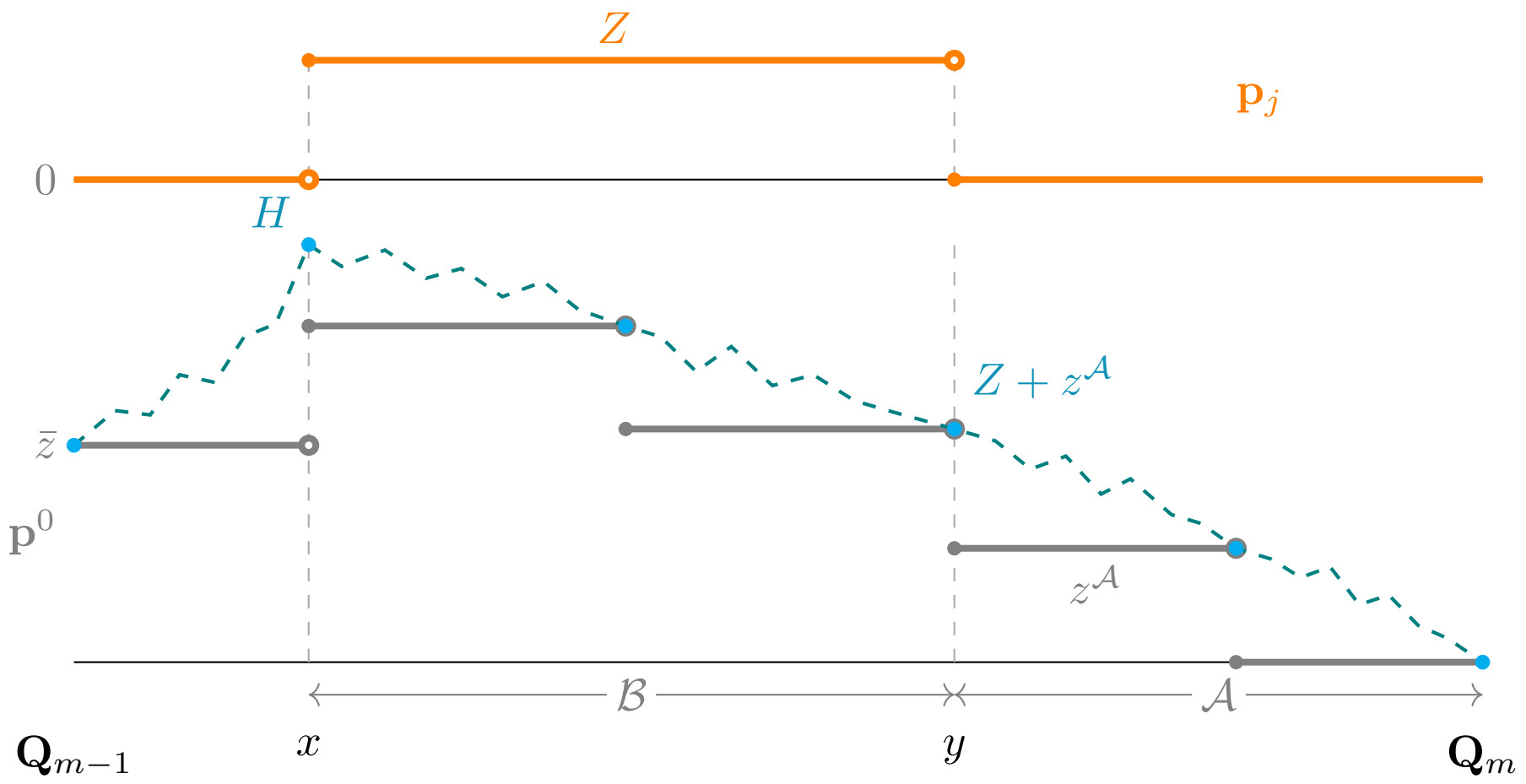
\begin{figure}[t]
\centering
\begin{tikzpicture}[x=.92cm,y=.85cm,font=\small,
 block/.style={orange,line width=1.5pt},
 barrier/.style={gray,line width=1.5pt},
 brownian/.style={teal,thick,dashed},
 guide/.style={gray!65,dashed,thin}]
\path[use as bounding box] (-1.0,-.95) rectangle (13.2,6.05);
\draw[thin] (0,0)--(12,0);
\draw[thin] (0,4.45)--(12,4.45);
\node[left,text=gray] at (0,4.45) {$0$};
\node[left,text=gray] at (0,1.2) {$\bp^0$};
\foreach \xx in {2.0,7.5} {
 \draw[guide] (\xx,0)--(\xx,3.95);
 \draw[guide] (\xx,4.45)--(\xx,5.65);}
\draw[block] (0,4.45)--(2,4.45);
\draw[block,fill=white] (2,4.45) circle (1.6pt);
\fill[orange] (2,5.55) circle (1.6pt);
\draw[block] (2,5.55)--(7.5,5.55);
\draw[block,fill=white] (7.5,5.55) circle (1.6pt);
\fill[orange] (7.5,4.45) circle (1.6pt);
\draw[block] (7.5,4.45)--(12,4.45);
\node[above,text=orange] at (4.6,5.55) {$Z$};
\node[text=orange] at (10.1,5.18) {$\bp_j$};
\draw[barrier] (0,2.0)--(2,2.0);
\draw[barrier,fill=white] (2,2.0) circle (1.6pt);
\fill[gray] (2,3.1) circle (1.6pt);
\draw[barrier] (2,3.1)--(4.7,3.1);
\draw[barrier,fill=white] (4.7,3.1) circle (1.6pt);
\fill[gray] (4.7,2.15) circle (1.6pt);
\draw[barrier] (4.7,2.15)--(7.5,2.15);
\draw[barrier,fill=white] (7.5,2.15) circle (1.6pt);
\fill[gray] (7.5,1.05) circle (1.6pt);
\draw[barrier] (7.5,1.05)--(9.9,1.05);
\draw[barrier,fill=white] (9.9,1.05) circle (1.6pt);
\fill[gray] (9.9,0) circle (1.6pt);
\draw[barrier] (9.9,0)--(12,0);
\draw[brownian] plot coordinates {
 (0,2.0) (.35,2.32) (.65,2.28) (.9,2.65) (1.2,2.58)
 (1.45,3.0) (1.72,3.12) (2,3.85)
 (2.28,3.65) (2.65,3.8) (3,3.54) (3.3,3.63) (3.65,3.37)
 (4,3.51) (4.3,3.25) (4.7,3.1)
 (5,3.0) (5.3,2.68) (5.6,2.91) (5.95,2.55) (6.3,2.65)
 (6.65,2.41) (7,2.3) (7.5,2.15)
 (7.85,2.04) (8.15,1.78) (8.45,1.9) (8.75,1.55)
 (9.0,1.69) (9.35,1.36) (9.6,1.28) (9.9,1.05)
 (10.2,.95) (10.45,.78) (10.7,.88) (10.95,.53)
 (11.2,.62) (11.45,.34) (11.7,.22) (12,0)};
\foreach \p in {(0,2),(2,3.85),(4.7,3.1),(7.5,2.15),(9.9,1.05),(12,0)}
 \fill[cyan] \p circle (1.65pt);
\node[above left,text=cyan!70!black] at (2,3.85) {$H$};
\node[left,text=gray] at (0,2.0) {$\bar z$};
\node[above right,text=cyan!70!black,yshift=3pt] at (7.5,2.15) {$Z+z^{\mathcal A}$};
\node[below,text=gray] at (8.7,1.05) {$z^{\mathcal A}$};
\draw[<->,gray] (2,-.3)--(7.5,-.3)
 node[midway,fill=white,inner sep=2pt] {$\mathcal B$};
\draw[<->,gray] (7.5,-.3)--(12,-.3)
 node[midway,fill=white,inner sep=2pt] {$\mathcal A$};
\node[below] at (0,-.52) {$\bQ_{m-1}$};
\node[below] at (2,-.52) {$x$};
\node[below] at (7.5,-.52) {$y$};
\node[below] at (12,-.52) {$\bQ_m$};
\end{tikzpicture}
\caption{Keeping the jump at $y$ unintegrated, shown in Case 1
with one index in $\mathcal B$ and one in $\mathcal A$. The orange function changes
to $Z$ at $x$ and vanishes at $y$.}
\label{Figure_marked_return}
\end{figure}

\begin{remark}\label{Remark_cut_endpoint}
The proof of Lemma~\ref{Lemma_cut_integral} also computes the
two factors in its stated expression separately.
For fixed orders $\pi$ and $\rho$, the part up to $y$, including
the jump at $x$ and the jumps indexed by $\mathcal B$, is
\begin{equation}\label{eq_cut_left_part}
 \frac{(-1)^{|\mathcal B|}\sigma}{2}
 \bF_0(x-\bQ_{m-1};H-\bar z)
 \bF_0(y-x;H-Z-z^{\mathcal A}).
\end{equation}
The part starting with the jump at $y$, including the later
jumps indexed by $\mathcal A$, is
\[
 \frac{(-1)^{|\mathcal A|+1}}{2}
 \bF_0(\bQ_m-y;Z+z^{\mathcal A}).
\]
There is no jump at $y$ included in \eqref{eq_cut_left_part} and this formula still applies if $y$ is the right endpoint
of the interval under consideration and the functions with values
$Z$ and $z_i$, $i\in\mathcal A$, remain at these values up
to $y-$. The prescribed height there is $Z+z^{\mathcal A}$.
In this situation we use only \eqref{eq_cut_left_part}; the
subsequent interval, including any change at $y$, is treated
separately. In
Section~\ref{Section_tree_one_area}, we  apply this
observation with $\bQ_{m-2}$ in place of $\bQ_{m-1}$ and
$y=\bQ_{m-1}$. Figure~\ref{Figure_marked_return} shows the
two parts when $y$ is a jump time.
\end{remark}

\begin{proof}[Proof of Lemma~\ref{Lemma_cut_integral} and
Remark~\ref{Remark_cut_endpoint}]
We first verify the admissible ranges and the sign $\sigma$.
For this calculation, set $z_{i'}=0$ in Cases 1 and 2, where
there is no incoming value assigned to $\bp_m$, and set
$z_{i''}=0$ in Cases 1 and 3, where the function changing at $x$
starts from $0$. These are conventions for the absent choices,
not additional indices. Since $x$ is the first jump, we have
\[
 \bp^0(x-)=\bar z-z_{i'},\qquad
 \bp^0(x)=\bar z-z_{i'}-z_{i''}+Z.
\]
The fifth point in Definition~\ref{Definition_block_height}
gives $H\ge\bH(\bQ_{m-1})=\bar z$.  By the fourth point in the same definition, $Z-z_{i''}$ and
$H-\bar z+z_{i'}-Z$ must have the same sign. Thus the
possibilities are
\[
 \begin{array}{ll}
 \text{an upward jump:}&z_{i''}<Z<H-\bar z+z_{i'},\\[2pt]
 \text{a downward jump:}&H-\bar z+z_{i'}<Z<z_{i''}.
 \end{array}
\]
The remaining height condition is automatic. For an upward
jump, $H>\bar z-z_{i'}+Z\ge\bp^0(x)$; for a downward
jump, $\bp^0(x)<\bp^0(x-)=\bar z-z_{i'}\le\bar z<H$.
The $\Xi_3$ factor following $x$ has sign $+1$ in the first
case and $-1$ in the second. The signed indicator is therefore
\[
 \don[z_{i''}<Z<H-\bar z+z_{i'}]
 -\don[H-\bar z+z_{i'}<Z<z_{i''}]
 =\don[Z<H-\bar z+z_{i'}]-\don[Z<z_{i''}],
\]
which gives exactly \eqref{eq_sigma_four_cases}.

The indices in $\mathcal A$ and $\mathcal B$, together with
the distinguished indices $i'$ and $i''$ when present, partition
$\llbracket1,n\rrbracket$. Consequently,
\[
 \bp^0(x)=Z+z^{\mathcal A}+\sum_{i\in\mathcal B}z_i.
\]
On the support of $\sigma$ this quantity is smaller than $H$.
In particular, $0<Z<H-z^{\mathcal A}$, as claimed.

We now use the convolution argument from the proof of
Lemma~\ref{Lemma_first_matching_integral}. The first interval,
$(\bQ_{m-1},x)$, is of type $\Xi_2$, and all subsequent
intervals are of type $\Xi_3$. Since we take the constant term
in $\beta^{-1}$, their Brownian kernels are replaced by
$\bF_0$. The first interval contributes
$\bF_0(x-\bQ_{m-1};H-\bar z)$.
In $(x,y)$, the depths in the successive $\bF_0$ factors are
\[
 H-Z-z^{\mathcal A}-\sum_{i\in\mathcal B}z_i,
 \qquad z_i,\quad i\in\mathcal B\text{ in order }\pi.
\]
Their sum is $H-Z-z^{\mathcal A}$. Applying
Lemma~\ref{lem:FP_convolution} $|\mathcal B|$ times gives
$2^{|\mathcal B|}\bF_0(y-x;H-Z-z^{\mathcal A})$.
The jumps at $x$ and at these $|\mathcal B|$ times contribute
$2^{-|\mathcal B|-1}$ and the signed indicator
$(-1)^{|\mathcal B|}\sigma$. This proves
\eqref{eq_cut_left_part}.

After $y$, the first depth is
$\bH(y)-\bp^0(y)=Z$, followed by $z_i$, $i\in\mathcal A$,
in order $\rho$. Their sum is $Z+z^{\mathcal A}$.
The $|\mathcal A|$ convolutions give
$2^{|\mathcal A|}\bF_0(\bQ_m-y;Z+z^{\mathcal A})$.
The jump at $y$ and the $|\mathcal A|$ subsequent jumps each
contribute a factor $1/2$ and a negative sign, giving the second
factor in Remark~\ref{Remark_cut_endpoint}. Multiplying the
two factors proves Lemma~\ref{Lemma_cut_integral}.
No factorial is present because $\pi$ and $\rho$ were fixed.
If either set is empty, the corresponding calculation uses no
convolution and gives the same formula.

The calculation before $y$ uses only the height
$Z+z^{\mathcal A}$ at $y$ and the jumps in $(x,y)$.
It does not use a jump at $y$ or any subsequent data.
This proves the endpoint assertion of
Remark~\ref{Remark_cut_endpoint}.
\end{proof}

We next sum over the four cases --- our goal is to show that the result is $0$. We expand each difference in
\eqref{eq_sigma_four_cases} into two signed terms and observe cancellations between them.

\begin{lemma}\label{Lemma_uniform_cancellation}
Fix $n\ge1$, $z_1,\ldots,z_n>0$, and $H>\bar z$.
For each $\mathcal A\subseteq\llbracket1,n\rrbracket$,
sum the weights
\begin{equation} \label{eq_x6}
 (-1)^{|\mathcal A|+|\mathcal B|}\sigma
\end{equation}
over Cases 1--4, all allowed choices of distinguished indices
outside $\mathcal A$, and all orders $\pi$ of $\mathcal B$
and $\rho$ of $\mathcal A$.
Here $\mathcal B$ contains the remaining indices, and each
summand has domain $0<Z<H-z^{\mathcal A}$.

After expanding $\sigma$ into signed indicators as in
\eqref{eq_sigma_four_cases}, the following cancellations hold:
\begin{enumerate}
\item
For each fixed $\mathcal A$ with $n-|\mathcal A|\ge2$,
the terms cancel in pairs at the same value of $Z$.
\item
For $\mathcal A=\llbracket1,n\rrbracket$ and
$\mathcal A=\llbracket1,n\rrbracket\setminus\{i\}$,
$1\le i\le n$, taken together, the terms cancel after
integration over their $Z$-domains.
\end{enumerate}
The same conclusions hold for the weights \eqref{eq_x6} replaced with
$(-1)^{|\mathcal A|+|\mathcal B|}\sigma
f(Z+z^{\mathcal A})$, for any function $f$ for which
the resulting integrals exist.
\end{lemma}

\begin{proof}
First fix $\mathcal A$ and an order $\rho$, and suppose
$n-|\mathcal A|\ge2$. Use the following correspondences:
\begin{enumerate}
\item
Pair the $u_0$ term of Case 1 with the $u_0$ term of
Case 2. Let $i$ be the first index in the order $\pi$
in Case 1, and let $\pi'$ be the order obtained by
deleting $i$. In Case 2, choose $i''=i$ and use the
order $\pi'$.

\item
Pair the $u_{i'}$ term of Case 3 with the $u_{i'}$ term
of Case 4, keeping $i'$ fixed. Let $i$ be the first
index in the order $\pi$ in Case 3, and let $\pi'$
be the order obtained by deleting $i$. In Case 4,
choose $i''=i$ and use the order $\pi'$.

\item
Pair the $-v_{i''}$ term of Case 2 with the $-v_{i''}$
term of Case 4, keeping $i''$ fixed. Let $j$ be the
first index in the order $\pi$ in Case 2, and let
$\pi'$ be the order obtained by deleting $j$.
In Case 4, choose $i'=j$ and use the order $\pi'$.
\end{enumerate}
Each correspondence is inverted by inserting the removed
index at the beginning of the order.
The indicator and the domain $0<Z<H-z^{\mathcal A}$
are unchanged, while $|\mathcal B|$ changes by one and
reverses the sign. The three correspondences are disjoint
and exhaust all expanded terms. The required first index
exists because $n-|\mathcal A|\ge2$.

Next, take
$\mathcal A=\llbracket1,n\rrbracket\setminus\{i\}$.
The first correspondence still pairs Case 1 with the
$u_0$ part of Case 2. The remaining terms are $-v_i$
in Case 2 and $u_i$ in Case 3.
Since $z^{\mathcal A}=\bar z-z_i$, the domain is
$0<Z<H-\bar z+z_i$, where $u_i=1$.
The two terms cancel for $Z<z_i$. Thus, apart from
endpoints of measure zero, one Case 3 term remains on
\[
 z_i<Z<H-\bar z+z_i
\]
for each order $\rho$ of $\mathcal A$, with sign
$(-1)^{n-1}$.

If $\mathcal A=\llbracket1,n\rrbracket$, only Case 1
is present. Its domain is $0<Z<H-\bar z$, where $u_0=1$.
For each order $\rho$ of $\llbracket1,n\rrbracket$,
let $i$ be its first index and let $\rho'$ be the order
obtained by deleting $i$. Pair this Case 1 term with
the remaining Case 3 term having
\[
 \mathcal A=\llbracket1,n\rrbracket\setminus\{i\},
 \qquad i'=i,
\]
and order $\rho'$, using the change of variable
$Z\mapsto Z+z_i$.
This translation maps $(0,H-\bar z)$ bijectively onto
$(z_i,H-\bar z+z_i)$. Moreover,
\[
 Z+\bar z=(Z+z_i)+(\bar z-z_i),
\]
so $Z+z^{\mathcal A}$ is preserved.
The two signs are $(-1)^n$ and $(-1)^{n-1}$.
The inverse subtracts $z_i$ and inserts $i$ at the
beginning of the order, proving the second assertion.

All pairings preserve $Z+z^{\mathcal A}$, so multiplication
by $f(Z+z^{\mathcal A})$ preserves these cancellations.
\end{proof}

\begin{proof}[Proof of
Proposition~\ref{Proposition_unicycle_zero_area}]

Fix $0<\epsilon<\min_j\bk_j$.
By Proposition~\ref{Proposition_tree_zero_area}, we may
assume that $\Gamma$ has a loop $m\to m$.
Lemma~\ref{Lemma_m0_loop_condition} reduces
$\bL^\Gamma_{[0];\,(\epsilon)}(\vec\bk,\vec\ttt)$
to the integral over blocks with $\bb_m=\varnothing$
whose first jump $x$ in the last interval is their unique
jump not to $0$; the contributions of all other blocks
cancel. These remaining blocks are partitioned into
Cases 1--4 above.

Fix the block data before $\bQ_{m-1}$, leaving open
the choice of which incoming function becomes $\bp_m$,
and fix $x,y,H$. The remaining labels $j>m$ are determined
by the order of first jumps, as in Definition~\ref{defn:blp}.
The incoming values $z_1,\ldots,z_n$ are now fixed, and
connectivity of $\Gamma$ and $m\ge2$ imply $n\ge1$.

For each case, choice of distinguished indices, sets
$\mathcal A,\mathcal B$, and orders $\pi,\rho$,
integrate over the times of all jumps to $0$ other than
the jump at $y$. By Lemma~\ref{Lemma_cut_integral},
the resulting expression is the weight
$(-1)^{|\mathcal A|+|\mathcal B|}\sigma$
multiplied by a common factor depending on $Z$ and
$\mathcal A$ only through $Z+z^{\mathcal A}$.
All factors coming from the earlier intervals are fixed.
Thus, Lemma \ref{Lemma_uniform_cancellation} gives zero
after summing over the four cases, distinguished indices,
subsets, and orders, and integrating over $Z$.

These cancellations take place within the fixed graph
$\Gamma$: the earlier edges are unchanged, every incoming
edge still ends at $m$, and the last interval contributes
the same single loop $m\to m$.

The cutoff in the last interval requires
$x>\bQ_{m-1}+\epsilon$, since $x$ is its first jump.
The cancellations leave $x$ and the earlier block data
unchanged, so they preserve this restriction.
Integrating over the remaining variables therefore gives
$\bL^\Gamma_{[0];\,(\epsilon)}(\vec\bk,\vec\ttt)=0.$
Sending $\epsilon\to0+$ in \eqref{eq_L_graph_r}
proves the proposition.
\end{proof}

\subsection{Proof of Proposition~\ref{Proposition_tree_one_area}}
\label{Section_tree_one_area}

The vanishing of the constant term in Propositions \ref{Proposition_tree_zero_area} and \ref{Proposition_unicycle_zero_area} means that $\bL^\Gamma_{\beta}(\vec\bk, \vec\ttt)$ has the form of a $\beta^{-1}$ prefactor multiplied by a power series in $\beta^{-1}$. In other words, the arguments of Sections \ref{Section_tree_constant}, \ref{Section_unicycle_zero_area} allowed us to gain one factor of $\beta^{-1}$. In contrast, Proposition \ref{Proposition_tree_one_area} requires us to find two such factors: since the constant term is equal to $0$ by Proposition \ref{Proposition_tree_zero_area}, vanishing of $\beta^{-1}$ term means divisibility by $\beta^{-2}$.

We note that the cancellations in the previous two sections were isolated in a single $(\bQ_{\ell-1},\bQ_{\ell})$ interval (in fact, we had $\ell=m$). In this section we reuse the same two mechanisms, but now apply them in two separate $(\bQ_{\ell-1},\bQ_{\ell})$ intervals.

There are two different situations. The first one corresponds to the first point of Definition \ref{Definition_group_of_graphs} and uses the cancellations of Section \ref{Section_tree_constant} twice.

\begin{lemma} \label{Lemma_two_cancellations}
 Suppose that a connected tree $\Gamma$ on $\llbracket 1,m\rrbracket$ has a vertex $\ell<m$ with no edges joining it with $\{\ell+1,\ell+2,\dots,m\}$. Then
 $\bL^\Gamma_{\beta;\,(\epsilon)}(\vec\bk,\vec\ttt)$ is divisible by $\beta^{-2}$ and consequently
      $\bL^\Gamma_{[1];\,(\epsilon)}(\vec\bk,\vec\ttt)=0$.
\end{lemma}
\begin{proof}
  In the interval $(\bQ_{\ell-1},\bQ_{\ell})$, all jumps are jumps to zero, and therefore the arguments of Section \ref{Section_tree_constant} apply to this interval equally well as to $(\bQ_{m-1},\bQ_m)$. (Here we subtract from the height process the total
value of the block functions which remain constant
throughout this interval.)  The matchings of that section can be run independently over the two segments, because they involve different block functions ($\bb_\ell=\varnothing$ matched with $\bb_\ell>0$; $\bb_m=\varnothing$ matched with $\bb_m>0$), and each of them gains us a factor of $\beta^{-1}$: we replace $\exp(\beta^{-1}(\cdot))$ by $1$ in the argument to see cancellation and $\exp\bigl(\beta^{-1}(\cdot)\bigr)-1$ is divisible by $\beta^{-1}$. Combining the two matchings together we get the desired $\beta^{-2}$ factor.
\end{proof}

In the rest of the section we assume that $m-1$ and $m$ are joined by an edge, which is always true if the second or third case in Definition \ref{Definition_group_of_graphs} hold, but the group does not satisfy the condition in the first case. We are going to exploit cancellations in the two last segments $(\bQ_{m-2},\bQ_{m-1})$ and $(\bQ_{m-1},\bQ_m)$. We are going to use the $\bb_m=\varnothing$ vs $\bb_m>0$ cancellation of Section \ref{Section_tree_constant} for the latter, and the four case cancellation of Section \ref{Section_unicycle_zero_area} for the former, thus, gaining two factors of $\beta^{-1}$. The new aspect is that the summation over the four cases now forces us sometimes to change the graph $\Gamma$. Hence the vanishing in Proposition \ref{Proposition_tree_one_area} needs a sum $\sum_{\Gamma\in G} \bL^\Gamma_{[1]}(\vec\bk, \vec\ttt)=0$ rather than an individual graph $\Gamma$. This is caused by the second cancellation in Lemma \ref{Lemma_uniform_cancellation} requiring a change in $\mathcal A$.

Fix a group $G$ and $0<\epsilon<\min_j\bk_j$.
We may assume that, in every tree of $G$, each vertex $\ell<m$
has an edge joining it to a vertex with a larger label. Indeed, replacing
$i\to m$ by $i\to(m-1)$ does not change which vertices have
such an edge. Thus, if this assumption fails, Lemma~6.12
applies to every tree in $G$ and proves the proposition.
There are $m-1$ edges altogether, so each vertex $\ell<m$
has exactly one edge to a larger vertex. By
Definition~\ref{Definition_connectivity_graph}, this means
that each interval $(\bQ_{\ell-1},\bQ_\ell)$, $\ell<m$,
contains exactly one jump which is not a jump to $0$.

The following analogue of Lemma~\ref{Lemma_m0_loop_condition}
reduces the penultimate interval to the same four cases as
in Section~6.2.

\begin{lemma}\label{Lemma_tree_first_jump}
For each $\Gamma\in G$, let $x$ be the unique jump not to $0$
in $(\bQ_{m-2},\bQ_{m-1})$. The total contribution to
$\bL^\Gamma_{[1];\,(\epsilon)}(\vec\bk,\vec\ttt)$
of all blocks except those satisfying
\[
 \bb_{m-1}=\varnothing,
 \qquad
 x\text{ is the first jump in }(\bQ_{m-2},\bQ_{m-1})
\]
is equal to zero.
\end{lemma}
\begin{proof}
Use the same matching as in the proof of
Lemma~\ref{Lemma_m0_loop_condition}, now in
$(\bQ_{m-2},\bQ_{m-1})$.
It pairs blocks with $\bb_{m-1}>0$ with blocks having
$\bb_{m-1}=\varnothing$ and a first jump to $0$ before $x$.
These pairs exhaust the blocks excluded in the statement.
As in Section~6.1, adding their signed contributions produces
a factor $\exp(\beta^{-1}z\bb_{m-1})-1$, where
$z=\bp_{m-1}(\bQ_{m-2}-)$ in the case with
$\bb_{m-1}>0$.

The matching only changes a function which vanishes before $x$,
so it does not change the functions entering the last interval.
We can therefore also apply the matching of Section~6.1 in
$(\bQ_{m-1},\bQ_m)$. Together the two matchings give a factor
divisible by $\beta^{-2}$, and hence no linear term in
$\beta^{-1}$. Both matchings preserve $\Gamma$ and the cutoff.
\end{proof}

It remains to integrate over the blocks satisfying the
conditions of Lemma~\ref{Lemma_tree_first_jump}.
Fix the block data before $\bQ_{m-2}$, leaving unspecified
which functions are assigned to $\bp_{m-1}$ and $\bp_m$.
Let
\[
 0<z_1<\cdots<z_n,\qquad
 \bar z=\sum_{i=1}^n z_i
\]
be all non-zero values of $\{\bp_j(\bQ_{m-2}-)\}_{j\ge 1}$.
Fix $x$ and $H=\bH(x)$, and let $Z=\bp_j(x)>0$ be the value immediately after
the jump at $x$. This function stays at $Z$ up to
$\bQ_{m-1}-$, since another jump before $\bQ_{m-1}$ would
produce a loop at $m-1$.

We use the same four cases and the same notation as in
Section \ref{Section_unicycle_zero_area}, with $\bp_m$ replaced by $\bp_{m-1}$.
Here $\mathcal A$ indexes the other incoming functions which
remain unchanged and positive up to $\bQ_{m-1}-$, and
$\mathcal B$ indexes those which jump to $0$ in
$(x,\bQ_{m-1})$. Thus, with
$z^{\mathcal A}=\sum_{i\in\mathcal A}z_i$,
\[
 \bH(\bQ_{m-1})=Z+z^{\mathcal A}.
\]
The distinguished indices $i',i''$ and the signed indicator
$\sigma$ are defined exactly as in
\eqref{eq_uniform_indicators}--\eqref{eq_sigma_four_cases}.
The admissible domain is again $H>\bar z$ and
$0<Z<H-z^{\mathcal A}$. We now use the endpoint version of
Remark~\ref{Remark_cut_endpoint}, with $\bQ_{m-2}$ in place
of $\bQ_{m-1}$ and $y=\bQ_{m-1}$.

To apply Lemma~\ref{Lemma_uniform_cancellation}, we identify
which choices of $\mathcal A$ belong to the group $G$.

\begin{lemma}\label{Lemma_tree_choices}
Delete the vertices $m-1,m$ and all their incident edges
from $\Gamma\in G$. The remaining graph is the same for
every $\Gamma\in G$ and has $n\ge1$ connected components,
one for each incoming value $z_i$.
The component corresponding to $z_i$ is attached by a single
edge to $m$ if $i\in\mathcal A$, and to $m-1$ otherwise.

Consequently, the trees in $G$ correspond to the following
choices:
\begin{enumerate}
\item
If $G$ is a singleton in the second case of
Definition~\ref{Definition_group_of_graphs}, then
$\mathcal A$ is fixed and $n-|\mathcal A|\ge2$.
\item
If $G=\{\Gamma_0,\ldots,\Gamma_n\}$ is in the third case,
then $\mathcal A=\llbracket1,n\rrbracket$ for $\Gamma_0$
and $\mathcal A=\llbracket1,n\rrbracket\setminus\{i\}$
for $\Gamma_i$, indexing these trees by the component
whose attachment is changed.
\end{enumerate}
\end{lemma}
\begin{proof}
By Definition~\ref{Definition_group_of_graphs}, the trees
in $G$ differ only in edges incident to $m-1$ or $m$.
Deleting these vertices and their incident edges therefore
gives the same graph for every $\Gamma\in G$.

Consider a connected component of this remaining graph,
with $v$ vertices. Since $\Gamma$ is a tree, the component
has $v-1$ internal edges. Each of its vertices has exactly
one edge to a larger vertex, so precisely one edge joins
the component to $\{m-1,m\}$. The block function
corresponding to this edge is positive at $\bQ_{m-2}-$,
identifying the component with its incoming value $z_i$.
Thus there are exactly $n$ components, and $n\ge1$
because $m\ge3$.

If $i\in\mathcal A$, the next jump of this function is at
$\bQ_{m-1}$ or later, so the edge ends at $m$.
If $i\notin\mathcal A$, the next jump is at $\bQ_{m-2}$,
at $x$, or at a jump to $0$ in $(x,\bQ_{m-1})$;
in each case the edge ends at $m-1$.

Hence $n-|\mathcal A|$ is the number of edges joining
$m-1$ to vertices in $\llbracket1,m-2\rrbracket$,
which gives the first assertion.
In the third case of Definition~\ref{Definition_group_of_graphs},
all components attach to $m$ in $\Gamma_0$,
while in $\Gamma_i$ only the component corresponding
to $z_i$ attaches to $m-1$. This gives the second assertion.

Figure \ref{Figure_tree_group} illustrates the second
assertion for $n=2$: moving the attachment of the
component corresponding to $z_i$ from $m$ to $m-1$
changes $\mathcal A=\{1,2\}$ to
$\mathcal A=\{1,2\}\setminus\{i\}$.
\end{proof}

\begin{proof}[Proof of Proposition \ref{Proposition_tree_one_area}]
By Lemma \ref{Lemma_two_cancellations}, we may assume that every vertex $\ell<m$
in each $\Gamma\in G$ has an edge to a larger vertex.
Lemma~\ref{Lemma_tree_first_jump} then reduces
$\bL^\Gamma_{[1];\,(\epsilon)}(\vec\bk,\vec\ttt)$
to the integral over blocks with $\bb_{m-1}=\varnothing$
whose first jump $x$ in $(\bQ_{m-2},\bQ_{m-1})$
is their unique jump not to $0$ in this interval.

Fix the block data before $\bQ_{m-2}$, leaving unspecified
which functions become $\bp_{m-1}$ and $\bp_m$.
The positive values at $\bQ_{m-2}-$ of incoming block functions are $z_1,\ldots,z_n$.
Choose one of the four cases from Figure \ref{Figure_four_cases}, its distinguished indices,
and the sets $\mathcal A,\mathcal B$ defined above:
$\mathcal A$ indexes the unchanged values surviving to
$\bQ_{m-1}-$, and $\mathcal B$ indexes the jumps to $0$
in $(x,\bQ_{m-1})$.
Keep $x,Z,H$ fixed and integrate over the times of these
jumps to $0$, in a fixed order $\pi$ of $\mathcal B$.
By Remark~\ref{Remark_cut_endpoint}, with
$y=\bQ_{m-1}$, the constant term of the factors from
the penultimate interval integrates to
\begin{equation}
\label{eq_x7}
 \frac{(-1)^{|\mathcal B|}\sigma}{2}
 \bF_0(x-\bQ_{m-2};H-\bar z)
 \bF_0(\bQ_{m-1}-x;H-Z-z^{\mathcal A}).
\end{equation}

The last interval has $|\mathcal A|+1$ incoming values:
$Z$ and $z_i$, $i\in\mathcal A$.
Apply Lemma~\ref{Lemma_first_matching_integral} to these
values, arranged in increasing order.
After summing over all choices in the last interval and
integrating over its jump times and virtual block, its
contribution has zero constant term and linear term
\begin{equation}
\label{eq_x8}
 \beta^{-1}(-1)^{|\mathcal A|}|\mathcal A|!
 J_\epsilon(Z+z^{\mathcal A}).
\end{equation}
Hence, to compute the coefficient of $\beta^{-1}$ in the
whole integral, all other factors can be replaced by their
constant terms.

Taking the coefficient of $\beta^{-1}$ in the product
of (6.11) and (6.12), and including the factor involving
$\vec\ttt$ at $\bQ_{m-1}$, we obtain, for each order $\pi$,
\begin{equation}
\label{eq_x9}
 \begin{split}
 &(-1)^{|\mathcal A|+|\mathcal B|}|\mathcal A|!\sigma
 \cdot\frac12
 \exp\!\left(
     \frac{\ttt_{m-1}-\ttt_m}{2}(Z+z^{\mathcal A})\right)
 \\
 &\qquad\times
 \bF_0(x-\bQ_{m-2};H-\bar z)
 \bF_0(\bQ_{m-1}-x;H-Z-z^{\mathcal A})
 J_\epsilon(Z+z^{\mathcal A}).
 \end{split}
\end{equation}
The factor $|\mathcal A|!$ is exactly the number of orders
$\rho$ of $\mathcal A$ in Lemma~\ref{Lemma_uniform_cancellation}.
Thus this expression is the sum over $\rho$ of the weights
in that lemma multiplied by the same function of
$Z+z^{\mathcal A}$. All factors from the earlier intervals
are fixed.

Lemma~\ref{Lemma_tree_choices} identifies the choices of
$\mathcal A$ with the two cancellations in
Lemma~\ref{Lemma_uniform_cancellation}.
For a singleton group, the first cancellation applies with
$\mathcal A$ fixed. For a group of $n+1$ trees, the second
cancellation applies after summing over the trees and
integrating over $Z$.
In both cases, the total contribution is zero.

The cutoff in the penultimate interval only requires
$x>\bQ_{m-2}+\epsilon$, since $x$ is its first jump,
and the cutoff in the last interval is already included in
$J_\epsilon$. These restrictions are preserved by the
cancellations. Integrating over $x,H$ and the remaining block data
before $\bQ_{m-2}$ gives
\[
 \sum_{\Gamma\in G}
 \bL^\Gamma_{[1];\,(\epsilon)}(\vec\bk,\vec\ttt)=0.
\]
Sending $\epsilon\to0+$ proves
the proposition.
\end{proof}

\section{The Gaussian limit of Laplace transforms}
\label{sec:covariance}

For any $\bk>0$, $\ttt\in\R$, $\eta>0$, define the centered Laplace transform of the Airy$_\beta$ line ensemble:
\begin{equation}\label{eq:X}
 X_\eta(\ttt,\bk)
 =
 \sum_{i=1}^{\infty}
 \exp\left(\frac{\eta\bk}{2}
 \cA_i^\beta(\eta^{-1/2}\ttt)\right)
 -
 \E\sum_{i=1}^{\infty}
 \exp\left(\frac{\eta\bk}{2}
 \cA_i^\beta(\eta^{-1/2}\ttt)\right).
\end{equation}
For any $(\bk_1,\ttt_1),(\bk_2,\ttt_2)\in(0,\infty)\times\R$,
we denote
\begin{equation}\label{eq:C}
 C(\bk_1,\ttt_1; \bk_2,\ttt_2)
 :=
 \frac{1}{\pi\sqrt{\bk_1\bk_2}}
 \int_0^\infty g
 \exp\left(
 -\frac{|\ttt_1-\ttt_2|g}{2}
 -\frac{g^2}{2\bk_1}
 -\frac{g^2}{2\bk_2}
 \right)\dd g.
\end{equation}

\begin{theorem}\label{Theorem_Laplace_CLT}
For any $m\in\N$ and
$(\vec\bk,\vec\ttt)\in\R_{>0}^m\times\R^m$, as
$\eta\to0+$, the vector
$\bigl(X_\eta(\ttt_i,\bk_i)\bigr)_{i=1}^m$
converges in all mixed moments and in distribution to a
centered Gaussian vector with covariance
\[
 \left(\beta^{-1}C(\bk_i,\ttt_i;\bk_j,\ttt_j)\right)_{i,j=1}^m.
\]
\end{theorem}
\begin{remark}
 The asymptotic variance is computed by plugging $\ttt_1=\ttt_2$, $\bk_1=\bk_2$ into \eqref{eq:C} and explicitly evaluating the integral. The result depends on neither $\ttt$, nor $\bk$:
 \begin{equation}
 \label{eq_variance_Laplace}
  \lim_{\eta\to 0+}  \E\bigl[X_\eta^2(\ttt,\bk)\bigr]=\frac{1}{2\pi\beta}.
 \end{equation} 
\end{remark}
\begin{proof}[Proof of Theorem \ref{Theorem_Laplace_CLT}]
Recall the moment-cumulant formula \eqref{eq_moments_cumulants} and that the convergence of moments implies the distributional convergence for Gaussian limits. Hence, it suffices to show that for every $m\ge2$ and fixed
$(\vec\bk,\vec\ttt)\in(0,\infty)^m\times\R^m$, we have
\begin{equation}   \label{eq:Cumetalim}
 \E^{(\textnormal{c})}\bigl[ X_\eta(\ttt_1,\bk_1),\ldots,
           X_\eta(\ttt_m,\bk_m)\bigr]
 =\begin{cases}
 \beta^{-1}C(\bk_1,\ttt_1;\bk_2,\ttt_2)+O(\eta^{3/2}),&m=2,\\
 O(\eta^{3/2}),&m\ge3,
 \end{cases}
 \qquad \eta\to 0+.
\end{equation}
By symmetry we may assume that $\ttt_1\le \cdots \le \ttt_m$. In the rest of the proof we establish \eqref{eq:Cumetalim}.

We note that the removal of the deterministic second term $\E$ in \eqref{eq:X} does not change the cumulants and use the graph expansion \eqref{eq_expected_Laplace_cum_rescaled} of Corollary~\ref{Corollary_cumulant_expansion}.
Every connected graph has $e\ge m-1$.
Proposition~\ref{Proposition_tree_zero_area} removes the
$r=0$, $e=m-1$ terms, and
Proposition~\ref{Proposition_unicycle_zero_area} removes
the $r=0$, $e=m$ terms.

For $m\ge3$,
Proposition~\ref{Proposition_tree_one_area} then further removes the $r=1$, $e=m-1$ terms, since all trees can be split into disjoint groups presented in Definition \ref{Definition_group_of_graphs}.
Thus, when $m\ge3$, all remaining terms have $e+r\ge m+1$ and the factor $\eta^{\frac32(e+r-m)}$ tends to $0$ at speed at least $\eta^{3/2}$ in each non-zero term of \eqref{eq_expected_Laplace_cum_rescaled}. Since the series at any $0<\eta<1$ is dominated by the absolutely convergent $\eta=1$ series, the dominated convergence theorem yields that the sum is $O(\eta^{3/2})$, as claimed in \eqref{eq:Cumetalim}.

\medskip

For $m=2$, the terms with $e+r>2$ contribute through $\eta^{\frac32(e+r-m)}$ to the $O(\eta^{3/2})$ remainder in \eqref{eq:Cumetalim}. The only remaining term with $e+r\le2$ is
$$
\beta^{-1}\bL^{\Gamma_{12}}_{[1]}(\vec\bk,\vec\ttt),
$$
where $\Gamma_{12}$ is the graph on $\{1,2\}$ with a single edge
$1\to2$. Therefore, it remains to prove that
\begin{equation}
\label{eq_x14}
\bL^{\Gamma_{12}}_{[1]}(\vec\bk,\vec\ttt)
=C(\bk_1,\ttt_1;\bk_2,\ttt_2).
\end{equation}
Let us describe the blocks with $\Gamma(\vec\bp)=\Gamma_{12}$.
Each block function with index $j>2$ has at least two jumps, so
$\bdel\ge2\bu$. Since $\bdel-\bu=1$, there are only two
possibilities:
\begin{enumerate}
\item
$\bu=0$, $\bdel=1$. The only nonzero block function is
$\bp_2=g\don_{[x,\bQ_1)}$, for some $0<x<\bk_1$ and $g>0$.
By Lemma~\ref{Lemma_no_loop_vanishing}, we have
$\bb_2=v$ for some $0<v<\bk_2$.
\item
$\bu=1$, $\bdel=2$. The only nonzero block function is
$\bp_3=g\don_{[x,\bQ_1+v)}$, for some $0<x<\bk_1$,
$0<v<\bk_2$, and $g>0$. Both virtual coordinates are empty.
\end{enumerate}
In both cases $\bb_1=\varnothing$, and the graph has the same
single edge $1\to2$. The prescribed heights are
\[
 \bH(0)=\bH(\bQ_2)=0,\qquad
 \bH(\bQ_1)=\bH(\bQ_1+v)=g.
\]
The remaining height can be written as $\bH(x)=g+h$, with
$h>0$.

We observe that these are precisely the two cases matched in the proof of
Proposition~\ref{Proposition_tree_zero_area}.
We now consider the integrand of $\bL^{\Gamma_{12}}_{[1];\,(\epsilon)}$, from \eqref{eq_L_graph_cutoff_r}.
We start with $\bI_{\beta,\vec\bk}[\vec\bp,\vec\bb, \bH]$ from Definition \ref{def:xibi}.
The factor
$2^{-\bdel-|\{\ell:\bb_\ell\ne\varnothing\}|}$ in
\eqref{eq:defnibk} is $1/4$ in both cases. The first case has
the additional factor $\exp(\beta^{-1}gv)$, while the second
has a minus sign from the jump to zero at $\bQ_1+v$.
Thus by pairing up terms in both cases with the same $x, v, g, h$, for the integrand we get
\[
 \frac14\exp\left(-\frac{\ttt_2-\ttt_1}{2}g\right)
 \bI_0(x;g+h)\bI_0(\bk_1-x;h)\bI_{0,0}(v)\bI_0(\bk_2-v;g)
 (\exp(\beta^{-1}gv) - 1).
\]
By taking the coefficient of $\beta^{-1}$ in its expansion, we get the paired integrand in $\bL^{\Gamma_{12}}_{[1];\,(\epsilon)}$ as
\[
 \frac14\exp\left(-\frac{\ttt_2-\ttt_1}{2}g\right)
 \bF_0(x;g+h)\bF_0(\bk_1-x;h)\bF_{0,0}(v)\bF_0(\bk_2-v;g) gv.
\]
For the integration domain, the matching preserves the cutoff, which requires that
$x, v>\epsilon$.
Thus by further taking
$\epsilon\to0+$, we get
\begin{multline}\label{eq:covariance-integral}
 \bL^{\Gamma_{12}}_{[1]}(\vec\bk,\vec\ttt)
 =\frac14
 \int_{\substack{0<x<\bk_1,\ 0<v<\bk_2\\g,h>0}}
 \exp\left(-\frac{\ttt_2-\ttt_1}{2}g\right) \\
 \times
 \bF_0(x;g+h)\bF_0(\bk_1-x;h)\bF_{0,0}(v)\bF_0(\bk_2-v;g)
 gv
 \dd x\dd v\dd g\dd h.
\end{multline}
It now remains to evaluate this integral.
By Lemma \ref{lem:FP_convolution},
\[
 \bF_0(x;g+h)
 =\frac12\int_0^x\bF_0(s;g)\bF_0(x-s;h)\dd s.
\]
Plugging this into \eqref{eq:covariance-integral}, we
then integrate over $h$, as
\begin{align*}
 \int_0^\infty
 \bF_0(x-s;h)\bF_0(\bk_1-x;h)\dd h
 &=
 \frac{2}{\pi[(x-s)(\bk_1-x)]^{3/2}}
 \int_0^\infty h^2
 \exp\left(-\frac{\bk_1-s}{2(x-s)(\bk_1-x)}h^2\right)\dd h\\
 &=
 \frac{2}{\sqrt{2\pi(\bk_1-s)^3}}
 =\bF_{0,0}(\bk_1-s).
\end{align*}
The integrand no longer depends on $x$. Thus the integration of $x$ over
$s<x<\bk_1$ gives another factor of $\bk_1-s$. Changing variables by
$t_1=\bk_1-s$ and $t_2=v$, we now obtain
\begin{equation}
\label{eq_cov_factorized}
 \bL^{\Gamma_{12}}_{[1]}(\vec\bk,\vec\ttt)
 =\frac18\int_0^\infty
 g\exp\left(-\frac{\ttt_2-\ttt_1}{2}g\right) \prod_{j=1}^2\left[\int_0^{\bk_j}
 t_j\bF_{0,0}(t_j)\bF_0(\bk_j-t_j;g)\dd t_j\right]\dd g.
\end{equation}
We note that from \eqref{eq_I0} and \eqref{eq_I00}, for any $\bk, g>0$,
\begin{multline*}
 \int_0^\bk t\bF_{0,0}(t)\bF_0(\bk-t;g)\dd t
 =
 \frac{2g}{\pi}\int_0^\bk
 t^{-1/2}(\bk-t)^{-3/2}
 \exp\left(-\frac{g^2}{2(\bk-t)}\right)\dd t\\
 =
 \frac{4g}{\pi\bk}\exp\left(-\frac{g^2}{2\bk}\right)
 \int_0^\infty
 \exp\left(-\frac{g^2u^2}{2\bk}\right)\dd u
 =
 \sqrt{\frac{8}{\pi\bk}}\exp\left(-\frac{g^2}{2\bk}\right),
\end{multline*}
where the second equality is due to a substitution $u=\sqrt{t/(\bk-t)}$.
Plugging this into \eqref{eq_cov_factorized} proves \eqref{eq_x14}.
\end{proof}

We present an alternative expression for the covariance \eqref{eq:C}, which is helpful in the identification with the Gaussian free field.

\begin{lemma}\label{lem:pullback-energy}

For any $(\bk_1, \bk_2, \ttt_1, \ttt_2) \in \R_{>0}^2\times \R^2$,
\begin{equation}\label{eq_covariance_GFF}
C(\bk_1,\ttt_1;\bk_2,\ttt_2)=
 \int_{-\infty}^{+\infty}\int_{-\infty}^{+\infty}
 \frac{\bk_1}{2}\exp\left(\frac{\bk_1 x}{2}\right) \cdot \frac{\bk_2}{2}\exp\left(\frac{\bk_2 y}{2}\right) \cdot
 K(\ttt_1,x;\, \ttt_2,y)\dd x\,\dd y,
\end{equation}
where
$$
  K(\ttt_1,x;\, \ttt_2,y)
 =\begin{cases}
 \displaystyle\frac{1}{2\pi^2}
 \log\frac{(\ttt_1-\ttt_2)^2/4+(\sqrt{-x}+\sqrt{-y})^2}
          {(\ttt_1-\ttt_2)^2/4+(\sqrt{-x}-\sqrt{-y})^2},
       &x,y<0,\\[3mm]
 0,    &x\ge0\text{ or }y\ge0.
 \end{cases}
$$
\end{lemma}

\begin{proof}
Set $a=|\ttt_1-\ttt_2|/2$.  We consider only the case $a>0$: $a=0$ is obtained by sending $\ttt_1\to\ttt_2$ in both sides of \eqref{eq_covariance_GFF}. Changing variables by
$x=-u^2$, $y=-v^2$, the right-hand side of \eqref{eq_covariance_GFF} becomes
\begin{equation}
\label{eq_x15}
 \frac{\bk_1\bk_2}{2\pi^2}
 \int_0^\infty\int_0^\infty
 uv\exp\left(-\frac{\bk_1u^2+\bk_2v^2}{2}\right)
 \log\frac{a^2+(u+v)^2}{a^2+(u-v)^2}\,\dd u\,\dd v.
\end{equation}
We use the Fourier-sine representation of the logarithmic kernel: for  any $u,v\in\R$ we have
\begin{multline*}
  4\int_0^\infty e^{-ag}
 \frac{\sin(gu)\sin(gv)}{g}\dd g= 2\int_0^\infty e^{-ag}
 \frac{\cos((u-v)g)-\cos((u+v)g)}{g}\dd g
 \\ =2\int_0^\infty
 \int_{u-v}^{u+v} e^{-ag} \sin(tg)\, \dd t\dd g= 2  \int_{u-v}^{u+v} \frac{t}{a^2+t^2}\, \dd t=
 \log\frac{a^2+(u+v)^2}{a^2+(u-v)^2}.
\end{multline*}
We plug this  into \eqref{eq_x15} and integrate in $u$ and $v$ using integration by parts and Gaussian integrals:
\begin{equation}
\label{eq:covariance-Gaussian-sine}
 \int_0^\infty \bk u
 \exp\left(-\frac{\bk u^2}{2}\right)\sin(gu)\dd u
 =
 g\int_0^\infty
 \exp\left(-\frac{\bk u^2}{2}\right)\cos(gu)\dd u
 =
 g\sqrt{\frac{\pi}{2\bk}}
 \exp\left(-\frac{g^2}{2\bk}\right).
\end{equation}
As a result, \eqref{eq_x15} becomes
\begin{multline*}
   \frac{2}{\pi^2}\int_0^\infty  \left[
 \int_0^\infty
 \bk_1 u \exp\left(-\frac{\bk_1u^2}{2}\right)
\sin(gu)\, \dd u \int_0^\infty \bk_2 v \sin(gv)  \exp\left(-\frac{\bk_2v^2}{2}\right) \,\dd v\right]
e^{-ag}
 \frac{\dd g}{g}
\\=   \frac{2}{\pi^2}\int_0^\infty  \left[  g\sqrt{\frac{\pi}{2\bk_1}}
 \exp\left(-\frac{g^2}{2\bk_1}\right) \cdot g\sqrt{\frac{\pi}{2\bk_2}}
 \exp\left(-\frac{g^2}{2\bk_2}\right)\right]
e^{-ag}
 \frac{\dd g}{g},
\end{multline*}
which matches \eqref{eq:C}.
\end{proof}

\section{GFF limit for the Airy surface}
\label{sec:pullback-proof}

In this section we prove Theorem \ref{thm:surface}.

\subsection{Preliminaries on the Gaussian free field}
\begin{definition} \label{Definition_GFF_informal}
The Gaussian free field on the upper half-plane $\U=\{z\in\C:\operatorname{Im}z>0\}$ is a mean $0$ generalized Gaussian random field $\GFF_{\U}(z)$, $z\in\U$, with covariance kernel
\begin{equation}\label{eq_GFF_U_kermel}
 G_{\U}(z,w)
 =-\frac{1}{2\pi}\log\left|\frac{z-w}{z-\overline w}\right|.
\end{equation}
\end{definition}
Note that $G_{\U}(z,z)=\infty$ and the values of $\GFF_{\U}(z)$ at individual points are not defined. Instead, we can pair it with sufficiently smooth signed test-measures: if $\nu$ is such a measure on $\U$, then $\int_{\mathbb U} \GFF_{\U}(z) \dd \nu$ is a mean $0$ Gaussian random variable. For two test-measures $\nu$ and $\mu$ the covariance between the corresponding variables is
$$
 \E\left[ \int_{\mathbb U}  \GFF_{\U}(z) \dd \nu \int_{\mathbb U}  \GFF_{\U}(z) \dd \mu\right]=\iint_{\U\times \U} G_\U(z,w) \dd \nu \dd \mu.
$$
An alternative to using pairings as a definition is to define a functional space to which the field $\GFF_{\U}$ belongs; \cite{Sheffield,WernerPowell2022GFF} contain further details on the definition of the Gaussian free field.

\medskip

In order to describe the fluctuations of the Airy surface through the Gaussian free field, we identify its coordinates $(\ttt,x)\in \mathbb R^2$ with $\U$. This is done by first restricting to the half-plane $x<0$ and using the map:
\begin{equation}
\label{eq_Omega_map}
 \Omega: \mathbb R\times \mathbb R_{<0}\to \U:\qquad  \Omega(\ttt,x)=\frac{\ttt}{2}+\ii\sqrt{-x}.
\end{equation}
Pairings of the pullback of $\GFF_\U$, denoted $\GFF_\U\circ \Omega$, with test-measures are computed as follows.
\begin{definition} Consider a signed test-measure $\nu$ on $\mathbb R\times \mathbb R_{<0}$, such that the integral \eqref{eq_GFF_push_cov} for $\mu=\nu$ is finite. The pairing of $\GFF_\U\circ \Omega$ with $\nu$ is a mean $0$ Gaussian random variable. For two such test-measures, $\nu$ and $\mu$, the covariance of the pairings is computed as
 \begin{multline}
 \label{eq_GFF_push_cov}
  \iint_{\U\times \U} G_\U(z,w)\, \dd (\Omega^* \nu) \dd (\Omega^* \mu) \\= \int_{\mathbb R\times \mathbb R_{<0}\times \mathbb R\times \mathbb R_{<0}} -\frac{1}{2\pi}\log\left|\frac{\Omega(\ttt_1,x_1)-\Omega(\ttt_2,x_2)}{\Omega(\ttt_1,x_1)-\overline \Omega(\ttt_2,x_2)}\right|
   \nu(\dd \ttt_1,\dd x_1) \mu(\dd \ttt_2,\dd x_2),
 \end{multline}
 where $\Omega^* \nu$ and  $\Omega^* \mu$ are pushforwards of the measures $\nu$ and $\mu$, respectively.
\end{definition}

The Airy surface is a random field on the plane $\mathbb R\times \mathbb R$, while $\GFF_\U\circ \Omega$ is a random field on the half-plane $\mathbb R\times \mathbb R_{<0}$. We fix this discrepancy by extending the field with zeros on the complementary half-plane.

\begin{definition} \label{Definition_extended_GFF}
 ${\overline{\GFF_\U\circ \Omega}}$ is a mean $0$ generalized Gaussian random field on $\mathbb R^2$ with coordinates $(\ttt,x)$
 with covariance given using \eqref{eq_Omega_map} by
\begin{equation}
 \mathrm{Cov}(\ttt_1,x_1;\, \ttt_2,x_2)=\begin{cases}  -\frac{1}{2\pi}\log\left|\frac{\Omega(\ttt_1,x_1)-\Omega(\ttt_2,x_2)}{\Omega(\ttt_1,x_1)-\overline \Omega(\ttt_2,x_2)}\right|, &x_1,x_2<0,\\
 0, & \text{otherwise.}
 \end{cases}
\end{equation}
For a signed test-measure $\nu$ on $\R^2$, we denote by $\langle \nu\rangle$ the pairing of $\nu$ with ${\overline{\GFF_\U\circ \Omega}}$. The covariance for the pairings is computed by restricting the test-measures onto $\R\times \R_{<0}$ and then using \eqref{eq_GFF_push_cov}.
\end{definition}

\subsection{Test-measures}

We are going to prove the convergence of the centered Airy surface to GFF on two separate classes of signed test-measures.

\begin{definition} \label{Definition_class_I_measure} Take $t\in\mathbb R$ and a continuous function $\phi:\R\to\R$ such that there exist
constants $B,\lambda,\Lambda>0$ for which
\begin{equation}\label{eq:admissible-spatial-test}
 |\phi(x)|\le B\bigl(\exp(\lambda x)+\exp(\Lambda x)\bigr),
 \qquad x\in\R.
\end{equation}
We define a \emph{class I test-measure} $\nu_{t;\phi}(\dd \ttt,\dd x)$ on  $\R^2$ as a measure supported on the line $\ttt=t$ of density $\phi(x)$ with respect to the linear Lebesgue measure.
\end{definition}

\begin{definition} \label{Definition_class_II_measure} Consider a function $f:\R^2\to\R$ which is Borel measurable, with
$f(\ttt,\cdot)$ continuous for every $\ttt\in\R$, and such that for a nonnegative function $B\in L^1(\R)$ and constants $\lambda,\Lambda>0$ we have 
\begin{equation}\label{eq:admissible-spacetime-test}
 |f(\ttt,x)|\le B(\ttt)
 \bigl(\exp(\lambda x)+\exp(\Lambda x)\bigr), \qquad \ttt,x\in\mathbb R.
\end{equation}
We define a \emph{class II test-measure} $\nu_{f}(\dd \ttt,\dd x)$ on $\R^2$ as a measure of density $f(\ttt,x)$ with respect to the planar Lebesgue measure.
\end{definition}

\begin{lemma}
 The field ${\overline{\GFF_\U\circ \Omega}}$ of Definition \ref{Definition_extended_GFF} has well-defined pairings with signed measures of Definitions \ref{Definition_class_I_measure} and \ref{Definition_class_II_measure}. Take two such pairings, $\langle \nu_1\rangle$ and $\langle \nu_2\rangle$. These are mean $0$ Gaussian random variables, whose covariance is:
 \begin{itemize}
  \item If $\nu_1$ and $\nu_2$ are class I measures: $\nu_1=\nu_{t_1;\phi_1}$ and $\nu_2=\nu_{t_2;\phi_2}$, then
  \begin{equation}
  \label{eq_cov_I_I}
   \E \langle \nu_1\rangle \langle \nu_2\rangle=  -\frac{1}{2\pi} \iint_{\R_{<0}^2}\log\left|\frac{\Omega(t_1,x_1)-\Omega(t_2,x_2)}{\Omega(t_1,x_1)-\overline \Omega(t_2,x_2)}\right| \phi_1(x_1) \phi_2(x_2)\dd x_1 \dd x_2.
  \end{equation}
    \item If $\nu_1$ and $\nu_2$ are class II measures: $\nu_1=\nu_{f_1}$ and $\nu_2=\nu_{f_2}$, then
  \begin{equation}
  \label{eq_cov_II_II}
   \E \langle \nu_1\rangle \langle \nu_2\rangle=   -\frac{1}{2\pi} \int_{\mathbb R\times \mathbb R_{<0}\times \mathbb R\times \mathbb R_{<0}} \log\left|\frac{\Omega(\ttt_1,x_1)-\Omega(\ttt_2,x_2)}{\Omega(\ttt_1,x_1)-\overline \Omega(\ttt_2,x_2)}\right|
   f_1(\ttt_1, x_1) f_2(\ttt_2,x_2)\, \dd \ttt_1 \dd x_1 \dd \ttt_2 \dd x_2.
  \end{equation}
 \item If $\nu_1$ is a class I measure, $\nu_1=\nu_{t,\phi}$, while $\nu_2$ is class II, $\nu_2=\nu_f$, then
  \begin{equation}
  \label{eq_cov_I_II}
   \E \langle \nu_1\rangle \langle \nu_2\rangle=  -\frac{1}{2\pi} \int_{\mathbb R_{<0}\times \mathbb R\times \mathbb R_{<0}} \log\left|\frac{\Omega(t,x_1)-\Omega(\ttt_2,x_2)}{\Omega(t,x_1)-\overline \Omega(\ttt_2,x_2)}\right|
   \phi(x_1) f(\ttt_2,x_2)\, \dd x_1 \dd \ttt_2 \dd x_2.
  \end{equation}
 \end{itemize}
\end{lemma}
We omit the proof, which directly follows from the definitions.

\subsection{Detailed statement of the convergence theorem}
The Airy$_\beta$ surface is defined via
\[
 \cA^\beta(\ttt,x)
 =\#\{i\in\N:\cA_i^\beta(\ttt)>x\},
\]
and its centered, rescaled version
\[
 H_\eta(\ttt,x)
 =\cA^\beta(\eta^{-1/2}\ttt,\eta^{-1}x)
 -\E[\cA^\beta(\eta^{-1/2}\ttt,\eta^{-1}x)],
 \qquad \eta>0.
\]
We call $H_\eta$ the centered \emph{height function} of the Airy$_\beta$ line ensemble, by analogy with random matrix and 2d statistical mechanics literature, see \cite{Borodin,BG,GorinBook}.

\begin{theorem} \label{Theorem_main_detailed}
 Fix $\beta>0$. As $\eta\to 0 +$ the random field $\Bigl[\sqrt{\tfrac{\beta \pi}{2}}H_\eta(\ttt,x)\Bigr]_{\ttt,x\in\R}$, converges to the $\Omega$--pullback of the Gaussian free field in the upper half-plane, i.e.\ to the field  ${\overline{\GFF_\U\circ \Omega}}$  of Definition \ref{Definition_extended_GFF}, in the following sense. For any $m\in\N$ and $m$ signed test-measures $\nu^1,\dots,\nu^m$ from either of the classes of Definitions \ref{Definition_class_I_measure}, \ref{Definition_class_II_measure}, the $m$--dimensional random vector
 \begin{equation}
 \label{eq_x17}
  \left( \sqrt{\tfrac{\beta \pi}{2}} \iint_{\R^2} H_\eta(\ttt,x) \nu^j(\dd \ttt,\dd x) \right)_{j=1}^m
 \end{equation}
 converges in distribution and in moments to a mean $0$ Gaussian vector whose covariance is given by \eqref{eq_cov_I_I},\eqref{eq_cov_II_II},\eqref{eq_cov_I_II}.
\end{theorem}

In the rest of this section, we prove Theorem \ref{Theorem_main_detailed}. We start from test-measures with exponential density in Section \ref{Section_exponential_test}, re-interpreting Theorem \ref{Theorem_Laplace_CLT}. In Section \ref{Section_class_I} we extend convergence to all class I test-measures by using a comparison inequality, based on the positive association in the Airy$_\beta$ line ensemble. In Section \ref{Section_class_II} we reach class II measures by integration in $\ttt$ direction.

\subsection{Convergence for test-measures with exponential density}
\label{Section_exponential_test}

Take $t_1,\dots,t_m\in\R$, $s_1,\dots,s_m\in \R_{>0}$ and suppose that the test-measures are of class I with
$$
 \nu^j=\nu_{t_j, \exp(s_j x)} \qquad 1\le j \le m.
$$
Then,
\begin{multline} \label{eq_x16}
  \iint_{\R^2} H_\eta(\ttt,x) \nu^j(\dd \ttt,\dd x)=\int_{-\infty}^{\infty}  H_\eta(t_j,x) \exp(s_j x)  \dd x \\ =
  \int_{-\infty}^{\infty} \cA^\beta(\eta^{-1/2}t_j,\eta^{-1}x)\exp(s_j x)  \dd x
 -\int_{-\infty}^{\infty} \E[\cA^\beta(\eta^{-1/2}t_j,\eta^{-1}x)] \exp(s_j x)  \dd x
 \\ =
  \int_{-\infty}^{\infty} \sum_{i=1}^{\infty} \don[x<\eta \cA_i^\beta(\eta^{-1/2}t_j)] \exp(s_j x)  \dd x
 -\E \int_{-\infty}^{\infty} \sum_{i=1}^{\infty} \don[x<\eta \cA_i^\beta(\eta^{-1/2}t_j)] \exp(s_j x)  \dd x
 \\= \frac{1}{s_j} \left[ \sum_{i=1}^{\infty} \exp\left(\eta s_j \cA_i^\beta(\eta^{-1/2}t_j)\right) - \E \sum_{i=1}^{\infty} \exp\left(\eta s_j \cA_i^\beta(\eta^{-1/2}t_j)\right)\right],
\end{multline}
where in the last equality we changed the order of integration and summation (which is justified by positivity) and then computed the integral of the exponential function. The last line of \eqref{eq_x16} can be identified with \eqref{eq:X} as
\begin{equation}
\label{eq_x25}
 \frac{1}{s_j} X_\eta\left(t_j, 2 s_j\right).
\end{equation}
Thus, applying Theorem \ref{Theorem_Laplace_CLT} and Lemma \ref{lem:pullback-energy}, the vector \eqref{eq_x17} converges to a mean $0$ Gaussian vector with covariance
\begin{multline}
  \frac{\pi}{2 s_1 s_2} \int_{-\infty}^{\infty}\int_{-\infty}^{\infty}
 s_1\exp(s_1 x) s_2\exp(s_2 y)
 K(t_1,x;\, t_2,y)\, \dd x\,\dd y
 \\=
  \frac{1}{4\pi} \int_{-\infty}^{0}\int_{-\infty}^{0}
 \exp(s_1 x) \exp(s_2 y)
 \log\left[\frac{(t_1-t_2)^2/4+(\sqrt{-x}+\sqrt{-y})^2}
          {(t_1-t_2)^2/4+(\sqrt{-x}-\sqrt{-y})^2}\right]\dd x\, \dd y.
\end{multline}
Using $\Omega(t,x)=\frac{t}{2}+\ii\sqrt{-x}$ to identify
$$
 \log\left[\frac{(t_1-t_2)^2/4+(\sqrt{-x}+\sqrt{-y})^2}
          {(t_1-t_2)^2/4+(\sqrt{-x}-\sqrt{-y})^2}\right]=-2\log\left|\frac{\Omega(t_1,x)-\Omega(t_2,y)}{\Omega(t_1,x)-\overline \Omega(t_2,y)}\right|
$$
we get a match with \eqref{eq_cov_I_I}, where $\phi_1(x)=\exp(s_1 x)$ and $\phi_2(y)=\exp(s_2 y)$.

\subsection{Convergence for class I test-measures}
\label{Section_class_I}

We next pass from  exponential to general densities of class I test-measures, i.e.\ we deal with
\begin{equation}
\label{eq_x26}
Y_\eta(\ttt;\phi)=\int\phi(x)H_\eta(\ttt,x)\dd x,
\end{equation}
where $\phi$ satisfies the bound \eqref{eq:admissible-spatial-test}. We need the following moment comparison.
\begin{lemma}\label{lem:surface-association}
Fix $\eta>0$ and $\ttt\in\R$. Let measurable $\phi,\psi:\R\to\R$
satisfy \eqref{eq:admissible-spatial-test}
and $|\phi|\le\psi$. Then
\begin{equation}
\label{eq_x20}
 \E \left[Y^2_\eta(\ttt;\phi)\right]
 \le\E\left[Y^2_\eta(\ttt;\psi)\right].
\end{equation}
Moreover, for every positive integer $r$ we have\footnote{The constant $2^{2r}$ in \eqref{eq_x21} is not optimal.} 
\begin{equation}
\label{eq_x21}
 \E\left[Y^{2r}_\eta(\ttt;\phi)\right]
 \le 2^{2r}\E\left[Y^{2r}_\eta(\ttt;\psi)\right].
\end{equation}
\end{lemma}
For related bounds on fluctuations of linear statistics, see, e.g., the Poincaré and logarithmic Sobolev inequalities for Gaussian $\beta$-ensembles in \cite{CL20} and the concentration inequalities for determinantal ensembles in \cite[Section 3]{BD14}.

Our proof of Lemma \ref{lem:surface-association} is based on the positive association property of the Airy particles, which is a version of the FKG inequality \cite{FKG} and can be deduced from the log-supermodularity of the density of the Gaussian $\beta$ ensemble.

\begin{lemma}[{\cite[Corollary 5(2)]{BKM}}]\label{lem:airy-association}
Fix $\beta>0$, $k\in\N$, and $\ttt\in\R$, and set
$Z=(\cA_i^\beta(\ttt))_{i=1}^k$.
Then, for any increasing Borel sets $E,F\subset\R^k$,
\[
 \PP(Z\in E\cap F)\ge\PP(Z\in E)\PP(Z\in F).
\]
Here a set $E$ is increasing if $x\in E$ and $y\ge x$
coordinatewise imply $y\in E$.
\end{lemma}
\begin{proof}[Proof of Lemma \ref{lem:surface-association}]
For any non-negative function $\varphi:\R\to \R$ satisfying
\eqref{eq:admissible-spatial-test}, the
definition of $H_\eta$ gives
\[
 Y_\eta(\ttt;\varphi)
 =\sum_{i=1}^\infty
   \int_{-\infty}^{\eta\cA_i^\beta(\eta^{-1/2}\ttt)}\varphi(x)\dd x -\E\left[\sum_{i=1}^\infty
   \int_{-\infty}^{\eta\cA_i^\beta(\eta^{-1/2}\ttt)}\varphi(x)\dd x\right].
\]
For every $a\in\R$, we have an identity of events
\[
 \left\{\sum_{i=1}^\infty
 \int_{-\infty}^{\eta\cA_i^\beta(\eta^{-1/2}\ttt)}\varphi(x)\dd x>a\right\}
 =\bigcup_{k=1}^\infty
 \left\{\sum_{i=1}^k
 \int_{-\infty}^{\eta\cA_i^\beta(\eta^{-1/2}\ttt)}\varphi(x)\dd x>a\right\}.
\]
The events on the right-hand side are increasing events of the positions of the first $k$
particles of the Airy point process and increase with $k$. Thus applying Lemma \ref{lem:airy-association} and sending $k\to\infty$, for any $\varphi, \zeta\ge  0$ satisfying \eqref{eq:admissible-spatial-test}, and $a,b\in\R$, we have
\begin{equation}
\label{eq_x18}
\PP(Y_\eta(\ttt;\varphi)>a,\;Y_\eta(\ttt;\zeta)>b) \ge \PP(Y_\eta(\ttt;\varphi)>a)\PP(Y_\eta(\ttt;\zeta)>b).
\end{equation}
By Hoeffding’s covariance identity (see, e.g., \cite{Lehmann1966}), for any random variables $X$, $Y$ and increasing continuously differentiable functions $f$, $g$ such that $\E f^2(X)$ and $\E g^2(Y)$ are finite,
\begin{equation}
\label{eq_x19}
 \Cov \bigl(f(X),\, g(Y) \bigr)=\int_{-\infty}^{\infty}\int_{-\infty}^{\infty} f'(a) g'(b)
 \bigl[\PP(X>a,\;Y>b) - \PP(X>a)\PP(Y>b)\bigr] \dd a\dd b.
\end{equation}
Applying \eqref{eq_x19} with $X=Y_\eta(\ttt;\varphi)$, $Y=Y_\eta(\ttt;\zeta)$,
$f(x)=x^{2p-1}$, and $g(x)=x^{2q-1}$ for integers $p,q\ge1$, the integrand is non-negative by \eqref{eq_x18} and we conclude that 
\begin{equation}   \label{eq:ghrcor}
\Cov\bigl(Y_\eta(\ttt;\varphi)^{2p-1}, Y_\eta(\ttt;\zeta)^{2q-1}\bigr)\ge 0.
\end{equation}
In particular, since $\E Y_\eta(\ttt;\zeta)=0$, we obtain for  non-negative $\varphi, \zeta$ the inequality
\begin{equation}   \label{eq:ghrcor2}
\E\left[Y_\eta(\ttt;\varphi)^{2p-1} \, Y_\eta(\ttt;\zeta)\right]\ge 0, \qquad p=1,2,\dots.
\end{equation}
Since $\psi-\phi$ and $\psi+\phi$
are nonnegative, \eqref{eq:ghrcor2} with $p=1$ gives \eqref{eq_x20}:
\[
 \E\bigl[Y^2_\eta(\ttt;\psi)\bigr]-\E\bigl[Y_\eta^2(\ttt;\phi)\bigr]=\E\bigl[Y_\eta(\ttt;\psi+\phi)\, Y_\eta(\ttt;\psi-\phi)\bigr]\ge0.
\]
For the higher moments in \eqref{eq_x21}, we first consider the case $\phi\ge 0$. In this situation
\begin{multline}
\label{eq_x23}
 \E\left[Y^{2r}_\eta(\ttt;\psi)\right]-\E\left[Y^{2r}_\eta(\ttt;\phi)\right]=\int_0^1 \frac{\partial}{\partial s}  \E\left[Y^{2r}_\eta(\ttt;\, s\psi+(1-s) \phi)\right]\,\dd s
 \\= 2r \int_0^1\E\Bigl[Y^{2r-1}_\eta(\ttt;\, s\psi+(1-s) \phi)\, Y_\eta(\ttt;\, \psi-\phi) \Bigr]\, \dd s\ge 0,
\end{multline}
where the integrand in the last line is non-negative by \eqref{eq:ghrcor2} with $p=r$, $\varphi=s\psi+(1-s) |\phi|\ge 0$, $\zeta=\psi-\phi\ge 0$. The general $\phi$ case is reduced to the non-negative $\phi$ case at the cost of $2^{2r}$ prefactor. Indeed, let $\phi_\pm=\max\{\pm\phi,0\}\ge 0$ be positive and (minus) negative parts of $\phi$, so that $\phi_+-\phi_-=\phi$ and $\phi_++\phi_-=|\phi|$. Applying $(f-g)^{2r}\le 2^{2r-1}(f^{2r}+g^{2r})$ and twice \eqref{eq_x23} we obtain
$$
  \E\left[Y^{2r}_\eta(\ttt;\phi)\right]\le 2^{2r-1}\E\left[Y^{2r}_\eta(\ttt;\phi_+)+Y^{2r}_\eta(\ttt;\phi_-)\right]
  \le 2^{2r-1}\E\left[Y^{2r}_\eta(\ttt;\psi)+Y^{2r}_\eta(\ttt;\psi)\right]. \qedhere
$$
\end{proof}

\begin{proof}[Proof of Theorem \ref{Theorem_main_detailed} for class I test-measures]
Our plan is to approximate an arbitrary function $\phi$ by linear combinations of exponential functions, apply the results of Section \ref{Section_exponential_test} for them and control the error terms using Lemma \ref{lem:surface-association}.

\noindent{\bf Step 1. Polynomial approximation.} Choose a function
$\phi$ satisfying \eqref{eq:admissible-spatial-test} with
constants $B,\lambda,\Lambda>0$. Let $\alpha=\min(\lambda,\Lambda)/2$ and fix $R>0$. Since $|\phi(x)|\le2B\exp(2\alpha x)$ for $x\le0$, the function
\[
 q\mapsto \begin{cases} q^{-1}\phi(\alpha^{-1}\log q),&
 0<q\le\exp(\alpha R),\\
 0,& q=0.
 \end{cases}
\]
is continuous on $[0,\exp(\alpha R)]$.
For each $\delta>0$, approximating this function
by a polynomial, multiplying by $q$, and changing the variables $q=\exp(\alpha x)$, gives a function $P_\delta(x)$ such that
\[
 P_\delta(x)=\sum_{j=1}^n a_j\exp(j\alpha x)
 \qquad\text{and}\qquad |\phi(x)-P_\delta(x)|\le\delta\exp(\alpha x)\quad \text{ for all } x\le R.
\]
Since both $P_\delta$ and $\phi$ have at most exponential growth as $x\to+\infty$, we can extend the approximation to all $\R$ in the following form: for each $\delta>0$ there exists $P_\delta(x)$ and a constant $C_{\delta}>\alpha$ such that
\begin{equation}
\label{e_x24}
 P_\delta(x)=\sum_{j=1}^n a_j\exp(j\alpha x)
 \qquad\text{and}\qquad |\phi(x)-P_\delta(x)|\le\delta\exp(\alpha x)+\delta\exp(C_{\delta} x) \quad \text{ for all } x\in\R.
\end{equation}

\smallskip

\noindent{\bf Step 2. Behavior of the limit under approximation.} We claim that as $\delta\to 0$, the pairing of ${\overline{\GFF_\U\circ \Omega}}$ with $\nu_{t,P_\delta}$ of Definition \ref{Definition_class_I_measure} approximates the pairing with $\nu_{t,\phi}$. Indeed, recalling the notation $\langle \nu\rangle$ for the pairing of $\nu$ with ${\overline{\GFF_\U\circ \Omega}}$  in Definition \ref{Definition_extended_GFF}, we have:
\begin{multline}
\label{eq_x22}
  \E \left[\bigl(\langle \nu_{t,\phi}\rangle-\langle\nu_{t,P_\delta}\rangle\bigr)^2\right]=\E \left[\langle \nu_{t,\phi-P_\delta}\rangle^2\right]
  \\=\frac{1}{2\pi} \iint_{\R_{<0}^2}\log\left|\frac{\Omega(t,x_1)-\overline \Omega(t,x_2)}{\Omega(t,x_1)-\Omega(t,x_2)}\right| \bigl(\phi(x_1)-P_\delta(x_1)\bigr)\bigl(\phi(x_2)-P_\delta(x_2)\bigr)\dd x_1 \dd x_2.
\end{multline}
The logarithm in the last formula is positive. Hence, using the bound \eqref{e_x24}  and noting that $0<\exp\bigl(C_{\delta}x\bigr)<\exp(\alpha x)$ for $x<0$, we upper-bound \eqref{eq_x22}  by
$$
\frac{2\delta^2}{\pi} \iint_{\R_{<0}^2}\log\left|\frac{\Omega(t,x_1)-\overline \Omega(t,x_2)}{\Omega(t,x_1)-\Omega(t,x_2)}\right| \exp(\alpha x_1) \exp(\alpha x_2)\dd x_1 \dd x_2.
$$
The last expression tends to $0$ as $\delta\to 0$. Since both $\langle \nu_{t,P_\delta}\rangle$ and $\langle \nu_{t,\phi}\rangle$ are centered Gaussian random variables, we conclude that all their moments become arbitrarily close as $\delta\to 0$. 

\smallskip

\noindent{\bf Step 3. Approximation of all moments for finite $\eta$.} 
We fix $k=2,\dots$ and compare the $k$--th moments of $Y_\eta(\ttt;\phi)$ and $Y_\eta(\ttt; P_\delta)$. We introduce the shortcut notations $Y$ and $\tilde Y$ for the former and latter, respectively. Using Hölder's inequality, we have:
\begin{multline}
\label{eq_x24}
 \left|\E Y^k - \E \tilde Y^k\right|= 
 \left|\E \bigl(Y - \tilde Y\bigr)\bigl(Y^{k-1}+Y^{k-2}\tilde Y+\dots + \tilde Y^{k-1}\bigr)\right|\\
 \le \sqrt{\E \bigl(Y - \tilde Y\bigr)^2}
 \sqrt{\E \bigl(Y^{k-1}+Y^{k-2}\tilde Y+\dots + \tilde Y^{k-1}\bigr)^2}\\ \le k  \sqrt{\E \bigl(Y - \tilde Y\bigr)^2}  \sqrt{\E\bigl[ 
 Y^{2k-2}\bigr] + \E\bigl[\tilde Y^{2k-2}\bigr]}
\end{multline} 
We claim that as $\eta\to 0$ and then $\delta\to 0$, the first factor in the last line goes to $0$, while the second one stays bounded. Indeed, using \eqref{eq_x20} and \eqref{e_x24},
\begin{multline}
 \E \bigl(Y - \tilde Y\bigr)^2= \E Y^2_\eta(t; \phi-P_\delta) \le \E Y^2_\eta \bigl(t;\, \delta\exp(\alpha x)+\delta\exp(C_{\delta} x) \bigr)\\ \le 2\delta ^2  \E Y^2_\eta \bigl(t;\,\exp(\alpha x)\bigr) +2\delta^2 \E Y^2_\eta\bigl(t;\,\exp(C_{\delta}x)\bigr).
\end{multline}
We first send $\eta\to 0$, using \eqref{eq_x16}, \eqref{eq_x25}, and \eqref{eq_variance_Laplace}, resulting in
$$
 \limsup_{\eta\to 0+} \E \bigl(Y - \tilde Y\bigr)^2\le \frac{\delta^2}{\pi \beta} \left( \frac{1}{\alpha^2}+\frac{1}{C_\delta^2}\right). 
$$
Since $\alpha$ is fixed and $C_\delta>\alpha$, the last expression tends to $0$ as $\delta\to 0$. 

We turn to the second factor in the last line of \eqref{eq_x24}. For the first term $\E\bigl[ 
 Y^{2k-2}\bigr]$, we produce an estimate using the bound on $\phi$ \eqref{eq:admissible-spatial-test} and comparison inequality \eqref{eq_x21} with $r=k-1$: 
 $$\E\bigl[ 
 Y^{2k-2}\bigr]\le 2^{2k-2}\E Y_\eta^{2k-2}(t,\, B\bigl(\exp(\lambda x)+\exp(\Lambda x)\bigr)).$$
 The right-hand side converges to a Gaussian moment, and therefore stays bounded by Section \ref{Section_exponential_test}.
  
  For the moment of $\tilde Y= Y_\eta(t; P_\delta)$, we note that \emph{a priori} the coefficients $a_j$ in \eqref{e_x24} might have been large. However, we can use the closeness of $P_\delta$ to $\phi$ in the form
 $$
  Y_\eta^{2k-2}(t; P_\delta)= \bigl(Y_\eta(t; \phi)+Y_\eta(t;\, P_\delta-\phi)\bigr)^{2k-2}\le 2^{2k-3} \left[
   Y_\eta^{2k-2}(t; \phi)+Y_\eta^{2k-2}(t;\,P_\delta-\phi)\right].
 $$
 Then using \eqref{e_x24}, \eqref{eq_x21} with $r=k-1$, and Section \ref{Section_exponential_test} we upper-bound
 \begin{multline*}
  \limsup_{\eta\to 0}\E Y_\eta^{2k-2}(t; P_\delta)\le  2^{2k-3} \limsup_{\eta\to 0}\E Y_\eta^{2k-2}(t; \phi) \\
  + 2^{(2k-3)+(2k-3)+(2k-2)} \delta^{2k-2} \E \left[\left(\sqrt{\tfrac{2}{\beta \pi}} \langle \nu_{t,\exp(\alpha x)}\rangle\right)^{2k-2} + \left(\sqrt{\tfrac{2}{\beta \pi}} \langle \nu_{t,\exp(C_\delta x)}\rangle\right)^{2k-2}\right].
 \end{multline*}
 Since $C_\delta>\alpha$ and the pairings are centered Gaussian random variables with variance given by \eqref{eq_cov_I_I}, the last sum remains bounded as $\delta\to 0$. We conclude that \eqref{eq_x24} tends to $0$.

\smallskip

\noindent{\bf Step 4. Finishing the proof.} We now combine the results of the three steps together. For a fixed $\delta>0$ and $t\in\R$, we choose $P_\delta$ from Step 1 and conclude that $\sqrt{\tfrac{\beta \pi}{2}} Y_\eta(t;\, P_\delta)$ converges in moments to $\langle \nu_{t;P_\delta}\rangle$ by Section \ref{Section_exponential_test} applied to a linear combination of exponential functions. By Step 2, for each $k=1,2,\dots$, the $k$th moment of $\langle \nu_{t;P_\delta}\rangle$ approaches the $k$-th moment of $\langle \nu_{t;\phi}\rangle$ as $\delta\to 0$. Simultaneously, by Step 3,
$$
 \lim_{\delta\to 0+} \limsup_{\eta\to 0+} |\E Y_\eta^k(t;\, \phi) - \E  Y_\eta^k(t;\, P_\delta)|=0.
$$
Therefore,
$$
 \lim_{\eta\to 0+} \E \left[\left(\sqrt{\tfrac{\beta \pi}{2}}  Y_\eta(t;\, \phi)\right)^k\right] = \E \bigl[\langle \nu_{t;\phi}\rangle^k\bigr].
$$
This proves Theorem \ref{Theorem_main_detailed} for class I test-measures and $m=1$. For $m>1$ the argument is the same by approximating joint moments instead of individual moments and replacing the identity in the first line of \eqref{eq_x24} with a telescoping sum:
$$
 \prod_{j=1}^k Z_j-\prod_{j=1}^k \tilde Z_j =\sum_{\ell=1}^k \left[ (Z_\ell-\tilde Z_\ell) \prod_{j<\ell} \tilde Z_j \prod_{j>\ell} Z_j\right].\qedhere 
$$
\end{proof}

\subsection{Convergence for class II test-measures}
\label{Section_class_II}

The extension of Theorem \ref{Theorem_main_detailed} from class I to class II signed test-measures follows by expressing the latter as integrals of the former. For a function $f$ as in Definition \ref{Definition_class_II_measure}, we write using the notation \eqref{eq_x26}:
\begin{equation}
   \sqrt{\tfrac{\beta \pi}{2}} \iint_{\R^2} H_\eta(\ttt,x) \nu_f(\dd \ttt,\dd x) = \sqrt{\tfrac{\beta \pi}{2}} \ \int_{\R} Y_\eta\bigl(\ttt;\, f(\ttt,\cdot)\bigr) \dd\ttt. 
\end{equation}
For several functions $f_1,\dots,f_k$, their joint mixed moment is computed as a $k$--fold integral:
\begin{equation}
\label{eq_x27}
   \left(\sqrt{\tfrac{\beta \pi}{2}}\right)^k \E\left[\prod_{j=1}^k \iint_{\R^2} H_\eta(\ttt,x) \nu_{f_j}(\dd \ttt,\dd x)\right] = 
   \left(\sqrt{\tfrac{\beta \pi}{2}}\right)^k \int_{\R^k} \E \left[\prod_{j=1}^k Y_\eta\bigl(\ttt_j;\, f_j(\ttt_j,\cdot)\bigr)\right] \dd\ttt_1\cdots \dd\ttt_k. 
\end{equation}
The same is true for the limiting pairings with GFF in Definition \ref{Definition_extended_GFF}:
\begin{equation}
\label{eq_x28}
  \E\left[\prod_{j=1}^k \langle \nu_{f_j}\rangle\right] = 
  \int_{\R^k} \E \left[\prod_{j=1}^k \langle \nu_{\ttt_j;\, f_j(\ttt_j,\cdot)}\rangle\right] \dd\ttt_1\cdots \dd\ttt_k.  
\end{equation}
The convergence as $\eta\to 0$ of integrals \eqref{eq_x27} to integrals \eqref{eq_x28} follows from the convergence of the integrands established in Section \ref{Section_class_I} by using the dominated convergence theorem. The bound \eqref{eq:admissible-spacetime-test} and comparison \eqref{eq_x21} provide the necessary upper-bound.

More generally, for $k$ class II signed test-measures and $\ell$ class I test-measures, the expectation of the product of integrals of  $H_\eta(\ttt,x)$ is computed by a $k$--fold integral of products of $Y_\eta$ random variables. This leads to the $\eta\to 0+$ asymptotics for all joint moments of the vector \eqref{eq_x17} in Theorem \ref{Theorem_main_detailed}. Since the limiting vector is Gaussian, convergence of all mixed moments also implies convergence in distribution. This finishes the proof of Theorems \ref{thm:surface} and \ref{Theorem_main_detailed}.

\appendix

\section{Convergence of the graph expansions}
\label{app:graph-beta-expansion}

In this appendix, we first fix a graph $\Gamma$, and justify the definitions of the graph-restricted
integrals $\bL^\Gamma_{\beta;\,(\epsilon)}$ in \eqref{eq_L_graph_cutoff} and
$\bL^\Gamma_{[r];\,(\epsilon)}$ in \eqref{eq_L_graph_cutoff_r}.  We then establish their limits as
$\epsilon\to0+$ to justify the definitions of $\bL^\Gamma_\beta$
and $\bL^\Gamma_{[r]}$, and prove 
Lemma \ref{Lemma_Graph_expansion_beta}.

Throughout, $\vec\bk,\vec\ttt$ are fixed as in Section~\ref{sec:basic}.
Constants may depend on the fixed graph and parameters.

\paragraph{The integrals.}
First, we fix $\Gamma$ to be a graph on $\llbracket1,m\rrbracket$, and justify the definitions of
$\bL^\Gamma_{\beta;\,(\epsilon)}$ and
$\bL^\Gamma_{[r];\,(\epsilon)}$ for $\epsilon>0$.

The integrals in Definition~\ref{defn:bIs} are nonnegative and have
finite exponential moments of every positive order.
Then as a function of $\beta^{-1}$, each of $\bI(x;h,g)$, $\bI_0(x;h)$, and $\bI_{0,0}(x)$ extends to a holomorphic function in $\beta^{-1}\in \C$.
Then so is $\beta^{-1}\mapsto \bI_{\beta,\vec\bk}[\vec\bp,\vec\bb, \bH]$ from Definition \ref{def:xibi},
for any $(\vec\bp,
\vec\bb, \bH)\in \sK[\vec\bk]$. Moreover, for every $R>0$,
\begin{equation}
 \sup_{|\beta^{-1}|\leq R}
 \left|\bI_{\beta,\vec\bk}
       [\vec\bp,\vec\bb,\bH]\right|
 \leq
 \left|\bI_{R^{-1},\vec\bk}
       [\vec\bp,\vec\bb,\bH]\right|.
 \label{eq:appg-positive-majorant}
\end{equation}
For $r\geq0$, Cauchy's formula applied to
\eqref{eq_I_r} gives
\begin{equation}\label{eq:appg-kernel-Cauchy}
 \bI_{\vec\bk}^r[\vec\bp,\vec\bb,\bH]
 =
 \frac{1}{2\pi\ii}\oint_{|w|=R}
 \frac{\bI_{w^{-1},\vec\bk}[\vec\bp,\vec\bb,\bH]
       }{w^{r+1}}\dd w,
\end{equation}
where the contour (and all other contours in this appendix) is positively oriented.
Then \eqref{eq:appg-positive-majorant}
implies
\begin{equation}   \label{eq:RIRk}
 \left|\bI_{\vec\bk}^r[\vec\bp,\vec\bb,\bH]\right|
 \leq R^{-r}
 \left|\bI_{R^{-1},\vec\bk}[\vec\bp,\vec\bb,\bH]\right|.
\end{equation}
The fact that the integral of $\bL^\Gamma_{\beta;\,(\epsilon)}$  in \eqref{eq_L_graph_cutoff} is absolutely convergent is a direct consequence of \cite[Proposition~4.42]{GXZ}, by restricting to $\Gamma(\vec\bp)=\Gamma$. Note that by \eqref{eq:appg-positive-majorant}, such absolute convergence holds for any $\beta^{-1}\in \C$, and the finiteness stated in \cite[Proposition~4.42]{GXZ} implies that
\begin{equation}  \label{eq:Lgammabeb}
\sup_{|\beta^{-1}|\le R}\left| \bL^\Gamma_{\beta;\,(\epsilon)} \right|<\infty.
\end{equation}

By \eqref{eq:RIRk} and again using \cite[Proposition~4.42]{GXZ}, we also have that the integral of $\bL^\Gamma_{[r];\,(\epsilon)}$ in \eqref{eq_L_graph_cutoff_r} is absolutely convergent.
Moreover, Fubini's theorem permits interchange of the block
integral and the contour integral, and yields
\begin{equation}\label{eq:appg-Cauchy}
 \bL^\Gamma_{[r];\,(\epsilon)}
 =
 \frac{1}{2\pi\ii}\oint_{|w|=R}
 \frac{\bL^\Gamma_{w^{-1};\,(\epsilon)}}{w^{r+1}}\dd w.
\end{equation}

\paragraph{The limits as $\epsilon\to0+$.}
Second, we continue working with a fixed $\Gamma$, and justify the definitions of $\bL^\Gamma_\beta$
and $\bL^\Gamma_{[r]}$.
The estimate needed is the following:
there is $\epsilon_0>0$ such that, for every $R>0$,
\begin{equation}
 \sup_{|\beta^{-1}|\leq R}
 \left|\bL^\Gamma_{\beta;\,(\epsilon)}
       -\bL^\Gamma_{\beta;\,(\epsilon')}\right|
 \leq C(\sqrt\epsilon+\sqrt{\epsilon'}),
 \qquad 0<\epsilon,\epsilon'\leq\epsilon_0.
 \label{eq:appg-cutoff-Cauchy}
\end{equation}
Here and in the rest of this appendix, we use $C>0$ to denote a large constant which may depend on $R$ and $\vec\bk,\vec\ttt$, and its value can change from line to line.

This bound follows from the proof of
\cite[Proposition~4.43]{GXZ}, and here we outline the key steps.

We work with the discrete walks, as done in \cite[Section 4]{GXZ}.
For any walk $\vec r$ as in \cite[Definition~4.7]{GXZ}, define
a graph $\Gamma(\vec r)$ on $\llbracket1,m\rrbracket$ as follows,
using $i_\ell,Q_\ell,\Delta_{j,\ell},J,U$ from
\cite[Section~4.4]{GXZ}.
For each $j\in J\cup U$, let $a<b$ range over consecutive
elements of the ordered set
$\{0,m+1\}\cup\{\ell\in\llbracket1,m\rrbracket:i_\ell=j\}$.
List the jumps in $\bigcup_{a<\ell<b}\Delta_{j,\ell}$ in time
order, assigning label $\ell$ to each jump in $\Delta_{j,\ell}$.
If $a>0$ and $r_j(Q_a)>0$, prepend the label $a$; if $b\leq m$
and $r_j(Q_{b-1})>0$, append the label $b$.
Each pair of consecutive entries in each resulting list contributes
an undirected edge between its labels; all multiplicities and loops
are retained.  In particular, when $b=a+1$, the list contributes
the edge $(a, b)$ exactly when $r_j(Q_a)>0$, even though there
are no intervening jumps.
For the regularized walks, this is the graph of the discrete
blocks in \cite[Lemma~4.35]{GXZ}, with the same endpoint convention
as Definition~\ref{Definition_connectivity_graph}.

Then one can prove a variant of \cite[Proposition~4.22]{GXZ}, where (following the notations there)
\begin{equation}   \label{eq:rBsum}
\left| \sum_{\vec r \in \mathscr{B}\setminus \mathscr{B}_\epsilon, \Gamma(\vec r)=\Gamma} w(\vec r) \right| <  C\sqrt{\epsilon} \prod_{\ell=1}^m (2\sqrt{NN_\ell})^{k_\ell} N^{-|J|},
\end{equation}
for any $\hat{C}N^{-2/3}<\epsilon<\hat{C}$, and $|\beta^{-1}|\le R$.
Here $\hat{C}>0$ is a universal constant.
The proof follows that of \cite[Proposition~4.22]{GXZ}.
Indeed, the walks are grouped according to their images under the map
$\mathcal{Q}$ from \cite[Definition~4.19]{GXZ}, and cancellation is used
only within each group, as in \cite[Lemma~4.23]{GXZ}.
All remaining estimates are obtained by taking absolute values.
Since $\Gamma(\mathcal{Q}(\vec r))=\Gamma(\vec r)$ for every
$\vec r\in\mathscr{B}$, restricting to a fixed graph $\Gamma$ preserves
these groups, so the same cancellation estimate applies.
To allow complex $\beta$, we replace the factor $|w(\vec r)|$ in the
bound of \cite[Lemma~4.23]{GXZ} by its value with $\beta$ replaced by
$|\beta|$.
The remaining estimates then yield \eqref{eq:rBsum}.

One can further prove a variant of \cite[Proposition~4.25]{GXZ}, where (following the notations there)
\begin{equation}   \label{eq:rBsumep}
\left|
N^m\sum_{\vec r \in \tilde{\mathscr{B}}_\epsilon^*,
\Gamma(\vec r)=\Gamma} w(\vec r) \right|
< C\prod_{\ell=1}^m
(2\sqrt{NN_\ell})^{k_\ell} (\sqrt{\epsilon}+N^{-1/12}),
\end{equation}
for $\epsilon$ and $\beta$ as above.
Here $\tilde{\mathscr{B}}_\epsilon^*$ is defined the same way as $\mathscr{B}_\epsilon^*$ in \cite[Proposition~4.25]{GXZ}, with the condition $\min\{\vartheta_\ell,\dot\vartheta_\ell\}>\epsilon N^{2/3}$ replaced by $\min\{\vartheta_\ell,\dot\vartheta_\ell\}\le\epsilon N^{2/3}$ for at least one $\ell \in J_{\mathrm I}\cup J_{\mathrm{II}}\cup J_{\mathrm{IV}}$.
This follows from \eqref{eq:rBsum}, plus a variant of \cite[Lemma 4.29]{GXZ}, which is restricted to fixed $\Gamma(\vec r)$ and is extended to all $|\beta^{-1}|\le R$.
One notes that the map $\mathcal{M}$ from \cite[Definition 4.26]{GXZ} also leaves the graph invariant: we always have $\Gamma(\mathcal{M}(\vec r))=\Gamma(\vec r)$ for any $\vec r \in \mathscr{B}$.

By passing $N\to \infty$ in \eqref{eq:rBsumep}, using arguments as in the proof of \cite[Proposition 4.44]{GXZ} (which can be adapted to fixed $\Gamma(\vec r)$ and all $|\beta^{-1}|\le R$), we get \eqref{eq:appg-cutoff-Cauchy}.
This justifies the definition of $\bL^\Gamma_\beta=\lim_{\epsilon\to0+}\bL^\Gamma_{\beta;\,(\epsilon)}$, uniformly for $|\beta^{-1}|\le R$.
By further using \eqref{eq:appg-Cauchy}, for any $r\ge 0$,
this justifies the definition of $\bL^\Gamma_{[r]}=\lim_{\epsilon\to0+}\bL^\Gamma_{[r];\,(\epsilon)}$, with
\begin{equation}\label{eq:appg-limit-Cauchy}
 \bL^\Gamma_{[r]}
 =
 \frac{1}{2\pi\ii}\oint_{|w|=R}
 \frac{\bL^\Gamma_{w^{-1}}}{w^{r+1}}\dd w.
\end{equation}

Finally, the contour identities above give the expansions
and their absolute convergence.

\begin{proof}[Proof of Lemma~\ref{Lemma_Graph_expansion_beta}]
We first prove the absolute convergence of \eqref{eq_beta_expansion}.
From the Cauchy integrals \eqref{eq:appg-Cauchy} and \eqref{eq:appg-limit-Cauchy}, to justify the expansions and show that the series are absolutely convergent,
it suffices to bound $|\bL^\Gamma_{[r];\,(\epsilon)}|$ uniformly in $0<\epsilon\le \epsilon_0$, and $|\bL^\Gamma_{[r]}|$.

Choose any $R>0$. By \eqref{eq:Lgammabeb} and \eqref{eq:appg-cutoff-Cauchy}, we have
\[
\sup_{|\beta^{-1}|\le R}\left| \bL^\Gamma_{\beta;\,(\epsilon)} \right| \le C(\sqrt{\epsilon} + \sqrt{\epsilon_0}) + \sup_{|\beta^{-1}|\le R}\left| \bL^\Gamma_{\beta;\,(\epsilon_0)} \right| < C.
\]
Thus by \eqref{eq:appg-Cauchy}, for $r\ge 0$ we have
\[
|\bL^\Gamma_{[r];\,(\epsilon)}| < CR^{-r}.
\]
Note that this is uniform in $0<\epsilon\le \epsilon_0$, by sending $\epsilon\to0+$ we get
\[
|\bL^\Gamma_{[r]}| < CR^{-r}.
\]
These imply the expansions \eqref{eq_beta_expansion}, and that the series are absolutely convergent, for any $|\beta^{-1}|<R$.
Since $R>0$ is arbitrary, we have the expansions with absolute convergence for any $\beta>0$.

We next prove the absolute convergence of \eqref{eq_beta_expansion_total}.
All the arguments so far with a fixed graph $\Gamma$ can be replaced by the summation over a set of graphs.
Namely, taking $\mathcal{C}$ to be a set of graphs on $\llbracket 1, m\rrbracket$, the same arguments so far give, for any $R>0$,
\[
\left| \sum_{\Gamma\in\mathcal{C}} \bL^\Gamma_{[r];\,(\epsilon)} \right| < CR^{-r}, \qquad
\left| \sum_{\Gamma\in\mathcal{C}} \bL^\Gamma_{[r]} \right| < CR^{-r}.
\]
Note that the constant $C$ may depend on $R$, but is independent of the set $\mathcal{C}$. Thus we get
\[
\sum_{\Gamma\text{ is a graph on }\llbracket1,m\rrbracket}|\bL^\Gamma_{[r];\,(\epsilon)}| < CR^{-r}, \qquad \sum_{\Gamma\text{ is a graph on }\llbracket1,m\rrbracket}|\bL^\Gamma_{[r]}| < CR^{-r}.
\]
Note that this gives the absolute convergence of the series \eqref{eq_beta_expansion_total}.

To get the equality to $\bL_{\beta}(\vec\bk, \vec\ttt)$ and $\bL_{\beta;\, (\epsilon)}(\vec\bk, \vec\ttt)$, given \eqref{eq_beta_expansion}, we still need to justify \eqref{eq_Graph_expansion}.
However, the absolute convergence of the series \eqref{eq_beta_expansion_total} already implies that the series in \eqref{eq_Graph_expansion} are absolute convergent.
Then from definitions we get that the second equality in \eqref{eq_Graph_expansion} holds.
To send $\epsilon\to0+$ on both sides, we note that, the same arguments for \eqref{eq:appg-cutoff-Cauchy} would prove that,
\[
 \sup_{|\beta^{-1}|\leq R}
 \left|\sum_{\Gamma\in\mathcal{C}} \bL^\Gamma_{\beta;\,(\epsilon)}
       -\bL^\Gamma_{\beta;\,(\epsilon')}\right|
 \leq C(\sqrt\epsilon+\sqrt{\epsilon'}),
 \qquad 0<\epsilon,\epsilon'\leq\epsilon_0,
\]
for any $R>0$ and  $\mathcal{C}$ being a set of graphs on $\llbracket 1, m\rrbracket$, with $C$ depending on $R$ but independent of $\mathcal{C}$.
Thus we get
\[
 \sup_{|\beta^{-1}|\leq R}
 \sum_{\Gamma\text{ is a graph on }\llbracket1,m\rrbracket} | \bL^\Gamma_{\beta;\,(\epsilon)}
       -\bL^\Gamma_{\beta;\,(\epsilon')}|
 \leq C(\sqrt\epsilon+\sqrt{\epsilon'}),
 \qquad 0<\epsilon,\epsilon'\leq\epsilon_0.
\]
Then we can pass $\epsilon\to0+$ in both sides of the second equality in \eqref{eq_Graph_expansion}, and get the first equality in \eqref{eq_Graph_expansion}.
\end{proof}

\bibliographystyle{alpha}
\bibliography{references}

@article{AH,
  author = {Aggarwal, Amol and Huang, Jiaoyang},
  title = {Edge statistics for lozenge tilings of polygons, {II}: {Airy} line ensemble},
  journal = {Forum of Mathematics, Pi},
  volume = {13},
  pages = {e2},
  year = {2025},
  doi = {10.1017/fmp.2024.16},
  eprint = {2108.12874},
  archivePrefix = {arXiv},
  url = {https://arxiv.org/abs/2108.12874}
}

@article{AHCharacterization,
  author        = {Aggarwal, Amol and Huang, Jiaoyang},
  title         = {Strong characterization for the {Airy} line ensemble},
  journal       = {Invent. Math.},
  volume        = {243},
  number        = {2},
  pages         = {293--605},
  year          = {2026},
  doi           = {10.1007/s00222-025-01381-6},
  eprint        = {2308.11908},
  archivePrefix = {arXiv},
  primaryClass  = {math.PR}
}

@article{BG,
  author = {Borodin, Alexei and Gorin, Vadim},
  title = {General {$\beta$}-{Jacobi} corners process and the {Gaussian} free field},
  journal = {Communications on Pure and Applied Mathematics},
  volume = {68},
  number = {10},
  pages = {1774--1844},
  year = {2015},
  doi = {10.1002/cpa.21546},
  eprint = {1305.3627},
  archivePrefix = {arXiv},
  url = {https://arxiv.org/abs/1305.3627}
}

@article{Borodin,
  author = {Borodin, Alexei},
  title = {{CLT} for spectra of submatrices of {Wigner} random matrices},
  journal = {Moscow Mathematical Journal},
  volume = {14},
  number = {1},
  pages = {29--38},
  year = {2014},
  doi = {10.17323/1609-4514-2014-14-1-29-38},
  eprint = {1010.0898},
  archivePrefix = {arXiv},
  url = {https://arxiv.org/abs/1010.0898}
}

@incollection{BorodinDynamic,
  author = {Borodin, Alexei},
  title = {{CLT} for spectra of submatrices of {Wigner} random matrices, {II}: Stochastic evolution},
  booktitle = {Random Matrix Theory, Interacting Particle Systems and Integrable Systems},
  editor = {Deift, Percy and Forrester, Peter},
  series = {Mathematical Sciences Research Institute Publications},
  volume = {65},
  publisher = {Cambridge University Press},
  pages = {57--69},
  year = {2014},
  eprint = {1011.3544},
  archivePrefix = {arXiv},
  url = {https://arxiv.org/abs/1011.3544}
}

@article{CH,
  author = {Corwin, Ivan and Hammond, Alan},
  title = {{Brownian Gibbs} property for {Airy} line ensembles},
  journal = {Inventiones Mathematicae},
  volume = {195},
  number = {2},
  pages = {441--508},
  year = {2014},
  doi = {10.1007/s00222-013-0462-3},
  eprint = {1108.2291},
  archivePrefix = {arXiv},
  url = {https://arxiv.org/abs/1108.2291}
}

@misc{DFZ,
  author = {Dimitrov, Evgeni and Fu, Alex and Zhou, Zhengye},
  title = {From the {Airy} line ensemble to the {Gaussian} free field},
  year = {2026},
  eprint = {2607.10725},
  archivePrefix = {arXiv},
  note = {\href{https://arxiv.org/abs/2607.10725}{arXiv:2607.10725}},
  url = {https://arxiv.org/abs/2607.10725}
}

@article{GS,
  author = {Gorin, Vadim and Shkolnikov, Mykhaylo},
  title = {Stochastic {Airy} semigroup through tridiagonal matrices},
  journal = {The Annals of Probability},
  volume = {46},
  number = {4},
  pages = {2287--2344},
  year = {2018},
  eprint = {1601.06800},
  archivePrefix = {arXiv},
  url = {https://arxiv.org/abs/1601.06800}
}

@misc{GXZ,
  author = {Gorin, Vadim and Xu, Jiaming and Zhang, Lingfu},
  title = {{Airy}$_\beta$ line ensemble and its {Laplace} transform},
  year = {2024},
  eprint = {2411.10829},
  archivePrefix = {arXiv},
  note = {\href{https://arxiv.org/abs/2411.10829}{arXiv:2411.10829}},
  url = {https://arxiv.org/abs/2411.10829}
}

@book{GFB,
  author = {Gorin, Vadim and Ferrari, Patrik and Borodin, Alexei},
  title = {Introduction to integrable probability},
  year = {2026},
  note = {\url{https://www.stat.berkeley.edu/~vadicgor/Book_IIP.pdf}}
}

@article{HZ,
  author = {Huang, Jiaoyang and Zhang, Lingfu},
  title = {A convergence framework for {Airy}$_\beta$ line ensemble via pole evolution},
  journal = {Advances in Mathematics},
  volume = {498},
  pages = {111028},
  year = {2026},
  doi = {10.1016/j.aim.2026.111028},
  eprint = {2411.10586},
  archivePrefix = {arXiv},
  url = {https://arxiv.org/abs/2411.10586}
}

@article{Janson,
  author = {Janson, Svante},
  title = {{Brownian} excursion area, {Wright}'s constants in graph enumeration, and other {Brownian} areas},
  journal = {Probability Surveys},
  volume = {4},
  pages = {80--145},
  year = {2007},
  doi = {10.1214/07-PS104},
  eprint = {0704.2289},
  archivePrefix = {arXiv},
  url = {https://arxiv.org/abs/0704.2289}
}

@article{Kenyon,
  author = {Kenyon, Richard},
  title = {Dominos and the {Gaussian} free field},
  journal = {The Annals of Probability},
  volume = {29},
  number = {3},
  pages = {1128--1137},
  year = {2001},
  doi = {10.1214/aop/1015345599},
  eprint = {math-ph/0002027},
  archivePrefix = {arXiv},
  url = {https://arxiv.org/abs/math-ph/0002027}
}

@article{KenyonHoneycomb,
  author = {Kenyon, Richard},
  title = {Height fluctuations in the honeycomb dimer model},
  journal = {Communications in Mathematical Physics},
  volume = {281},
  number = {3},
  pages = {675--709},
  year = {2008},
  doi = {10.1007/s00220-008-0511-8},
  eprint = {math-ph/0405052},
  archivePrefix = {arXiv},
  url = {https://arxiv.org/abs/math-ph/0405052}
}

@article{KRV,
  author        = {Krishnapur, Manjunath and Rider, Brian
                   and Vir{\'a}g, B{\'a}lint},
  title         = {Universality of the stochastic {Airy} operator},
  journal       = {Comm. Pure Appl. Math.},
  volume        = {69},
  number        = {1},
  pages         = {145--199},
  year          = {2016},
  doi           = {10.1002/cpa.21573},
  eprint        = {1306.4832},
  archivePrefix = {arXiv},
  primaryClass  = {math.PR}
}

@misc{Landon,
  author = {Landon, Benjamin},
  title = {Edge scaling limit of {Dyson Brownian} motion at equilibrium for general {$\beta\ge1$}},
  year = {2020},
  eprint = {2009.11176},
  archivePrefix = {arXiv},
  note = {\href{https://arxiv.org/abs/2009.11176}{arXiv:2009.11176}},
  url = {https://arxiv.org/abs/2009.11176}
}

@article{Petrov,
  author = {Petrov, Leonid},
  title = {Asymptotics of uniformly random lozenge tilings of polygons. {Gaussian} free field},
  journal = {The Annals of Probability},
  volume = {43},
  number = {1},
  pages = {1--43},
  year = {2015},
  doi = {10.1214/12-AOP823},
  eprint = {1206.5123},
  archivePrefix = {arXiv},
  url = {https://arxiv.org/abs/1206.5123}
}

@article{PS,
  author = {Pr{\"a}hofer, Michael and Spohn, Herbert},
  title = {Scale invariance of the {PNG} droplet and the {Airy} process},
  journal = {Journal of Statistical Physics},
  volume = {108},
  number = {5--6},
  pages = {1071--1106},
  year = {2002},
  doi = {10.1023/A:1019791415147},
  eprint = {math/0105240},
  archivePrefix = {arXiv},
  url = {https://arxiv.org/abs/math/0105240}
}

@book{PT,
  author = {Peccati, Giovanni and Taqqu, Murad S.},
  title = {{Wiener} Chaos: Moments, Cumulants and Diagrams: A Survey with Computer Implementation},
  series = {Bocconi \& Springer Series},
  volume = {1},
  publisher = {Springer},
  address = {Milan},
  year = {2011},
  doi = {10.1007/978-88-470-1679-8},
  url = {https://doi.org/10.1007/978-88-470-1679-8}
}

@article{RRV,
  author = {Ram{\'i}rez, Jos{\'e} A. and Rider, Brian and Vir{\'a}g, B{\'a}lint},
  title = {Beta ensembles, stochastic {Airy} spectrum, and a diffusion},
  journal = {Journal of the American Mathematical Society},
  volume = {24},
  number = {4},
  pages = {919--944},
  year = {2011},
  doi = {10.1090/S0894-0347-2011-00703-0},
  eprint = {math/0607331},
  archivePrefix = {arXiv},
  url = {https://arxiv.org/abs/math/0607331}
}

@article{Sheffield,
  author = {Sheffield, Scott},
  title = {{Gaussian} free fields for mathematicians},
  journal = {Probability Theory and Related Fields},
  volume = {139},
  number = {3--4},
  pages = {521--541},
  year = {2007},
  doi = {10.1007/s00440-006-0050-1},
  eprint = {math/0312099},
  archivePrefix = {arXiv},
  url = {https://arxiv.org/abs/math/0312099}
}

@book{WernerPowell2022GFF,
  author    = {Werner, Wendelin and Powell, Ellen},
  title     = {Lecture Notes on the {G}aussian Free Field},
  series    = {Cours Sp\'ecialis\'es},
  volume    = {28},
  publisher = {Soci\'et\'e Math\'ematique de France},
  address   = {Paris},
  year      = {2022},
  isbn      = {978-2-85629-952-4},
  issn      = {1284-6090},
}

@article{Sodin,
  author = {Sodin, Sasha},
  title = {A limit theorem at the spectral edge for corners of time-dependent {Wigner} matrices},
  journal = {International Mathematics Research Notices},
  volume = {2015},
  number = {17},
  pages = {7575--7607},
  year = {2015},
  doi = {10.1093/imrn/rnu180},
  eprint = {1312.1007},
  archivePrefix = {arXiv},
  url = {https://arxiv.org/abs/1312.1007}
}

@misc{BorotGG,
  author = {Borot, Ga{\"e}tan and Gorin, Vadim and Guionnet, Alice},
  title = {Macroscopic asymptotics in discrete beta-ensembles and random tilings},
  year = {2026},
  eprint = {2601.16377},
  archivePrefix = {arXiv},
  note = {\href{https://arxiv.org/abs/2601.16377}{arXiv:2601.16377}},
  url = {https://arxiv.org/abs/2601.16377}
}

@misc{HuangHeight,
  author = {Huang, Jiaoyang},
  title = {Height fluctuation for lozenge tilings of polygons},
  year = {2026},
  eprint = {2609.04789},
  archivePrefix = {arXiv},
  note = {\href{https://arxiv.org/abs/2609.04789}{arXiv:2609.04789}},
  url = {https://arxiv.org/abs/2609.04789}
}

@misc{BKM,
  author = {Baslingker, Jnaneshwar and Krishnapur, Manjunath and Madiman, Mokshay},
  title = {Log-concavity in one-dimensional {Coulomb} gases and related ensembles},
  year = {2024},
  eprint = {2412.15116},
  archivePrefix = {arXiv},
  note = {\href{https://arxiv.org/abs/2412.15116v2}{arXiv:2412.15116v2}},
  url = {https://arxiv.org/abs/2412.15116v2}
}

@article{FKG,
  author = {Fortuin, C. M. and Kasteleyn, P. W. and Ginibre, J.},
  title = {Correlation inequalities on some partially ordered sets},
  journal = {Communications in Mathematical Physics},
  volume = {22},
  number = {2},
  pages = {89--103},
  year = {1971},
  doi = {10.1007/BF01651330},
  url = {https://doi.org/10.1007/BF01651330}
}

@article{ES,
  author = {Edelman, Alan and Sutton, Brian D.},
  title = {From random matrices to stochastic operators},
  journal = {Journal of Statistical Physics},
  volume = {127},
  number = {6},
  pages = {1121--1165},
  year = {2007},
  doi = {10.1007/s10955-006-9226-4},
  eprint = {math-ph/0607038},
  archivePrefix = {arXiv},
  url = {https://doi.org/10.1007/s10955-006-9226-4}
}

@misc{BGS,
  author = {Bykhovskaya, Anna and Gorin, Vadim and Sodin, Sasha},
  title = {How weak are weak factors? {U}niform inference for signal strength in signal plus noise models},
  year = {2026},
  eprint = {2507.18554},
  archivePrefix = {arXiv},
  note = {\href{https://arxiv.org/abs/2507.18554v2}{arXiv:2507.18554v2}},
  url = {https://arxiv.org/abs/2507.18554v2}
}

@article{CCAiry,
  author = {Charlier, Christophe and Claeys, Tom},
  title = {Large gap asymptotics for {Airy} kernel determinants with discontinuities},
  journal = {Communications in Mathematical Physics},
  volume = {375},
  number = {2},
  pages = {1299--1339},
  year = {2020},
  doi = {10.1007/s00220-019-03538-w},
  eprint = {1812.01964},
  archivePrefix = {arXiv},
  url = {https://doi.org/10.1007/s00220-019-03538-w}
}

@incollection{CL20,
  author    = {Chafa{\"i}, Djalil and Lehec, Joseph},
  title     = {On {Poincar\'e} and logarithmic {Sobolev} inequalities
               for a class of singular {Gibbs} measures},
  booktitle = {Geometric Aspects of Functional Analysis},
  editor    = {Klartag, Bo'az and Milman, Emanuel},
  series    = {Lecture Notes in Mathematics},
  volume    = {2256},
  pages     = {219--246},
  publisher = {Springer},
  address   = {Cham},
  year      = {2020},
  doi       = {10.1007/978-3-030-36020-7_10},
  eprint    = {1805.00708},
  archivePrefix = {arXiv},
  primaryClass  = {math.PR}
}

@article{BD14,
  author    = {Breuer, Jonathan and Duits, Maurice},
  title     = {The {Nevai} condition and a local law of large numbers
               for orthogonal polynomial ensembles},
  journal   = {Advances in Mathematics},
  volume    = {265},
  pages     = {441--484},
  year      = {2014},
  doi       = {10.1016/j.aim.2014.07.026},
  eprint    = {1301.2061},
  archivePrefix = {arXiv},
  primaryClass  = {math.PR}
}

@article{Lehmann1966,
  author  = {Lehmann, E. L.},
  title   = {Some Concepts of Dependence},
  journal = {The Annals of Mathematical Statistics},
  year    = {1966},
  volume  = {37},
  number  = {5},
  pages   = {1137--1153},
  doi     = {10.1214/aoms/1177699260}
}

@book{GorinBook,
  author = {Gorin, Vadim},
  title = {Lectures on Random Lozenge Tilings},
  series = {Cambridge Studies in Advanced Mathematics},
  volume = {193},
  publisher = {Cambridge University Press},
  address = {Cambridge},
  year = {2021},
  doi = {10.1017/9781108921183},
  isbn = {9781108843966},
  url = {https://doi.org/10.1017/9781108921183}
}

@article{BEY,
  author = {Bourgade, Paul and Erd{\H o}s, L{\'a}szl{\'o} and Yau, Horng-Tzer},
  title = {Edge universality of beta ensembles},
  journal = {Communications in Mathematical Physics},
  volume = {332},
  number = {1},
  pages = {261--353},
  year = {2014},
  doi = {10.1007/s00220-014-2120-z},
  eprint = {1306.5728},
  archivePrefix = {arXiv},
  url = {https://arxiv.org/abs/1306.5728}
}

@article{BFG,
  author = {Bekerman, Florent and Figalli, Alessio and Guionnet, Alice},
  title = {Transport maps for {$\beta$}-matrix models and universality},
  journal = {Communications in Mathematical Physics},
  volume = {338},
  number = {2},
  pages = {589--619},
  year = {2015},
  doi = {10.1007/s00220-015-2384-y},
  eprint = {1311.2315},
  archivePrefix = {arXiv},
  url = {https://arxiv.org/abs/1311.2315}
}

@article{CorwinKPZ,
  author = {Corwin, Ivan},
  title = {The {Kardar--Parisi--Zhang} equation and universality class},
  journal = {Random Matrices: Theory and Applications},
  volume = {1},
  number = {1},
  pages = {1130001},
  year = {2012},
  doi = {10.1142/S2010326311300014},
  eprint = {1106.1596},
  archivePrefix = {arXiv},
  url = {https://arxiv.org/abs/1106.1596}
}

@misc{XuJack,
  author = {Xu, Jiaming},
  title = {{Airy} limit for the {Jack} process and topological expansion},
  year = {2026},
  eprint = {2609.03532},
  archivePrefix = {arXiv},
  note = {\href{https://arxiv.org/abs/2609.03532}{arXiv:2609.03532}},
  url = {https://arxiv.org/abs/2609.03532}
}

@article{BPR,
  author = {Berestycki, Nathana{\"e}l and Powell, Ellen and Ray, Gourab},
  title = {A characterisation of the {Gaussian} free field},
  journal = {Probability Theory and Related Fields},
  volume = {176},
  number = {3--4},
  pages = {1259--1301},
  year = {2020},
  doi = {10.1007/s00440-019-00939-9},
  eprint = {1802.01195},
  archivePrefix = {arXiv},
  url = {https://arxiv.org/abs/1802.01195}
}

@article{CKP,
  author = {Cohn, Henry and Kenyon, Richard and Propp, James},
  title = {A variational principle for domino tilings},
  journal = {Journal of the American Mathematical Society},
  volume = {14},
  number = {2},
  pages = {297--346},
  year = {2001},
  doi = {10.1090/S0894-0347-00-00355-6},
  eprint = {math/0008220},
  archivePrefix = {arXiv},
  url = {https://arxiv.org/abs/math/0008220}
}

@article{BufetovGorin,
  author = {Bufetov, Alexey and Gorin, Vadim},
  title = {Fluctuations of particle systems determined by {Schur} generating functions},
  journal = {Advances in Mathematics},
  volume = {338},
  pages = {702--781},
  year = {2018},
  doi = {10.1016/j.aim.2018.07.009},
  eprint = {1604.01110},
  archivePrefix = {arXiv},
  url = {https://arxiv.org/abs/1604.01110}
}

@article{JohanssonAztec,
  author = {Johansson, Kurt},
  title = {The arctic circle boundary and the {Airy} process},
  journal = {The Annals of Probability},
  volume = {33},
  number = {1},
  pages = {1--30},
  year = {2005},
  doi = {10.1214/009117904000000937},
  eprint = {math/0306216},
  archivePrefix = {arXiv},
  url = {https://arxiv.org/abs/math/0306216}
}

@article{PetrovEdge,
  author = {Petrov, Leonid},
  title = {Asymptotics of random lozenge tilings via {Gelfand--Tsetlin} schemes},
  journal = {Probability Theory and Related Fields},
  volume = {160},
  number = {3--4},
  pages = {429--487},
  year = {2014},
  doi = {10.1007/s00440-013-0532-x},
  eprint = {1202.3901},
  archivePrefix = {arXiv},
  url = {https://arxiv.org/abs/1202.3901}
}

@article{Spohn1998,
  author = {Spohn, Herbert},
  title = {{Dyson}'s model of interacting {Brownian} motions at arbitrary coupling strength},
  journal = {Markov Processes and Related Fields},
  volume = {4},
  number = {4},
  pages = {649--661},
  year = {1998},
  url = {https://www.math-mprf.org/journal/articles/id832/}
}

\end{document}